\documentclass[a4paper, 10pt]{amsart}
\usepackage{graphicx}
\usepackage{enumerate}
\usepackage{xspace}
\usepackage{boxedminipage}
\usepackage[ruled]{algorithm}
\usepackage{subfigure}
\usepackage{booktabs}
\usepackage{tabularx}
\usepackage[square,sort,numbers]{natbib}
\usepackage{hyperref} 
\usepackage{times}
\usepackage{amsmath,amssymb,amsxtra,amsthm}
\usepackage{subfigure}
\usepackage{geometry}                

\newcommand{\Proofname}{Proof}
\newenvironment{Proof}[1][{. }]
	       {\noindent{\bf \Proofname\ #1}}
	       {{\raggedright{{ }\hfill\qed}}} 
    
\newtheoremstyle
    {note}
    {}
    {}
    {}
    {}
    {\bfseries}
    {.}
    {1.0ex}
    {}
    
\swapnumbers
\numberwithin{equation}{section}
\newtheorem{The}[subsection]{Theorem}

\newtheorem{Lem}[subsection]{Lemma}
\newtheorem{Cor}[subsection]{Corollary}
\newtheorem{Prop}[subsection]{Proposition}

\theoremstyle{note}
\newtheorem{Example}[subsection]{Example}
\newtheorem{Rem}[subsection]{Remark}

\newtheorem{Defn}[subsection]{Definition}

\providecommand{\highlight}[1]{{\color{red}#1}}
\usepackage{ifthen}
         \provideboolean{showchanges}
         \setboolean{showchanges}{true}
         \provideboolean{newchanges}
         \setboolean{newchanges}{true}
         \providecommand{\changes}[1]{
           \ifthenelse{\boolean{showchanges}}{{\highlight{#1}}}{#1}
         }
 \providecommand{\newchanges}[1]{
           \ifthenelse{\boolean{newchanges}}{{\highlight{#1}}}{#1}
         }
         \providecommand{\changefromto}[3][replace with]{
           \ifthenelse{\boolean{showchanges}}
           {{\sout{#2}\margnote{#1}}{\highlight{#3}}}
           {#3\xspace}
         }
         \providecommand{\ChangePar}[2]{
           \ifthenelse{\boolean{showchanges}}
           {{\par$\mapsfrom$ \textcolor{red!20}{#1}}{\par$\mapsto$ \textcolor{blue}{#2}}}
           {\par #2}
         }
         \providecommand{\InsertPar}[1]{
           \ifthenelse{\boolean{showchanges}}
           {{\par$\mapsto$ \textcolor{blue}{#1}}}
           {\par #1}
         }

\provideboolean{showdelete}
\setboolean{showdelete}{false}
\newcommand{\delete}[1]{
  \ifthenelse{\boolean{showdelete}} {{\color{red}{#1}}}{}
}

\renewcommand{\vec}[1]{\ensuremath{\boldsymbol{#1}}}

\newcommand{\myall}{\ensuremath{\quad \forall}}
\newcommand{\pdt}{\ensuremath{\partial_t}}
\newcommand{\pdtau}{\ensuremath{\partial_\tau}}

\newcommand{\mdt}[1]{\ensuremath{\partial^{\bullet}_{#1}}}
\newcommand{\ntd}{\ensuremath{\partial^{\circ}}}
\newcommand{\mdth}[1]{\ensuremath{\partial^{\bullet}_{h,#1}}}
\newcommand{\mdthn}[1]{\ensuremath{\partial^{\bullet,\tau}_{h,#1}}}

\newcommand{\Reals}{\ensuremath{{\mathbb{R}}}}

\newcommand{\diff}{\ensuremath{{\operatorname{d}}}}

\newcommand{\normal}{\ensuremath{{\vec{\nu}}}}
\newcommand{\conormal}{\ensuremath{{\vec{\mu}}}}
\newcommand{\tangent}{\ensuremath{{\vec{\mathcal{T}}}}}
\newcommand{\surface}{\ensuremath{{\mathcal{M}}}}
\newcommand{\G}{\ensuremath{{\Gamma}}}
\newcommand{\Gt}{\ensuremath{{\G(t)}}}
\newcommand{\Gtn}[1]{\ensuremath{{\G^{#1}}}}
\newcommand{\Gc}{\ensuremath{{\G_h}}}
\newcommand{\Gct}{\ensuremath{{\G_h(t)}}}
\newcommand{\Gctn}[1]{\ensuremath{{\G_h^{#1}}}}
\newcommand{\calGT}{\ensuremath{\mathcal{G}_T}}
\newcommand{\calGTh}{\ensuremath{\mathcal{G}_{h,T}}}
\newcommand{\Lp}[1]{\ensuremath{\operatorname{L}_{#1}}}

\newcommand{\Hil}[1]{\ensuremath{\operatorname{H}^{#1}}}
\newcommand{\Cont}[1]{\ensuremath{\operatorname{C}^{#1}}}
\newcommand{\Sc}{\ensuremath{{\mathcal{S}_h}}}
\newcommand{\Scn}[1]{\ensuremath{{\mathcal{S}^{#1}_h}}}
\newcommand{\Sct}{\ensuremath{{\mathcal{S}_h(t)}}}
\newcommand{\Scl}{\ensuremath{{\mathcal{S}_h^l}}}
\newcommand{\Scln}[1]{\ensuremath{{\mathcal{S}_h^{#1,l}}}}
\newcommand{\Sclt}{\ensuremath{{\mathcal{S}_h^l(t)}}}
\newcommand{\ScT}{\ensuremath{{\mathcal{S}_h^T}}}
\newcommand{\ScTl}{\ensuremath{{(\mathcal{S}_h^T)^l}}}

\newcommand{\ltwon}[2]{\ensuremath{\left\|#1\right\|}_{\Lp{2}\left(#2\right)}}
\newcommand{\Hiln}[3]{\ensuremath{\left\|#1\right\|}_{\Hil{#2}(#3)}}

\newcommand{\lv}{\ensuremath{\left\vert}}
\newcommand{\rv}{\ensuremath{\right\vert}}
\usepackage{upgreek}
\newcommand{\lap}{\ensuremath{\Updelta}}

\newcommand{\T}{\ensuremath{{\mathcal{T}}}}

\newcommand{\Ritz}{\ensuremath{\operatorname{R}^h}}

\newcommand{\abil}[2]{\ensuremath{a\left({#1},{#2}\right)}}
\newcommand{\atbil}[3]{\ensuremath{\tilde{a}\left({#1},{#2};{#3}\right)}}
\newcommand{\mbil}[2]{\ensuremath{m\left({#1},{#2}\right)}}
\newcommand{\gbil}[3]{\ensuremath{g\left({#1},{#2};{#3}\right)}}
\newcommand{\bbil}[3]{\ensuremath{b\left({#1},{#2};{#3}\right)}}
\newcommand{\btbil}[4]{\ensuremath{\tilde{b}\left({#1},{#2};{#3};{#4}\right)}}
\newcommand{\ahbil}[2]{\ensuremath{a_h\left({#1},{#2}\right)}}
\newcommand{\ahtbil}[3]{\ensuremath{\tilde{a}_h\left({#1},{#2};{#3}\right)}}
\newcommand{\mhbil}[2]{\ensuremath{m_h\left({#1},{#2}\right)}}
\newcommand{\ghbil}[3]{\ensuremath{g_h\left({#1},{#2};{#3}\right)}}
\newcommand{\bhbil}[3]{\ensuremath{b_h\left({#1},{#2};{#3}\right)}}
\newcommand{\bhtbil}[4]{\ensuremath{\tilde{b}_h\left({#1},{#2};{#3};{#4}\right)}}
\newcommand{\Ltbil}[2]{\ensuremath{\mathcal{L}_2\left({#1},{#2}\right)}}
\newcommand{\Lthbil}[2]{\ensuremath{\mathcal{L}_{2,h}\left({#1},{#2}\right)}}
\newcommand{\bhn}[2]{\ensuremath{\underline{#1}^{#2}_h}}
\newcommand{\bn}[2]{\ensuremath{\underline{#1}^{#2}}}
\newcommand{\defB}[2]{\ensuremath{\mathcal{B}_{#2}{\left(#1\right)}}}
\newcommand{\defD}[2]{\ensuremath{\mathcal{D}_{#2}{\left(#1\right)}}}

\newcommand{\Proj}{\ensuremath{\vec{P}}}

\usepackage{color}

\definecolor{MyGreen}{rgb} {0.05,0.4,0.05}
\definecolor{RedViolet}{rgb} {0.1,0.1,0.75}

          \providecommand{\highlight}[1]{{\color{blue}#1}}
         \provideboolean{shownotes}
          \setboolean{shownotes}{false}
          \newcommand{\standout}[1]{\colorbox{red}{\textcolor{white}{#1}}}
          \newcounter{margnote}[page]
          \newcommand{\margnotemark}{{\standout{\footnotesize\upshape\texttt{\arabic{margnote}}}}}
          \newcommand{\margnote}[2][]{
            \ifthenelse{
              \boolean{shownotes}
            }{\stepcounter{margnote}\margnotemark\marginpar{
                \texttt{
                  \begin{minipage}{2cm}
                    \raggedright\tiny
                    \margnotemark{#1}: 
                    #2
                  \end{minipage}
            }}}{}
          }
          \provideboolean{showchanges}
          \setboolean{showchanges}{false}
          \providecommand{\chcolor}{\color{blue}}
          \providecommand{\changes}[1]{
            \ifthenelse{\boolean{showchanges}}
          	     {{\chcolor #1}}
          	     {#1
          }
          }
          \providecommand{\changefromto}[3][replace with]{
            \ifthenelse{\boolean{showchanges}}
            {{\sout{#2}\margnote{#1}}{\highlight{#3}}}
            {#3\xspace}
          }
          \providecommand{\ChangePar}[2]{
            \ifthenelse{\boolean{showchanges}}
            {{\par$\mapsfrom$ \textcolor{red!20}{#1}}{\par$\mapsto$ \textcolor{blue}{#2}}}
            {\par #2}
          }
          \providecommand{\InsertPar}[1]{
            \ifthenelse{\boolean{showchanges}}
            {{\par$\mapsto$ \textcolor{blue}{#1}}}
            {\par #1}
         }
	       
\usepackage{hyperref}
\usepackage[british]{babel}
\newcommand \beq{\begin{equation}}
\newcommand \eeq{\end{equation}}

\def\eps{\varepsilon}

\def\bpmat{\begin{pmatrix}}
\def\epmat{\end{pmatrix}}
\def\mathref#1{\ifmmode\mathrm{(\ref{#1})}\else{\rm(\ref{#1})}\fi} 
\def\nref#1{\ifmmode\mathrm{\ref{#1}}\else{\rm\ref{#1}}\fi}

\begin{document}
\bibliographystyle{unsrtnat}
\title{Error analysis for an ALE evolving surface finite element method}
\date{\today}
\author{Charles M. Elliott}
\address[C.M. Elliott]{Mathematics Institute, Zeeman Building, University of Warwick, Coventry, UK, CV4 7AL. }
\email[C.M. Elliott]{C.M.Elliott@warwick.ac.uk}

  \author{Chandrasekhar Venkataraman}
  \address[C. Venkataraman]{Department of Mathematics, University of Sussex, Falmer, UK,  BN1 9QH. }
 \email[C. Venkataraman]{c.venkataraman@sussex.ac.uk}
\begin{abstract}
We consider an arbitrary-Lagrangian-Eulerian  evolving surface finite element method for the numerical approximation 
of advection and diffusion of a conserved scalar quantity on a moving surface. We describe the method, prove optimal order error
bounds and present numerical simulations that agree with the theoretical results.
\end{abstract}
\maketitle
\section{Introduction}\label{sec:intro}
For each $t\in[0,T], T>0,$ let $\Gt$ be a smooth connected hypersurface in $\Reals^{m+1},m=1,2,3$, oriented by the normal vector field $\normal(\cdot,t)$, with $\G^0:=\G(0)$. We assume that there exists a diffeomorphism $\vec G(\cdot,t):\G^0\to\Gt,$ satisfying $\vec G\in\Cont{2}([0,T],\Cont{2}(\G^0))$. 
We set $\vec v(\vec G(\cdot,t),t)=\pdt \vec G(\cdot,t)$ with $\vec G(\cdot,0)=\vec {I}$ (the identity). Furthermore we assume that $\vec v(\cdot,t) \in \Cont{2}(\Gt)$.
The given velocity field $\vec v=\vec v_\normal+\vec v_\tangent$ may contain both normal $\vec v_\normal$ and  tangential $\vec v_\tangent$ components, i.e.,  $\vec v_\normal=\vec v\cdot\normal\normal$ and $\vec v_\tangent\cdot\normal=0$.

 We focus on the following linear parabolic partial differential equation on $\Gt$;
\beq\label{eqn:pde}
\mdt{\vec v} u+u\nabla_{\Gt}\cdot\vec v-\lap_{\Gt}u=0\quad\mbox{ on }\Gt,
\eeq
where, $\nabla_\Gt=\nabla-\nabla\cdot\normal\normal$ denotes the surface gradient, $\lap_\Gt=\nabla_\Gt\cdot\nabla_\Gt$ the Laplace Beltrami operator \changes{
 and    
\margnote{ref 2. pt 2.}
\[
\mdt{\vec v} u=\pdt u+\vec v\cdot\nabla u
=\pdt u+\vec v_\normal\cdot\nabla u+\vec v_\tangent \cdot\nabla_\Gt u
\]
is the material derivative {with respect to the velocity field $\vec v$}.
For simplicity we will assume that the boundary of $\Gt$ is empty and hence no boundary conditions are needed. The method is easily adapted to surfaces with boundary.
In the case that the surface has a boundary, under suitable assumptions (c.f., Remark \ref{rem:bdry}), our analysis is valid for  homogeneous Neumann boundary conditions, i.e.,
\margnote{ ref 2 pt 3.}
\beq\label{eqn:BCs}
\nabla_\Gt u\cdot\conormal=0\quad\text{ on }\partial\Gt,
\eeq
where $\conormal$ is the conormal to the boundary of the surface. The upshot is that the total mass is conserved i.e., $\frac{\diff }{\diff t}\int_\Gt u=0$. Note that the case $\normal(\vec x,t)$ being constant in space and time corresponds to the $n$-dimensional hypersurface $\Gt$ being flat. In this case \eqref{eqn:pde} is a standard bulk PDE.

We expect that our results apply also to the case of Dirichlet boundary conditions however in this setting one must also estimate the error due to boundary approximation which we neglect in this work. 
}
  
 The following variational formulation of \mathref{eqn:pde} was derived in \cite{dziuk2007finite} \changes{and makes use of the transport formula \eqref{eqn:transport_scalar} and the integration by parts formula on the evolving surface \cite[(2.2)]{dziuk2007finite}\margnote{ref 2. pt. 3.}};
 \beq\label{eqn:wf}
\frac{\diff}{\diff t}\int_\Gt u\varphi+\int_\Gt\nabla_\Gt u\cdot\nabla_\Gt \varphi=\int_\Gt u\mdt{\vec v} \varphi,
\eeq
where $\varphi$ is a sufficiently smooth test function defined on the space-time surface 
\[
\calGT:=\bigcup_{t\in[0,T]}\Gt\times\{t\}.
\]
In  \cite{dziuk2007finite}, a  finite element approximation was proposed for (\ref{eqn:wf}) using piecewise linear finite elements
on a triangulated surface interpolating (at the nodes) $\Gt$, the vertices of the triangulated surface were moved with
the material velocity (of points on $\Gt$) $\vec v$. In this work we adopt a similar setup, in that we propose a  finite element approximation using 
piecewise linear finite elements on a triangulated surface interpolating (at the nodes) $\Gt$, however we move the 
vertices of the triangulated surface with the velocity $\vec v_a=\vec v +\vec a_\tangent$, where $\vec a_\tangent$ is an
{\it arbitrary tangential velocity field} ($\vec a_\tangent\cdot\normal=0$). 
Furthermore we assume that $\vec v_a$ satisfies the same smoothness assumptions as the material velocity $\vec v$, i.e., there exits a diffeomorphism $\tilde{\vec G}(\cdot,t):\G^0\to\Gt,$ satisfying $\tilde{\vec G}\in\Cont{2}([0,T],\Cont{2}(\G^0))$ with $\vec v_a(\tilde{\vec G}(\cdot,t),t)=\pdt \tilde{\vec G}(\cdot,t)$ and with $\tilde{\vec G}(\cdot,0)=\vec {I}$ (the identity) and $\vec v_a(\cdot,t) \in \Cont{2}(\Gt)$.
\changes{
We remark, that if $\Gt$ has a boundary this assumption implies that the arbitrary velocity  $\vec a_\tangent$ has zero conormal component on the boundary, i.e., for $t\in[0,T]$,\margnote{ref 2. pt. 4.}
\beq\label{eqn:ALE_conormal}
(\vec v-\vec v_a)\cdot\conormal=0 \quad\text{ on }\partial\Gt.
\eeq
}

For a sufficiently smooth function $f$, we have that
\begin{align*}
\mdt{\vec v_a}f&=\pdt f+ \vec v_a\cdot\nabla f=\pdt f +\vec v_\normal\cdot\nabla f+(\vec a_\tangent  +\vec v_\tangent) \cdot\nabla_\Gt f=\mdt{\vec v} f +\vec a_\tangent\cdot\nabla_\Gt f.
\end{align*}
Thus we may write the following equivalent variational formulation to \mathref{eqn:pde}, which will form the basis for our finite element approximation
 \beq\label{eqn:ale_wf}
\frac{\diff}{\diff t}\int_\Gt u\varphi+\int_\Gt\nabla_\Gt u\cdot\nabla_\Gt \varphi=\int_\Gt\left( u\mdt{\vec v_a} \varphi-u\vec a_\tangent\cdot\nabla_\Gt \varphi\right), 
\eeq
where $\varphi$ is a sufficiently smooth test function defined on $\calGT$. We note that (\ref{eqn:ale_wf}) may be thought of as a weak formulation of an advection diffusion-equation on a surface with material velocity $\vec v_a$, in which the advection $\vec a_\tangent$ is governed by some external process other than material transport. Hence the results we present are also an analysis of a numerical scheme for an advection diffusion equation on an evolving surface with another source of advective transport other than that due to the material velocity.

The original  (Lagrangian) evolving surface finite element method (ESFEM) was proposed and analysed in \cite{dziuk2007finite}, where optimal  error bounds were shown for the error in the energy norm in the semidiscrete (space discrete) case. Optimal $\Lp{2}$ error bounds for the semidiscrete case were  shown in \cite{dziuk2010l2} and  an optimal error bound for the full discretisation was shown in \cite{doi:10.1137/110828642}. High order Runge-Kutta time discretisations and BDF timestepping schemes for the ESFEM were analysed in \citep{dziuk2011runge}  and \cite{lubich2013backward} respectively. There has also been recent work on the analysis of ESFEM approximations of    the Cahn-Hilliard equation on an evolving surface \cite{2013arXiv1310.4068E}, scalar conservation laws on evolving surfaces \cite{2013arXiv1307.1056D} and the wave equation on an evolving surface \cite{lubich2012variational}. For an overview of finite element methods for PDEs on fixed and evolving surfaces see \cite{dziuk2013finite}. Although the analytical results have thus far focussed on the case where the discrete velocity is an interpolation of the continuous material velocity, the Lagrangian setting, in many applications it proves computationally efficient to consider a mesh velocity which is different to the interpolated material velocity. In particular it appears that the arbitrary tangential velocity, that we consider in this study can be chosen such that beneficial mesh properties are observed in practice. This provides the motivation for this work in which we analyse an ESFEM where the material velocity of the mesh is different to (the interpolant of)  the material velocity of the surface, i.e., an arbitrary Lagrangian-Eulerian ESFEM (ALE-ESFEM). We refer to \cite{EllSty12} for extensive computational investigations of the ALE-ESFEM that we analyse in this study. For examples  in the numerical simulation  of mathematical models for cell motility and biomembranes, where the ALE approach proves computationally more robust than the Lagrangian approach, we refer to \cite{neilson2010use,elliott2012CiCP,2013arXiv1311.7602C,venk11chemotaxis}. 

 Our main results  are Theorems \ref{the:sd_convergence} and \ref{the:BDF2_fd_convergence} where we show optimal order error bounds for the semidiscrete (space discrete, time continuous) and fully discrete numerical schemes. The fully discrete bound is proved for a second order backward difference time discretisation. An optimal error bound is also stated for an implicit Euler time discretisation.  While the fully discrete bound is proved independently of the bound on the semidiscretisation, we believe that the analysis of the semidiscrete scheme may prove a useful starting point for the analysis of other time discretisations. We also observe that, under suitable assumptions on the evolution, the analysis holds for smooth flat surfaces, i.e, bulk domains with smooth boundary. Thus the analysis we present is also an analysis of ALE schemes for PDEs in evolving bulk domains.

 We report on numerical simulations of the fully discrete scheme that support our theoretical results and illustrate that the arbitrary tangential velocity may be chosen such that the meshes generated during the evolution are more suitable than in the Lagrangian case.   \changes{Proposing a general recipe for choosing the tangential velocity is a challenging task that is beyond the scope of this article and we do not address this issue. Moreover the choice of the tangential velocity and what constitutes a good computational mesh is likely to depend heavily on the specific application.\margnote{ref 2. major pt. 4.}} We also investigate numerically the long time behaviour of solutions to (\ref{eqn:pde}) with different initial data when the evolution of the surface is a periodic function of time. Our numerical results indicate that in the example we consider the solution converges to the same periodic solution for different initial data. 
  
The original ESFEM was formulated  for a surface with a smooth material velocity that had both normal and tangential components \cite{dziuk2007finite}. Hence many of the results from the literature are applicable in the present setting of a smooth arbitrary tangential velocity. 
\section{Setup}\label{sec:setup}
We start by introducing an abstract  notation in which we formulate the problem.  
\begin{Defn}[Bilinear forms]\label{def:bf}
\changes{
For $\varphi,\psi\in\Hil{1}{(\Gt)},\vec w\in\Cont{2}{(\Gt)}$ we define the following bilinear forms
\begin{align*}
\abil{\varphi(\cdot,t)}{\psi(\cdot,t)}&=\int_\Gt\nabla_\Gt\varphi(\cdot,t)\cdot\nabla_\Gt\psi(\cdot,t)\\
\mbil{\varphi(\cdot,t)}{\psi(\cdot,t)}&=\int_\Gt\varphi(\cdot,t)\psi(\cdot,t)\\
\gbil{\varphi(\cdot,t)}{\psi(\cdot,t)}{\vec w(\cdot,t)}&=\int_\Gt\varphi(\cdot,t)\psi(\cdot,t)\nabla_\Gt\cdot\vec w(\cdot,t)\\
\bbil{\varphi(\cdot,t)}{\psi(\cdot,t)}{\vec w(\cdot,t)}&=\int_\Gt\varphi(\cdot,t)\nabla_\Gt\psi(\cdot,t)\cdot\vec w(\cdot,t)\\
\atbil{\varphi(\cdot,t)}{\varphi(\cdot,t)}{\vec v^a(\cdot,t)}&=\int_\Gt\left(\nabla_\Gt\cdot\vec v^a(\cdot,t)-2\defD{\vec v^a(\cdot,t)}{\Gt}\right)\nabla_\Gt\varphi(\cdot,t)\cdot\nabla_\Gt\psi(\cdot,t)\\
\btbil{\varphi(\cdot,t)}{\psi(\cdot,t)}{\vec w(\cdot,t)}{\vec v^a(\cdot,t)}&=\int_\Gt\nabla_\Gt\cdot\vec v^a(\cdot,t)\left(\varphi(\cdot,t)\vec w(\cdot,t)\cdot\nabla_\Gt \psi(\cdot,t)\right)\\
&-\int_{\Gct}\varphi(\cdot,t)\vec w(\cdot,t)\cdot \defB{\vec v^a(\cdot,t)}{\Gt}\nabla_\Gt \psi(\cdot,t),
\end{align*}
\margnote{ref 2. pt. 10.}
with the deformation tensors $\defB{\cdot}{(\cdot)}$ and $\defD{\cdot}{(\cdot)}$ as defined in Lemma \ref{lem:transport}.
}
\end{Defn}
We may now write the equation (\ref{eqn:ale_wf}) as 
 \beq\label{eqn:ale_vf}
 \frac{\diff }{\diff t}\mbil{u}{\varphi}+\abil{u}{\varphi}=\mbil{u}{\mdt{\vec v_a}\varphi}-\bbil{u}{\varphi}{\vec a_\tangent}.
 \eeq
 In \cite{dziuk2007finite} the authors showed existence of a weak solution to (\ref{eqn:wf}) and hence a weak solution exists to the (reformulated) problem (\ref{eqn:ale_wf}), furthermore for sufficiently smooth initial data the authors proved the following estimate for the solution of (\ref{eqn:wf}) and hence of (\ref{eqn:ale_wf}) 
 \begin{align}
 \sup_{t\in(0,T)}\ltwon{u(\cdot,t)}{\Gt}^2+\int_0^T\ltwon{\nabla_\Gt u}{\Gt}^2\diff t\leq c\ltwon{u_0}{\G^0}^2,\\
 \int_{0}^T\ltwon{\mdt{\vec v}u}{\Gt}^2\diff t+\sup_{t\in(0,T)}\ltwon{\nabla_\G u}{\G}\leq c\Hiln{u_0}{1}{\G^0}^2.\label{eqn:cont_md_pde_lag}
 \end{align}
We immediately conclude that as $ \mdt{\vec v_a} u -\mdt{\vec v}u=\vec a_\tangent\cdot\nabla_\G u$ the bound (\ref{eqn:cont_md_pde_lag}) holds with the material derivative with respect to the material velocity replaced with the material derivative with respect to the ALE-velocity. See \cite{vierling2011control,olshanskii2013eulerian,AlpEllSti}
for further discussion on the well-posedness of the weak formulation of the continuous problem.
\section{Surface finite element discretisation}\label{sec:fe_disc}
\subsection{Surface discretisation}
 The smooth surface $\Gamma(t)$ is interpolated at nodes
$\vec X_j(t)\in\Gt$ ($j=1,\ldots,J$)
by a discrete evolving surface $\Gct$.  These nodes move with velocity $\diff \vec X_j(t)/\diff t=\vec v_a(\vec X_j(t),t)$ and hence the nodes of the discrete surface $\Gct$ remain on the surface $\Gt$ for all $t\in[0,T]$. The discrete surface, 
\begin{equation*}
\Gamma_h(t)=\bigcup_{K(t)\in\mathcal{T}_h(t)}  K(t)
\end{equation*}
is the union of $m$-dimensional  simplices $K(t)$ that is assumed to form an
admissible triangulation $\mathcal{T}_h(t)$; see \cite[\S 4.1]{dziuk2013finite} for details.  We assume that the maximum diameter 
 of the simplices is bounded uniformly in time and we denote this bound by $h$ which we refer to as the mesh size.

\changes{
 We assume that for each point $x$ on $\Gct$ the exists a unique point $\vec p(\vec x,t)$ on $\Gt$ such that for $t\in[0,T]$
 (see \cite[Lemma 2.8]{dziuk2013finite} for sufficient conditions such that this assumption holds)\margnote{ref 2. pt. 5}
 \beq\label{eqn:x_gct_p_gt}
\vec x=\vec p(\vec x,t)+d(\vec x,t)\normal(\vec p(\vec x,t),t),
\eeq
where $d$ is the oriented distance function to $\Gt$ (see \cite[\S 2.2]{dziuk2013finite} for details).

For a continuous function $\eta_h$ defined on $\Gct$  we define its  lift $\eta_h^l$ onto $\Gt$ by extending constantly in the normal direction   $\normal$ (to the continuous surface) as follows
\beq\label{eqn:lift}
\eta_h^l(\vec p,t)=\eta_h(\vec x(\vec p,t),t)\quad\text{for }\vec p\in\Gt,
\eeq
where $\vec x(\vec p,t)$ is defined by (\ref{eqn:x_gct_p_gt}). 

We assume that the triangulated and continuous surfaces are such that for each simplex $K(t)\in\mathcal{T}_h(t)$ there is a unique  $k(t)\subset\Gt$, whose edges are the unique projections of the edges  of $K(t)$ onto $\Gt$. The union of the $k(t)$ induces an exact triangulation of $\Gt$ with curved edges. We refer, for example, to \cite[\S 4.1]{dziuk2013finite} for further details. 

We also find it convenient to introduce the discrete space-time surface
\[
\calGTh:=\bigcup_{t\in[0,T]}\Gct\times\{t\}.
\]
}

\begin{Defn}[Surface finite element spaces]\label{def:FE_space} For each $t\in[0,T]$ we define the finite element spaces together
with their associated lifted finite element spaces 
\begin{align*}
\Sct&=\left\{\Phi\in C^0\left(\Gct\right)|\Phi|_K \mbox{ is linear affine for each }K\in\T_h(t)\right\},\\
\Sclt&=\left\{\varphi=\Phi^l|\Phi\in\Sct\right\}.
\end{align*}
\end{Defn}

Let $\chi_j(\cdot,t)$ ($j=1,\dots,N$) be the nodal basis of $\Sct$, so that, denoting by $\lbrace\vec X_j\rbrace_{j=1}^J$ the vertices of $\Gct$,
$\chi_j(\vec X_i(t),t)=\delta_{ji}$.
The discrete surface moves with the piecewise linear velocity $\vec V^a_h$ and by $\vec T^a_h$ we denote the interpolant of the arbitrary  tangential velocity $\vec a_\tangent$
\begin{align}
\label{eqn:disc_surface_mat_velocity}
\vec V^a_h(\vec x,t)&=\sum_{j=1}^J\vec v_a(\vec X_j(t),t)\chi_j(\vec x,t),\\
\label{eqn:disc_surface_tang_velocity}
\vec  T^a_h (\vec x,t)&=\sum_{j=1}^J\vec a_\tangent(\vec X_j(t),t)\chi_j(\vec x,t).
\end{align}

The discrete surface gradient is defined piecewise on each surface simplex $K(t)\in\mathcal{T}_h(t)$ as
$$\nabla_{\Gamma_h}g=\nabla g - \nabla g \cdot  \normal_h  \normal _h,$$
where $\normal _h$ denotes the normal to the discrete surface defined element wise.

We now relate the material velocity $\vec V^a_h$ of the triangulated surface $\Gc$ to the material velocity $\vec v^a_h$ of the smooth triangulated surface. For each $\vec X(t)$ on $\Gct$ there is a unique $\vec Y(t)=\vec p(\vec X(t),t)\in\Gt$
with
\margnote{ ref 2. pt 6.}
\changes{
\begin{align}\label{eqn:v_h_def}
\frac{\diff}{\diff t}{\vec Y}(t)=\pdt\vec p(\vec X(t),t)+\vec V^a_h(\vec X(t),t)\cdot\nabla\vec p(\vec X(t),t)=:\vec v^a_h(\vec p(\vec X(t),t),t),
\end{align}
}
where $\vec p$ is as in (\ref{eqn:x_gct_p_gt}). We note that $\vec v^a_h$ is not the interpolant of the velocity $\vec v_a$ into the space $\Scl$ (c.f., \cite[Remark 4.4]{dziuk2010l2}). We denote by $\vec t^a_h= (\vec T^a_h)^l$ the lift of the velocity $\vec T^a_h$ to the smooth surface.

\changes{
\begin{Rem}[Surfaces with boundary]\label{rem:bdry}\margnote{ref 2. major pt. 5 and pt. 41}
While the method is directly applicable to surfaces with boundary, for the analysis we require  the  lift of the triangulated surface to be the smooth surface i.e., $(\Gct)^l=\Gt$. Thus in general we must allow the faces of elements on the boundary of the triangulated surface to be  curved. For the  natural boundary conditions  we consider it is possible to define a conforming piecewise linear finite element space on a triangulation with curved boundary elements, see for example \cite{brenner2002mathematical}, assuming this setup and neglecting the variational crime committed in integrating over curved faces the analysis we present in the subsequent sections  remains valid. 
However as the surface is evolving a further requirement is that the continuous material velocity $\vec v$ and the material velocity of the smooth triangulated (lifted) surface $\vec v^a_h$ must satisfy
\beq\label{eqn:surf_boundary_assumption}
(\vec v-\vec v_h^a)\cdot\conormal=0 \quad\text{ on }\partial\Gt,
\eeq
where $\conormal$ is the conormal to the surface.
If \eqref{eqn:surf_boundary_assumption} does not hold, the additional  error due to domain approximation must also be estimated. We remark that this issue is not specific to the ALE scheme we consider in this study and   also arises if we take $\vec a_\tangent=\vec 0$, i.e., the Lagrangian setting.

 We note that although restrictive the above assumptions are satisfied for some nontrivial examples that are of interest in applications. In \S  \ref{Sec:examples} we present two such examples. In  Example \ref{eg:benchmark}, we present an example of a moving surface with boundary where the lift of the polyhedral surface (with straight boundary faces) is the smooth surface.  In Example \ref{eg:graph}, we present an example where the surface is the graph of a time dependent function over the unit disc. Here the boundary curve is fixed and the boundary edges of elements on the boundary of the triangulated surface are curved such that the triangulation of the boundary is exact.
 \end{Rem}

\begin{Rem}[ALE schemes for PDEs posed in bulk domains]\label{Rem:bdry_2}\margnote{ref 2. pt 42.}
As a special case of a surface with boundary, the method is applicable  to a moving domain in $n$  $(n=1,2,3)$ dimensions, that is when $\Gt$ is a flat (i.e., the normal to the surface is constant) $n$ dimensional hypersurface in $\Reals^{n+1}$ with smooth boundary.  We note that the formulae for the discrete schemes (\ref{eqn:sd_scheme}), (\ref{eqn:BDF2_fd_scheme}) and (\ref{eqn:fd_scheme}) are the same as in the case of a curved hypersurface. Under suitable assumptions on the velocity at the boundary, the analysis we present is valid in this setting. Specifically the analysis we present is valid for a flat three dimensional surface with zero normal velocity but nonzero tangential (and conormal) velocity (subject to \eqref{eqn:surf_boundary_assumption}). In this case, as the domain is flat the geometric errors we estimate in the subsequent sections are zero (as $\Gct=(\Gct)^l=\Gt$ since the lift is in the normal direction only). We note that this assumption necessitates curved boundary elements in this case.  Therefore, as a consequence of our analysis we get  an error estimate for an ALE scheme for a linear parabolic equation on an evolving three dimensional bulk domain in which all of the analysis is all carried out on the physical domain. 
  \end{Rem}
}


\subsection{Material derivatives}

\changes{
We introduce the normal time derivative $\ntd$ on a surface moving with material velocity $\vec v$ defined by
$$
\ntd \eta :=\pdt \eta +\vec v\cdot\normal\normal\cdot\nabla\eta,
$$
and define the space \margnote{Needed for time regularity ..}
$$\Hil{1}{(\calGT)}:=\left\{\eta\in\Lp{2}(\calGT)\lv\nabla_\G\eta\in\Lp{2}(\calGT)\rv\ntd \eta\in \Lp{2}(\calGT)\right\}.$$
}

We are now in a position to define material derivatives on the triangulated surfaces. Given the velocity field $\vec V^a_h\in(\Sc)^{m+1}$ and the associated velocity $\vec v^a_h$ on $\Gt$ we define discrete material derivatives on $\Gct$ and $\Gt$ element wise as follows,\changes{\margnote{ref 2. pt. 7} for $\Phi_h(\cdot,t)\in\Sct$ with $\ntd \Phi_h\vert_{K(t)}\in\Lp{2}(K(t))$ and $\varphi(\cdot,t)\in\Hil{1}{(\calGT)}$,}
\begin{align}
\mdth{\vec V^a_h}\Phi_h\vert_{K(t)}&=\left( \pdt \Phi_h+\vec V^a_h\cdot\nabla\Phi_h\right)\vert_{K(t)},\\
\mdth{\vec v^a_h}\varphi\vert_{k(t)}&=\left( \pdt \varphi+\vec v^a_h\cdot\nabla\varphi\right)\vert_{k(t)}.
\end{align}

\changes{ We find it convenient to introduce the spaces   \margnote{Needed for time regularity ..}
\[
\ScT:=\left\{\Phi_h\text{ and }\mdt{\vec V^a_h}\Phi_h\in\Cont{0}(\calGTh)\vert\Phi_h(\cdot,t)\in\Sct,t\in[0,T]\right\}
\]
and
\[
\ScTl:=\left\{\varphi_h\text{ and }\mdt{\vec v^a_h}\varphi_h\in\Cont{0}(\calGT)\vert\varphi_h(\cdot,t)\in\Sclt,t\in[0,T]\right\}.
\]
\margnote{Check second space}
}

The following transport property of the finite element basis  functions was shown in \cite[\S 5.2]{dziuk2007finite}
\beq\label{eqn:basis_transport}
\mdth{\vec V^a_h}\chi_j=0,\quad\mdth{\vec v^a_h}\chi_j^l=0,
\eeq
which implies that for $\Phi_h=\sum_j\Phi_j(t)\chi_j(\cdot,t)\in\Sct$ with $\varphi_h=\Phi_h^l\in\Sclt$
\beq
\mdth{\vec V^a_h}\Phi_h(\cdot,t)=\sum_{j=1}^J\left(\frac{\diff}{\diff t}{\Phi}_j(t)\right)\chi_j(\cdot,t),\quad \mdth{\vec v^a_h}\varphi_h(\cdot,t)=\sum_{j=1}^J\left(\frac{\diff}{\diff t}{\Phi}_j(t)\right)\chi_j(\cdot,t)^l.
\eeq

We now introduce the notation we need  to formulate and analyse the fully discrete scheme. 
Let $N$ be a positive integer, we define the uniform timestep $\tau=T/N$. For each $n\in\lbrace 0,\dots,N\rbrace$ we set $t^n=n\tau$. \changes{We also occasionally  use the same shorthand for time dependent objects, e.g., $\Gtn{n}:=\G(t^n)$ and $\Gctn{n}:=\Gc(t^n)$, $\left(\vec T_a^h\right)^n:=\vec T_a^h(\cdot,t^n)$ etc.
\margnote{ref 2. pts. 22., 25., 36.}
}
For a discrete time sequence $f^n$, $n\in\lbrace 0,\dots,N\rbrace$ we introduce the notation 
\beq
\pdtau f^n=\frac{1}{\tau}\left(f^{n+1}-f^{n}\right).
\eeq
For $n\in\lbrace 0,\dots,N\rbrace$ we denote by $\Scn{n}=\Sc(t^n)$ and by $\Scln{n}=\Scl(t^n)$. For $j=\lbrace1,\dots,J\rbrace$, we set  
\beq
\chi_j^n=\chi_j(\cdot,t^n), \quad \chi_j^{n,l}=\chi^l_j(\cdot,t^n)
\eeq
and employ the notation
\beq
\Phi_h^n=\sum_{j=1}^J\Phi_j^n\chi_j^n\in\Scn{n}, \quad  \varphi_h^n=\Phi^{n,l}\in\Scln{n}.
\eeq
Following  \citep{doi:10.1137/110828642}  we find it convenient to define for $\alpha=-1,0,1$ and $t\in[t^{n-1},t^{n+1}]$
\begin{align}\label{eqn:pbhn_defn}
\bhn{\Phi}{n+\alpha}(\cdot,t)=\sum_{j=1}^J\Phi_j^{n+\alpha}\chi_j(\cdot,t)\in\Sct,\\
\bhn{\varphi}{n+\alpha}(\cdot,t)=\left(\bhn{\Phi}{n+\alpha}(\cdot,t)\right)^l\in\Sclt.
\end{align}
We now introduce a concept of material derivatives for time discrete functions as defined in \cite{doi:10.1137/110828642}.
Given $\Phi_h^n\in\Scn{n}$ and  $\Phi_h^{n+1}\in\Scn{n+1}$ we define the time discrete material derivative as follows
\changes{
\margnote{ref 2. pt 9. notation changed}
\beq
\mdthn{\vec V^a_h}\Phi_h^n=\sum_{j=1}^J\pdtau\Phi_j^n\chi_j^n\in\Scn{n},\qquad\mdthn{\vec v^a_h}\varphi_h^n=\sum_{j=1}^J\pdtau\Phi_j^n\chi_j^{n,l}\in\Scln{n}.
\eeq
}
The following observations are taken from \cite[\S 2.2.3]{doi:10.1137/110828642}, for $n\in{0,\dotsc,N-1}$
\changes{
\margnote{md of time discrete object}
\beq\label{eqn:td_basis_transport}
\mdthn{\vec V^a_h}\chi_j^n=0,\qquad\mdthn{\vec v^a_h}\chi_j^{n,l}=0.
\eeq
}
On $[t^{n-1},t^{n+1}]$, for $\alpha=-1,0,1$
\beq\label{eqn:td_pb_transport}
\mdth{\vec V^a_h}\bhn{\Phi}{n+\alpha}=0,\quad\mdth{\vec v^a_h}\bhn{\varphi}{n+\alpha}=0,
\eeq
which implies 
\changes{
\beq\label{eqn:pbn_n}
\bhn{\Phi}{n+1}(\cdot,t^n)=\Phi_h^n+\tau\mdthn{\vec V^a_h}\Phi_h^n,\quad\bhn{\varphi}{n+1}(\cdot,t^n)=\varphi_h^n+\tau\mdthn{\vec v^a_h}\varphi_h^n.
\eeq
\margnote{ref 2. pt 8.}
We will also make use of the following notation, for $n\in\lbrace 0\dots,N-1\rbrace$ given $\Phi_h^n\in\Scn{n}$ and  $\Phi_h^{n+1}\in\Scn{n+1}$, with lifts $\varphi_h^n\in\Scln{n}$ and  $\varphi_h^{n+1}\in\Scln{n+1}$, we define $\Phi_h^L\in\Sct$ and $\varphi_h^L\in\Sclt$, $t\in[0,T]$ such that for  $t\in[t^n,t^{n+1}]$
\begin{align}\label{eqn:phi_hl}
\Phi_h^L(\cdot,t)=\frac{t^{n+1}-t}{\tau}\bhn{\Phi}{n}(\cdot,t)+\frac{t-t^{n}}{\tau}\bhn{\Phi}{n+1}(\cdot,t),\\
\varphi_h^L(\cdot,t)=\frac{t^{n+1}-t}{\tau}\bhn{\varphi}{n}(\cdot,t)+\frac{t-t^{n}}{\tau}\bhn{\varphi}{n+1}(\cdot,t).
\end{align}
We note that \eqref{eqn:td_pb_transport} implies
\begin{align}
\mdth{\vec V_h^a}\Phi_h^L(\cdot,t)=\frac{1}{\tau}\left(\bhn{\Phi}{n+1}(\cdot,t)-\bhn{\Phi}{n}(\cdot,t)\right),\\
\mdth{\vec v_h^a}\varphi_h^L(\cdot,t)=\frac{1}{\tau}\left(\bhn{\varphi}{n+1}(\cdot,t)-\bhn{\varphi}{n}(\cdot,t)\right).
\end{align}

}

%
\begin{Defn}[Discrete bilinear forms]\label{def:bf_gamma_h}
We define the analogous bilinear forms to those defined in Definition \ref{def:bf} as follows, for $\Phi_h\in\Sct$, $\Psi_h\in\Sct$ and
$\vec W_h\in(\Sct)^{m+1}$ 
\begin{align*}
\ahbil{\Phi_h(\cdot,t)}{\Psi_h(\cdot,t)}&=\int_\Gct\nabla_\Gct\Phi_h(\cdot,t)\cdot\nabla_\Gct\Psi_h(\cdot,t)\\
\mhbil{\Phi_h(\cdot,t)}{\Psi_h(\cdot,t)}&=\int_\Gct\Phi_h(\cdot,t)\Psi_h(\cdot,t)\\
\ghbil{\Phi_h(\cdot,t)}{\Psi_h(\cdot,t)}{\vec W_h(\cdot,t)}&=\int_\Gct\Phi_h(\cdot,t)\Psi_h(\cdot,t)\nabla_\Gct\cdot\vec W_h(\cdot,t)\\
\bhbil{\Phi_h(\cdot,t)}{\Psi_h(\cdot,t)}{\vec W_h(\cdot,t)}&=\int_\Gct\Phi_h(\cdot,t)\vec W_h(\cdot,t)\cdot\nabla_\Gct\Psi_h(\cdot,t)\\
\ahtbil{\Phi_h(\cdot,t)}{\Psi_h(\cdot,t)}{\vec V^a_h(\cdot,t)}&=\\
\int_\Gct\big(\nabla_\Gc\cdot\vec V^a_h(\cdot,t)-2&\defD{\vec V^a_h(\cdot,t)}{\Gc}\big)\nabla_\Gct\Phi_h(\cdot,t)\cdot\nabla_\Gct\Psi_h(\cdot,t)\\
\bhtbil{\Phi_h(\cdot,t)}{\Psi_h(\cdot,t)}{\vec W_h(\cdot,t)}{\vec V^a_h(\cdot,t)}&=
\int_\Gct\nabla_\Gc\cdot\vec V^a_h(\cdot,t)\left(\Phi\vec W_h(\cdot,t)\cdot\nabla_\Gc \Psi_h(\cdot,t)\right)\\
&-\int_{\Gct}\Phi(\cdot,t)\vec W_h(\cdot,t)\cdot \defB{\vec V^a_h(\cdot,t)}{\Gc}\nabla_\Gc \Psi_h(\cdot,t),
\end{align*}
\changes{
\margnote{ref 2. pt. 10.}
with the deformation tensors $\defB{\cdot}{(\cdot)}$ and $\defD{\cdot}{(\cdot)}$ as defined in Lemma \ref{lem:transport}.
}
\end{Defn}

\subsection{Transport formula}
We recall some results proved in \citep{dziuk2010l2} and \cite{doi:10.1137/110828642} that state (time continuous) transport formulas on the triangulated surfaces and define an adequate notion of discrete in time transport formulas and certain corollaries. The proofs of the transport formulas on the lifted surface (i.e., the smooth surface) follow from the formula given in Lemma \ref{lem:transport}, the corresponding proofs on the triangulated surface $\Gc$ follow once we note that we may apply the same    transport formula stated in Lemma \ref{lem:transport} (with the velocity of $\Gc$ replacing the velocity of $\Gt$) element by element.

We note the transport formula are shown for a triangulated surface with a material velocity that is the interpolant  of a velocity that has both normal and tangential components. Hence the formula may be applied directly  to the present setting where the velocity of the triangulated surface $\vec V^a_h$ is the interpolant of the velocity $\vec v_a$.

\begin{Lem}[Triangulated surface transport formula]\label{lem:transport_gamma_h}
\changes{
Let $\Gct$ be an evolving admissible triangulated surface with material velocity $\vec V^a_h$. Then 
for $\Phi_h,\Psi_h,\vec W_h\in\ScT\times\ScT\times(\ScT)^{m+1}$, \margnote{ref 2. pt. 11.}
\begin{align}
\frac{\diff}{\diff t}\mhbil{\Phi_h}{\Psi_h}&=\mhbil{\mdth{\vec V^a_h}\Phi_h}{\Psi_h}+\mhbil{\mdth{\vec V^a_h}\Psi_h}{\Phi_h}+\ghbil{\Phi_h}{\Psi_h}{\vec V^a_h}\label{eqn:transp_mh}\\
\frac{\diff}{\diff t}\ahbil{\Phi_h}{\Psi_h}&=\ahbil{\mdth{\vec V^a_h}\Phi_h}{\Psi_h}+\ahbil{\mdth{\vec V^a_h}\Psi_h}{\Phi_h}+\ahtbil{\Phi_h}{\Psi_h}{\vec V^a_h}\label{eqn:transp_ah}\\
\frac{\diff}{\diff t}\bhbil{\Phi_h}{\Psi_h}{\vec W_h}&=\bhbil{\mdth{\vec V^a_h}\Phi_h}{\Psi_h}{\vec W_h}+\bhbil{\Phi_h}{\mdth{\vec V^a_h}\Psi_h}{\vec W_h}\label{eqn:transp_bh}\\
\notag&+\bhbil{\Phi_h}{\Psi_h}{\mdt{\vec V^a_h}\vec W_h}+\bhtbil{\Phi_h}{\Psi_h}{\vec W_h}{\vec V^a_h}.
\end{align}
Let $\Gt$ be an evolving surface made up of curved elements $k(t)$ whose edges move with velocity $\vec v^a_h$. Then for
$\varphi,\psi,\vec w\in\Hil{1}{(\calGT)}\times\Hil{1}{(\calGT)}\times(\Cont{1}{(\calGT)})^{m+1}$,
\begin{align}
\frac{\diff}{\diff t}\mbil{\varphi}{\psi}&=\mbil{\mdth{\vec v^a_h}\varphi}{\psi}+\mbil{\varphi}{\mdth{\vec v^a_h}\psi}+\gbil{\varphi}{\psi}{\vec v^a_h}\label{eqn:transp_mhl}\\
\frac{\diff}{\diff t}\abil{\varphi}{\psi}&=\abil{\mdth{\vec v^a_h}\varphi}{\psi}+\abil{\varphi}{\mdth{\vec v^a_h}\psi}+\atbil{\varphi}{\psi}{\vec v^a_h}\label{eqn:transp_ahl}\\
\frac{\diff}{\diff t}\bbil{\varphi}{\psi}{\vec w}&=\bbil{\mdth{\vec v^a_h}\varphi}{\psi}{\vec w}+\bbil{\varphi}{\mdth{\vec v^a_h}\psi}{\vec w}\label{eqn:transp_bhl}\\
\notag&+\bbil{\varphi}{\psi}{\mdth{\vec v^a_h}\vec w}+\btbil{\varphi}{\psi}{\vec w}{\vec v^a_h}.
\end{align}
}
\end{Lem}
We  find it convenient to introduce the following notation
for $W_h\in\Sct$ and $w_h\in\Hil{1}{(\Gt)}, t\in[t^{n-1},t^{n+1}]$ and for a given $\Phi_h^{n+1}\in\Scn{n+1}$ and corresponding lift $\varphi_h^{n+1}\in\Scln{n+1}$
\begin{align}\label{eqn:Lh2}
\Lthbil{W_h}{\Phi^{n+1}}=&\frac{3}{2\tau}\Bigg(\mhbil{W_h(\cdot,t^{n+1})}{\bhn{\Phi}{n+1}(\cdot,t^{n+1})}-\mhbil{W_h(\cdot,t^{n})}{\bhn{\Phi}{n+1}(\cdot,t^{n})}\Bigg)\\
\notag
&-\frac{1}{2\tau}\Bigg(\mhbil{W_h(\cdot,t^{n})}{\bhn{\Phi}{n+1}(\cdot,t^{n})}-\mhbil{W_h(\cdot,t^{n-1)}}{\bhn{\Phi}{n+1}(\cdot,t^{n-1})}\Bigg),\\
\label{eqn:L2}
\Ltbil{w_h}{\varphi^{n+1}}=&\frac{3}{2\tau}\Bigg(\mbil{w_h(\cdot,t^{n+1})}{\bhn{\varphi}{n+1}(\cdot,t^{n+1})}-\mbil{w_h(\cdot,t^{n})}{\bhn{\varphi}{n+1}(\cdot,t^{n})}\Bigg)\\
\notag
&-\frac{1}{2\tau}\Bigg(\mbil{w_h(\cdot,t^{n})}{\bhn{\varphi}{n+1}(\cdot,t^{n})}-\mbil{w_h(\cdot,t^{n-1)}}{\bhn{\varphi}{n+1}(\cdot,t^{n-1})}\Bigg).
\end{align}
\changes{
\margnote{ref 2. pt 12.}
The following Lemma defines an adequate notion of discrete in time transport and follows easily from the transport formula (\ref{eqn:transp_mh}) and \eqref{eqn:td_pb_transport}.
}
\begin{Lem}[Discrete in time transport formula]\label{lem:discrete_in_time_transport}
For $W_h\in\Sct$ and $w_h\in\Hil{1}{(\Gt)}, t\in[t^{n},t^{n+1}]$ and for a given $\Phi_h^{n+1}\in\Scn{n+1}$ and corresponding lift $\varphi_h^{n+1}\in\Scln{n+1}$
\begin{align}
\label{eqn:transp_L2_h}
&\Lthbil{W_h}{\Phi_h^{n+1}}=\\
\notag
&\frac{3}{2\tau}\int_{t^n}^{t^{n+1}}\frac{\diff }{\diff t}\mhbil{W_h(\cdot,t)}{\bhn{\Phi}{n+1}(\cdot,t)}\diff t-\frac{1}{2\tau}\int_{t^{n-1}}^{t^{n}}\frac{\diff }{\diff t}\mhbil{W_h(\cdot,t)}{\bhn{\Phi}{n+1}(\cdot,t)}\diff t\\
\notag
&=\frac{3}{2\tau}\int_{t^n}^{t^{n+1}}\mhbil{\mdth{\vec V^a_h}W_h(\cdot,t)}{\bhn{\Phi}{n+1}(\cdot,t)}+\ghbil{W_h(\cdot,t)}{\bhn{\Phi}{n+1}(\cdot,t)}{\vec V^a_h(\cdot,t)}\diff t\\
\notag
&-\frac{1}{2\tau}\int_{t^{n-1}}^{t^{n}}\mhbil{\mdth{\vec V^a_h}W_h(\cdot,t)}{\bhn{\Phi}{n+1}(\cdot,t)}+\ghbil{W_h(\cdot,t)}{\bhn{\Phi}{n+1}(\cdot,t)}{\vec V^a_h(\cdot,t)}\diff t\\
\label{eqn:transp_L2}
&\Ltbil{w_h}{\varphi_h^{n+1}}=\\
\notag
&\frac{3}{2\tau}\int_{t^n}^{t^{n+1}}\frac{\diff }{\diff t}\mbil{w_h(\cdot,t)}{\bhn{\varphi}{n+1}(\cdot,t)}\diff t-\frac{1}{2\tau}\int_{t^{n-1}}^{t^{n}}\frac{\diff }{\diff t}\mbil{w_h(\cdot,t)}{\bhn{\varphi}{n+1}(\cdot,t)}\diff t\\
\notag
&=\frac{3}{2\tau}\int_{t^n}^{t^{n+1}}\mbil{\mdth{\vec v^a_h}w_h(\cdot,t)}{\bhn{\varphi}{n+1}(\cdot,t)}+\gbil{w_h(\cdot,t)}{\bhn{\varphi}{n+1}(\cdot,t)}{\vec v^a_h(\cdot,t)}\diff t\\
&-\frac{1}{2\tau}\int_{t^{n-1}}^{t^{n}}\mbil{\mdth{\vec v^a_h}w_h(\cdot,t)}{\bhn{\varphi}{n+1}(\cdot,t)}+\gbil{w_h(\cdot,t)}{\bhn{\varphi}{n+1}(\cdot,t)}{\vec v^a_h(\cdot,t)}\diff t
\notag
\end{align}
\end{Lem}

For $t\in[t^{n-1},t^{n+1}]$ and $\tau\leq\tau_0$ the following bounds hold. The result was proved for $t\in[t^{n},t^{n+1}]$ \changes{\margnote{ref 2. pt 13.} in \cite[Lemma 3.6]{doi:10.1137/110828642}}. The  proof may be extended for $t\in[t^{n-1},t^{n+1}]$  as $\mdth{\vec V^a_h}\bhn{\Phi}{n+1}=0$ and $\mdth{\vec v^a_h}\bhn{\varphi}{n+1}=0$,\changes{\margnote{ref 2. pt. 27}
\begin{align}
\label{eqn:mh_pb_td}
\lv\mhbil{\Phi_h^{n+1}}{\Phi_h^{n+1}}-\mhbil{\bhn{\Phi}{n+1}(\cdot,t^n)}{\bhn{\Phi}{n+1}(\cdot,t^n)}\rv&\leq c\tau\mhbil{\Phi_h^{n+1}}{\Phi_h^{n+1}},
\end{align}
and for $t\in[t^{n-1},t^{n+1}]$ and $\tau\leq\tau_0$\begin{align}
\label{eqn:pbn_pn_l2}
\ltwon{\bhn{\Phi}{n}(\cdot,t)}{\Gct}\leq c&\ltwon{\Phi_h^{n}}{\Gctn{n}},\quad\ltwon{\bhn{\varphi}{n}(\cdot,t)}{\Gt}\leq c\ltwon{\varphi_h^{n}}{\Gtn{n}},\\
\label{eqn:pbn_pn_grad}
\ltwon{\nabla_{\Gct}\bhn{\Phi}{n}(\cdot,t)}{\Gct}&\leq c\ltwon{\nabla_{\Gctn{n}}\Phi_h^{n}}{\Gctn{n}},\\\ltwon{\nabla_{\Gt}\bhn{\varphi}{n}(\cdot,t)}{\Gt}&\leq c\ltwon{\nabla_{\Gtn{n}}\varphi_h^{n}}{\Gtn{n}}\notag.
\end{align}
}
\changes{
The following Lemma proves useful in the analysis of the fully discrete scheme.
\begin{Lem}
If $\mdth{\vec V^a_h}{\Phi_h}=0$ and $\mdth{\vec V^a_h}{\Psi_h}=0$ then
\begin{align}
\label{eqn:taupdtau_mh}
\Big\vert&\mhbil{\Phi_h(\cdot,t^{k+1})}{\Psi_h(\cdot,t^{k+1})}-\mhbil{\Phi_h(\cdot,t^{k})}{\Psi_h(\cdot,t^{k})}\Big\vert\\
\notag
&\leq
c\int_{t^k}^{t^{k+1}}\mhbil{\Phi_h(\cdot,t)}{\Phi_h(\cdot,t)}^{1/2}\mhbil{\Psi_h(\cdot,t)}{\Psi_h(\cdot,t)}^{1/2}\diff t
\\
\label{eqn:taupdtau_ah}
\Big\vert&\ahbil{\Phi_h(\cdot,t^{k+1})}{\Phi_h(\cdot,t^{k+1})}-\ahbil{\Phi_h(\cdot,t^{k})}{\Phi_h(\cdot,t^{k})}\Big\vert\\
\notag
&\leq
c\int_{t^k}^{t^{k+1}}\ahbil{\Phi_h(\cdot,t)}{\Phi_h(\cdot,t)}\diff t\\
\label{eqn:taupdtau_bh}
\Big\vert&\bhbil{\Phi_h(\cdot,t^{k+1})}{\Psi_h(\cdot,t^{k+1})}{\vec T^a_h(\cdot,t^{k+1})}-\bhbil{\Phi_h(\cdot,t^{k})}{\Psi_h(\cdot,t^{k})}{\vec T^a_h(\cdot,t^{k})}\Big\vert\\
\notag
&\leq
c\int_{t^k}^{t^{k+1}}\mhbil{\Phi_h(\cdot,t)}{\Phi_h(\cdot,t)}^{1/2}\ahbil{\Psi_h(\cdot,t)}{\Psi_h(\cdot,t)}^{1/2}\diff t,
\end{align}
\end{Lem}
\begin{Proof}
The first two estimates (\ref{eqn:taupdtau_mh}) and (\ref{eqn:taupdtau_ah}) are proved in \cite[Lemma 3.7]{doi:10.1137/110828642}.  To prove (\ref{eqn:taupdtau_bh}) we use the transport formula \eqref{eqn:transp_bh} which yields
\begin{align*}
\Big\vert&\bhbil{\Phi_h(\cdot,t^{k+1})}{\Phi_h(\cdot,t^{k+1})}{\vec T^a_h(\cdot,t^{k+1})}-\bhbil{\Phi_h(\cdot,t^{k})}{\Phi_h(\cdot,t^{k})}{\vec T^a_h(\cdot,t^{k})}\Big\vert\\
&\leq\lv \int_{t^k}^{t^{k+1}}\bhbil{\Phi_h(\cdot,t)}{\Psi_h(\cdot,t)}{\mdt{\vec V^a_h(\cdot,t)}\vec T^a_h(\cdot,t)}+\bhtbil{\Phi_h(\cdot,t)}{\Psi_h(\cdot,t)}{\vec T^a_h(\cdot,t)}{\vec V^a_h(\cdot,t)}\diff t\rv\\
& \leq c\int_{t^k}^{t^{k+1}}\ltwon{\Phi_h}{\Gct}\ltwon{\nabla_{\Gct}\Psi_h}{\Gct}\diff t,
\end{align*}
which is the desired estimate.
\end{Proof}
}

\section{Semidiscrete ALE-ESFEM}\label{sec:sd}
\subsection{Semidiscrete scheme}
\changes{
Given $U_h^0\in\Sc(0)$ find $U_h\in \ScT$ such that  $U_h(\cdot,0)=U_h^0$ and for all $\Phi_h\in\ScT$ and $t\in(0,T]$
\margnote{ref 2. pt. 14.}
\begin{equation}\label{eqn:sd_scheme}
\frac{\diff }{\diff t}\mhbil{U_h}{\Phi_h}+\ahbil{U_h}{\Phi_h}=\mhbil{U_h}{\mdth{\vec V^a_h}\Phi_h}-\bhbil{U_h}{\Phi_h}{\vec T^a_h},
\end{equation}
}

By the transport property of the basis functions (\ref{eqn:basis_transport}) we have the equivalent definition
\begin{equation}\label{eqn:sd_scheme_basis_functions}
\frac{\diff }{\diff t}\mhbil{U_h}{\chi_j}+\ahbil{U_h}{\chi_j}=-\bhbil{U_h}{\chi_j}{\vec T^a_h},  \quad U_h(\cdot,0)=U_h^0, \text{ for }j=1,\dotsc,J.
\end{equation}
Thus a matrix vector formulation of the scheme is  given $\vec \alpha(0)$ find a coefficient vector $\vec \alpha(t), t\in(0,T]$
such that
\beq\label{eqn:mv_sds_scheme}
\frac{\diff}{\diff t}\left(\vec{M}(t)\vec \alpha(t)\right)+\left(\vec{S}(t)+\vec{B}(t)\right)\vec{\alpha}(t)=0,
\eeq
\changes{
where 
$\vec{M}(t),\vec{S}(t)$ and $\vec{B}(t)$ are time dependent mass, stiffness and nonsymmetric matrices with coefficients given by \margnote{ref 2. pt. 15.}
\begin{align}
\label{eqn:mass}
& M(t)_{ij}=\int_{\Gct}\chi_i(\cdot,t)\chi_j(\cdot,t),\quad S(t)_{ij}=\int_{\Gct}\nabla_\Gct\chi_i(\cdot,t)\nabla_\Gct\chi_j(\cdot,t),\\
& B(t)_{ij}=\int_{\Gct}\chi_i(\cdot,t)\vec T^a_h(\cdot,t)\cdot\nabla_\Gct\chi_j(\cdot,t)\notag.
\end{align}
Existence and uniqueness of the semidiscrete finite element solution follows easily as the mass matrix is positive definite, the stiffness matrix is positive semidefinite and the nonsymmetric matrix is bounded. \margnote{ref 2. major pt. 2}  
}
\begin{Lem}[Stability of the semidiscrete scheme]\label{Lem:sd_stability}
\changes{
The finite element solution $U_h$ to (\ref{eqn:sd_scheme}) satisfies the following bounds \margnote{ref 2. pt 16.}
\begin{align}
\sup_{t\in[0,T]}\ltwon{U_h}{\Gct}^2+\int_0^T\ltwon{\nabla_{\Gc(s)}U_h}{\Gc(s)}^2\diff s\leq c\ltwon{U_h}{\Gc(0)}^2,	\label{eqn:L2_stability_sd}\\
\sup_{t\in[0,T]}\ltwon{u_h}{\Gt}^2+\int_0^T\ltwon{\nabla_{\G(s)}u_h}{\G(s)}^2\diff s\leq c\ltwon{u_h}{\G^0}^2,	\label{eqn:L2_stability_sd_l}\\
\int_0^T\ltwon{\mdth{\vec V^a_h}U_h}{\Gc(s)}^2\diff s+\sup_{t\in[0,T]}\ltwon{\nabla_\Gct U_h}{\Gct}^2\leq c\Hiln{U_h}{1}{\Gc(0)}^2,\label{eqn:md_stability_sd}\\
\int_0^T\ltwon{\mdth{\vec v^a_h}u_h}{\G(s)}^2\diff s+\sup_{t\in[0,T]}\ltwon{\nabla_\Gt u_h}{\Gt}^2\leq c\Hiln{u_h}{1}{\G^0}^2\label{eqn:md_stability_sd_l}.
\end{align}
}
\end{Lem}
\begin{Proof}
\changes{
\margnote{IBP removed}
We start with (\ref{eqn:L2_stability_sd}), testing with $U_h$ in (\ref{eqn:sd_scheme}) and applying the transport formula (\ref{eqn:transp_mh}) as in \citep{dziuk2007finite}  yields 
\begin{align*}
\frac{1}{2}\frac{\diff }{\diff t}\mhbil{U_h}{U_h}+\ahbil{U_h}{U_h}&=-\bhbil{U_h}{U_h}{\vec T^a_h}-\frac{1}{2}\ghbil{U_h}{U_h}{\vec V^a_h}.
\end{align*}
 Using Young's inequality to bound the first term on the right hand side and Cauchy-Schwarz on the second term on the right,  we conclude 
\begin{align*}
\frac{1}{2}\frac{\diff }{\diff t}\ltwon{U_h}{\Gct}^2+\ltwon{\nabla_\Gct U_h}{\Gct}^2&\leq\frac{1}{2}\ltwon{\nabla_\Gct U_h}{\Gct}^2+c\ltwon{U_h}{\Gct}^2.
\end{align*}
A Gronwall argument implies the desired result.

For the proof of (\ref{eqn:md_stability_sd}) we  apply the transport formula (\ref{eqn:transp_mh}) to rewrite (\ref{eqn:sd_scheme}) as
\begin{equation*}
\mhbil{\mdth{\vec V^a_h}U_h}{\Phi_h}+\ahbil{U_h}{\Phi_h}=-\ghbil{U_h}{\Phi_h}{\vec V^a_h}-\bhbil{U_h}{\Phi_h}{\vec T^a_h},
\end{equation*}
testing with $\mdth{\vec V^a_h}U_h$ gives
\begin{align}\label{eqn:sd_md_stability_proof_1}
\ltwon{\mdth{\vec V^a_h}U_h}{\Gct}^2+&\ahbil{U_h}{\mdth{\vec V^a_h}U_h}=\\
\notag&-\bhbil{U_h}{\mdth{\vec V^a_h}U_h}{\vec T^a_h}-\ghbil{U_h}{\mdth{\vec V^a_h}U_h}{\vec V^a_h}.
\end{align}
From the transport formulae (\ref{eqn:transp_ah}) and (\ref{eqn:transp_bh}) and  we have
\begin{equation}\label{eqn:sd_md_stability_proof_2}
\ahbil{U_h}{\mdth{\vec V^a_h}U_h}=\frac{1}{2}\left(\frac{\diff}{\diff t}\ahbil{U_h}{U_h}-\ahtbil{U_h}{U_h}{\vec V^a_h}\right),
\end{equation}
and
\begin{align}\label{eqn:sd_md_stability_proof_2_2}
\bhbil{U_h}{\mdth{\vec V^a_h}U_h}{\vec T^a_h}=&\frac{\diff}{\diff t}\bhbil{U_h}{U_h}{\vec T^a_h}-\bhbil{\mdth{\vec V^a_h}U_h}{U_h}{\vec T^a_h}\\
&
-\bhbil{U_h}{U_h}{\mdth{\vec V^a_h}\vec T^a_h}
-\bhtbil{U_h}{U_h}{\vec T^a_h}{\vec V^a_h}.\notag
\end{align}
Using (\ref{eqn:sd_md_stability_proof_2}) and (\ref{eqn:sd_md_stability_proof_2_2}) in (\ref{eqn:sd_md_stability_proof_1}) gives
\begin{align}\label{eqn:sd_md_stability_proof_3}
&\ltwon{\mdth{\vec V^a_h}U_h}{\Gct}^2+\frac{1}{2}\frac{\diff}{\diff t}\ltwon{\nabla_\Gct U}{\Gct}^2+\frac{\diff}{\diff t}\bhbil{U_h}{U_h}{\vec T^a_h}\\
&\quad=\notag\frac{1}{2}\ahtbil{U_h}{U_h}{\vec V^a_h}+\bhbil{\mdth{\vec V^a_h}U_h}{U_h}{\vec T^a_h}\\
&\qquad+\bhbil{U_h}{U_h}{\mdth{\vec V^a_h}\vec T^a_h}+\bhtbil{U_h}{U_h}{\vec T^a_h}{\vec V^a_h}-\ghbil{U_h}{\mdth{\vec V^a_h}U_h}{\vec V^a_h}.\notag
\end{align}
The Cauchy-Schwarz inequality together with the smoothness of the velocity fields $\vec v_a, \vec a_\tangent$ (and hence $\vec V^a_h$ and $\vec T^a_h$),  yields the following estimates
\begin{align}
\ahtbil{U_h}{U_h}{\vec V^a_h}&\leq c\ltwon{\nabla_\Gct U_h}{\Gct}^2,\label{eqn:sd_md_stability_proof_4}\\
\bhbil{\mdth{\vec V^a_h}U_h}{U_h}{\vec T^a_h}&=\int_\Gct\mdth{\vec V^a_h}U_h\vec T^a_h\cdot\nabla_\Gct U_h\label{eqn:sd_md_stability_proof_5}\\
&\notag\leq c\ltwon{\mdth{\vec V^a_h}U_h}{\Gct}\ltwon{\nabla_\Gct U_h}{\Gct}\\
\bhbil{U_h}{U_h}{\mdth{\vec V^a_h}\vec T^a_h}&\leq c\ltwon{U_h}{\Gct}\ltwon{\nabla_\Gct U_h}{\Gct}\label{eqn:sd_md_stability_proof_6}\\
\bhtbil{U_h}{U_h}{\vec T^a_h}{\vec V^a_h}&\leq c\ltwon{U_h}{\Gct}\ltwon{\nabla_\Gct U_h}{\Gct}\label{eqn:sd_md_stability_proof_7}\\
\ghbil{U_h}{\mdth{\vec V^a_h}U_h}{\vec V^a_h}&\leq c\ltwon{U_h}{\Gct}\ltwon{\mdth{\vec V^a_h}U_h}{\Gct}\label{eqn:sd_md_stability_proof_8}
\end{align}
Applying estimates (\ref{eqn:sd_md_stability_proof_4})---(\ref{eqn:sd_md_stability_proof_8}) in (\ref{eqn:sd_md_stability_proof_3}) gives
\begin{align*}
\ltwon{\mdth{\vec V^a_h}U_h}{\Gct}^2&+\frac{1}{2}\frac{\diff}{\diff t}\ltwon{\nabla_\Gct U}{\Gct}^2+\frac{\diff}{\diff t}\bhbil{U_h}{U_h}{\vec T^a_h}\\
\leq&c\ltwon{\nabla_\Gct U_h}{\Gct}^2+c\ltwon{U_h}{\Gct}
\ltwon{\nabla_\Gct U_h}{\Gct} \\
&+c\ltwon{\mdth{\vec V^a_h}U_h}{\Gct}\left(\ltwon{U_h}{\Gct}+\ltwon{\nabla_\Gct U_h}{\Gct}\right)
\end{align*}
Integrating in time and applying (weighted)
Young's inequalities to bound  the third term on the left hand side and the terms on the third line yields  for $t\in[0,T]$,
\begin{align*}
\int_0^t&\ltwon{\mdth{\vec V^a_h}U_h}{\Gc(s)}^2\diff s+\ltwon{\nabla_\Gct U}{\Gct}^2\\
&\leq c\Hiln{U_h}{1}{\Gc(0)}^2+c(\eps)\Bigg(\ltwon{U_h}{\Gct}^2+\int_0^t\ltwon{U_h}{\Gc(s)}^2+
\ltwon{\nabla_{\Gc(s)} U_h}{\Gc(s)}^2\diff s\Bigg),
 \end{align*}
 the estimate (\ref{eqn:L2_stability_sd}) and a
Gronwall argument completes the proof of (\ref{eqn:md_stability_sd}).

Due to the equivalence of the $\Lp{2}$ norm and the $\Hil{1}$ seminorm on $\Gc$ and $\Gt$ (c.f., \cite{dziuk2007finite}), the estimates (\ref{eqn:L2_stability_sd}) and (\ref{eqn:md_stability_sd}) imply the estimates (\ref{eqn:L2_stability_sd_l}) and (\ref{eqn:md_stability_sd_l}) respectively.
}
\end{Proof}

\begin{The}[Error bound for the semidiscrete scheme]\label{the:sd_convergence}
Let $u$ be a sufficiently smooth solution of (\ref{eqn:pde}) and  let the geometry be sufficiently regular.
Furthermore let $u_h(t), t\in[0,T]$ denote the lift of the solution of the semidiscrete scheme (\ref{eqn:sd_scheme}). Furthermore, assume that initial data is sufficiently smooth and  approximation of the initial data is such that
\begin{equation}
\ltwon{u(\cdot,0)-\Ritz u(\cdot,0)}{\G^0}+\ltwon{\Ritz u(\cdot,0)-u_h(\cdot,0)}{\G^0}\leq ch^2,
\end{equation}
holds.
\changes{
Then for $0<h\leq h_0$ with $h_0$ dependent on the data of the problem, the following error bound holds \margnote{ref 2. pt 19.}
\begin{align}
\sup_{t\in(0,T)}\ltwon{u(\cdot,t)-u_h(\cdot,t)}{\Gt}^2+h^2\int_0^T\ltwon{\nabla_\G\left(u(\cdot,t)-u_h(\cdot,t)\right)}{\Gt}^2 \diff t\\ \leq ch^4\sup_{t\in(0,T)}\left(\Hiln{u}{2}{\Gt}^2+\Hiln{\mdt{\vec v_a}u}{2}{\Gt}^2\right).\notag
\end{align}
}
\end{The}
\subsection{Error decomposition}\label{subsec:semidisc_error_decomp}
It is convenient in the analysis to decompose the error as follows
\begin{equation}\label{eqn:semidisc_error_decomp}
u-u_h=\rho+\theta, \quad \rho:=u-\Ritz u,\quad \theta=\Ritz u-u_h\in\Scl,
\end{equation}
with $\Ritz$ the Ritz projection defined in (\ref{eqn:RP_definition}).
\changes{
\begin{Rem}[Applicability of the Ritz projection error bounds]\label{Rem:RP_mass}
In Lemma  \ref{Lem:RP_bounds} we state estimates of the error between a function and its Ritz projection for the case that the function has mean value zero. We note that the solution $u$ to \eqref{eqn:pde} satisfies $\int_\Gt u=\int_{\G^0} u^0$ and
 from the proof of \citep[Thm. 6.1 and Thm. 6.2]{dziuk2010l2} it is clear the bounds remain valid for a function that has a constant mean value (with the  Ritz projection defined by \eqref{eqn:RP_definition} with $\int_\G \Ritz u=\int_\G u$). More generally if we insert a source term $f$ in the right hand side of \eqref{eqn:pde} then the conservation reads $\int_\Gt u=\int_{\G^0} u^0+\int_0^t\int_{\G(s)}f(\cdot,s)\diff s$. Thus if the mean value of $f$ is smooth in time the bounds  remain valid and without loss of generality we may assume the mean value of $f$ is zero.   
  \margnote{ ref 2. pt 20.}
\end{Rem}
}
We shall prove some preliminary Lemmas before proving the Theorem.
\begin{Lem}[Semidiscrete error relation]
We have the following error relation between the semidiscrete solution and the Ritz projection.\changes{ For $\varphi_h\in\ScTl$\margnote{ref 2. pt. 21.}}
\begin{equation}\label{eqn:sd_error_relation}
\frac{\diff }{\diff t}\mbil{\theta}{\varphi_h}+\abil{\theta}{\varphi_h}-\mbil{\theta}{\mdth{\vec v^a_h}\varphi_h}+\bbil{\theta}{\varphi_h}{\vec t^a_h}=F_2(\varphi_h)-F_1(\varphi_h),
\end{equation}
where 
\begin{align}
F_1(\varphi_h)&=
\mbil{\mdth{\vec v^a_h}u_h}{\varphi_h}-\mhbil{\mdth{\vec V^a_h}U_h}{\Phi_h}
+\abil{u_h}{\varphi_h}-\ahbil{U_h}{\Phi_h}\\
&-\bhbil{U_h}{\Phi_h}{\vec T^a_h}+\bbil{u_h}{\varphi_h}{\vec t^a_h}
+\gbil{u_h}{\varphi_h}{\vec v^a_h}-\ghbil{U_h}{\Phi_h}{\vec V^a_h}\notag,\\
F_2(\varphi_h)&=\mbil{-\mdth{\vec v^a_h}\rho}{\varphi_h}-\gbil{\rho}{\varphi_h}{\vec v^a_h}-\bbil{\rho}{\varphi_h}{\vec t^a_h}\\
&+\mbil{u}{\mdt{\vec v_a}\varphi_h-\mdth{\vec v^a_h}\varphi_h}-\bbil{u}{\varphi_h}{\vec a_\tangent-\vec t^a_h}\notag.
\end{align}
\end{Lem}
\begin{Proof}
From the definition of the semidiscrete scheme (\ref{eqn:sd_scheme}) we have
\begin{align}\label{eqn:persd_1}
\frac{\diff }{\diff t}\mbil{u_h}{\varphi_h}+\abil{u_h}{\varphi_h}-\mbil{u_h}{\mdth{\vec v^a_h}\varphi_h}+\bbil{u_h}{\varphi_h}{\vec t^a_h}
&=\\
\frac{\diff }{\diff t}\mbil{u_h}{\varphi_h}+\abil{u_h}{\varphi_h}-\mbil{u_h}{\mdth{\vec v^a_h}\varphi_h}&+\bbil{u_h}{\varphi_h}{\vec t^a_h}\notag\\
-\frac{\diff }{\diff t}\mhbil{U_h}{\Phi_h}-\ahbil{U_h}{\Phi_h}+\mhbil{U_h}{\mdth{\vec V^a_h}\Phi_h}-&\bhbil{U_h}{\Phi_h}{\vec T^a_h}\notag\\
&=F_1(\varphi_h)\notag,
\end{align}
where we have used the transport formulas (\ref{eqn:transp_mh}) and (\ref{eqn:transp_mhl}) for the last step. Using the variational formulation of the continuous equation (\ref{eqn:ale_wf}) we have
\begin{align}\label{eqn:persd_2}
\frac{\diff }{\diff t}&\mbil{\Ritz u}{\varphi_h}+\abil{\Ritz u}{\varphi_h}-\mbil{\Ritz u}{\mdth{\vec v^a_h}\varphi_h}+\bbil{\Ritz u}{\varphi_h}{\vec t^a_h}
\\
=&\frac{\diff }{\diff t}\mbil{\Ritz u}{\varphi_h}+\abil{\Ritz u}{\varphi_h}-\mbil{\Ritz u}{\mdth{\vec v^a_h}\varphi_h}+\bbil{\Ritz u}{\varphi_h}{\vec t^a_h}\notag\\
&-\frac{\diff }{\diff t}\mbil{u}{\varphi_h}-\abil{u}{\varphi_h}+\mbil{u}{\mdt{\vec v_a}\varphi_h}-\bbil{u}{\varphi_h}{\vec a_\tangent}\notag\\
=&\frac{\diff }{\diff t}\mbil{-\rho}{\varphi_h}+\mbil{\rho}{\mdth{\vec v^a_h}\varphi_h}+\mbil{u}{\mdt{\vec v_a}\varphi_h-\mdth{\vec v^a_h}}\notag\\
&-\bbil{\rho}{\varphi_h}{\vec t^a_h}-\bbil{u}{\varphi_h}{\vec a_\tangent-\vec t^a_h}\notag\\
=&F_2(\varphi_h),\notag
\end{align}
where we have used (\ref{eqn:RP_definition}) in the second step and the transport theorem (\ref{eqn:transp_mhl}) in the final step. Subtracting (\ref{eqn:persd_1})
from (\ref{eqn:persd_2}) yields the desired error relation.
\end{Proof}

We estimate the two terms on the right hand side of (\ref{eqn:sd_error_relation}) as follows. From Lemma \ref{lem:geom_pert_errors} we have
\begin{align}
\lv F_1(\varphi_h)\rv
\leq
&ch^2\Big(\ltwon{\mdth{\vec v^a_h}u_h}{\Gt}\ltwon{\varphi_h}{\Gt}+\ltwon{\nabla_\Gt u_h}{\Gt}\ltwon{\nabla_\Gt \varphi_h}{\Gt}\\
&+\ltwon{u_h}{\Gt}\ltwon{\nabla_\Gt \varphi_h}{\Gt}+\Hiln{u_h}{1}{\Gt}\Hiln{\varphi_h}{1}{\Gt}\Big)\notag.
\end{align}
We apply Young's inequality to conclude that with $\epsilon>0$ a positive constant of our choice
\begin{align}\label{eqn:sd_error_bound_epsilon_1}
\lv F_1(\varphi_h)\rv
\leq
&c(\epsilon)h^4\Big(\ltwon{\mdth{\vec v^a_h}u_h}{\Gt}^2  +\Hiln{u_h}{1}{\Gt}^2\Big)
+c(\epsilon)\ltwon{\varphi_h}{\Gt}^2+\epsilon\ltwon{\nabla_\Gt\varphi_h}{\Gt}^2
\end{align}

For the term $F_2$ on the right hand side of (\ref{eqn:sd_error_relation}), we have
\begin{align}
\lv F_2(\varphi_h)\rv\leq& \lv\mbil{-\mdth{\vec v^a_h}\rho}{\varphi_h}\rv +\lv\gbil{\rho}{\varphi_h}{\vec v^a_h}\rv+\lv\bbil{\rho}{\varphi_h}{\vec t^a_h}\rv\\
&+\lv\mbil{u}{\mdt{\vec v_a}\varphi_h-\mdth{\vec v^a_h}\varphi_h}\rv+\lv\bbil{u}{\varphi_h}{\vec a_\tangent-\vec t^a_h}\rv\notag\\
:=&\lv I \rv+ \lv II\rv +\lv III\rv + \lv IV\rv +\lv V\rv.\notag
\end{align}
Using (\ref{eqn:MD_Ritz_bound}) we have
\beq\label{eqn:sd_f2_estim_1}
\lv I\rv\leq\ltwon{\mdth{\vec v^a_h}\rho}{\Gt}\ltwon{\varphi_h}{\Gt}\leq ch^2\left(\Hiln{u}{2}{\Gt}+\Hiln{\mdt{\vec v_a}u}{2}{\Gt}\right)\ltwon{\varphi_h}{\Gt}.
\eeq
We estimate the second and third terms with (\ref{eqn:Ritz_bound}) as follows
\beq\label{eqn:sd_f2_estim_2}
\lv II\rv\leq c\ltwon{\rho}{\Gt}\ltwon{\varphi_h}{\Gt}\leq ch^2\Hiln{ u }{2}{\Gt}\ltwon{\varphi_h}{\Gt},
\eeq
\beq\label{eqn:sd_f2_estim_3}
\lv III\rv\leq c\ltwon{\rho}{\Gt}\ltwon{\nabla_\Gt\varphi_h}{\Gt}\leq ch^2\Hiln{ u }{2}{\Gt}\ltwon{\nabla_\Gt\varphi_h}{\Gt}.
\eeq
For the next term we use (\ref{md_l2_bound}) to conclude
\beq\label{eqn:sd_f2_estim_4}
\lv IV\rv\leq \ltwon{u}{\Gt}\ltwon{\mdt{\vec v_a}\varphi_h-\mdth{\vec v^a_h}\varphi_h}{\Gt}\leq ch^2\ltwon{u}{\Gt}\Hiln{\varphi_h}{1}{\Gt}.
\eeq
Finally for the last term we apply (\ref{tang_velocity_bound}) which yields
\beq\label{eqn:sd_f2_estim_5}
\lv V \rv\leq ch^2\ltwon{u}{\Gt}\ltwon{\nabla_\Gt \varphi_h}{\Gt}.
\eeq
Combining the estimates (\ref{eqn:sd_f2_estim_1})-(\ref{eqn:sd_f2_estim_5}) we have
\begin{align}
\lv F_2(\varphi_h)\rv\leq& ch^2\Bigg( \left(\Hiln{u}{2}{\Gt}+\Hiln{\mdt{\vec v_a}u}{2}{\Gt}\right)\ltwon{\varphi_h}{\Gt}+\\
&\left(\ltwon{u}{\Gt}+\Hiln{u}{2}{\Gt}\right)\ltwon{\nabla_\Gt\varphi_h}{\Gt}+\ltwon{u}{\Gt}\Hiln{\varphi_h}{1}{\Gt}\Bigg)\notag.
\end{align}
We apply Young's inequality to conclude that with $\epsilon>0$ a positive constant of our choice
\beq\label{eqn:sd_error_bound_epsilon_2}
\lv F_2(\varphi_h)\rv\leq c(\epsilon)h^4\left(\Hiln{u}{2}{\Gt}^2+\Hiln{\mdt{\vec v_a}u}{2}{\Gt}^2\right)+c(\epsilon)\ltwon{\varphi_h}{\Gt}^2+\epsilon\ltwon{\nabla_\Gt \varphi_h}{\Gt}^2.
\eeq

\begin{Proof}[of Theorem \ref{the:sd_convergence}]
We test with $\theta$ in the error relation (\ref{eqn:sd_error_relation}) which gives
\beq
\frac{\diff }{\diff t}\mbil{\theta}{\theta}+\abil{\theta}{\theta}-\mbil{\theta}{\mdth{\vec v^a_h}\theta}+\bbil{\theta}{\theta}{\vec t^a_h}=F_2(\theta)-F_1(\theta).
\eeq
Applying the transport formula (\ref{eqn:transp_mhl}) we have
\beq
\frac{1}{2}\frac{\diff }{\diff t}\mbil{\theta}{\theta}+\abil{\theta}{\theta}=F_2(\theta)-F_1(\theta)-\gbil{\theta}{\theta}{\vec v^a_h}-\bbil{\theta}{\theta}{\vec t^a_h}.
\eeq
\changes{
\margnote{IBP removed}
Using a weighted Young's inequality to deal with the last term on the right hand side and application of the estimates (\ref{eqn:sd_error_bound_epsilon_1}) and (\ref{eqn:sd_error_bound_epsilon_2}) and  gives
\begin{align}
\frac{1}{2}&\frac{\diff }{\diff t}\ltwon{\theta}{\Gt}^2+(1-\epsilon)\ltwon{\nabla_\Gt \theta}{\Gt}^2\leq c(\epsilon)\ltwon{\theta}{\Gt}^2 \\
&+ c(\epsilon)h^4\left(\ltwon{\mdth{\vec v^a_h}u_h}{\Gt}^2  +\Hiln{u_h}{1}{\Gt}^2+\Hiln{u}{2}{\Gt}^2+\Hiln{\mdt{\vec v_a}u}{2}{\Gt}^2\right),\notag
\end{align}
with $\epsilon>0$ a positive constant of our choice.}
A Gronwall argument, the stability estimates in Lemma  \ref{Lem:sd_stability}, the error decomposition (\ref{eqn:semidisc_error_decomp}) and the estimates on the error in the Ritz projection (\ref{eqn:Ritz_bound}) complete the proof.
\end{Proof}
\section{Fully discrete ALE-ESFEM}\label{sec:fd}
We consider a second order time discretisation of the semidiscrete scheme (\ref{eqn:sd_scheme})
based on a (second order backward differentiation formula) BDF2 time discretisation defined as follows;
\subsection{Fully discrete BDF2 ALE-ESFEM scheme}
\changes{
\margnote{ref 2. pt. 22.}
Given $U^0_h\in\Scn{0}$ and $U^1_h\in\Scn{1}$ find $U_h^{n+1}\in\Scn{n+1},n\in\{1,\dots,N-1\rbrace$ such that for all $\Phi_h^{n+1}\in\Scn{n+1}$ 
and for $n\in\lbrace1,\dots,N-1\rbrace$
\beq\label{eqn:BDF2_fd_scheme}
\Lthbil{{U}_h^L}{\Phi_h^{n+1}}+\ahbil{U^{n+1}_h}{\Phi^{n+1}_h}=-\bhbil{U_h^{n+1}}{\Phi_h^{n+1}}{(\vec T^a_h)^{n+1}},
\eeq
where we have used the notation introduced in (\ref{eqn:phi_hl}).
}
For the basis functions we note that by definition for $\alpha=-1,0,1$,
 \beq
 \underline{\chi}^{n+1}_j(\cdot,t^{n+\alpha})=\chi_j^{n+\alpha}\in\Scn{n+\alpha}.
 \eeq
Therefore the matrix vector formulation of the scheme (\ref{eqn:BDF2_fd_scheme}) is for $n=\lbrace 1,\dots,N-1\}$ given $\vec U^n,\vec U^{n-1}$ find a coefficient vector $\vec U^{n+1}$ \beq\label{eqn:BDF2_mv_fd_scheme}
\left(\frac{3}{2}\vec{M}^{n+1}+\tau\left(\vec{S}^{n+1}+\vec{B}^{n+1}\right)\right)\vec{U}^{n+1}=2\vec{M}^n\vec{U}^n-\frac{1}{2}\vec M^{n-1}\vec{U}^{n-1},
\eeq
where 
$\vec{M}^n=\vec{M}(t^n),\vec{S}^n=\vec{S}(t^n)$ and $\vec{B}^n=\vec{B}(t^n)$ are time dependent mass, stiffness and nonsymmetric matrices (see (\ref{eqn:mass})).

\begin{Prop}[Solvability of the fully discrete scheme]
For $\tau<\tau_0$, where $\tau_0$ depends on the data of the problem and the arbitrary tangential velocity $\vec a_\tangent$, and for each $n\in\lbrace 2,\dotsc,N\rbrace$, the finite element solution $U_h^n$ to the scheme  (\ref{eqn:BDF2_fd_scheme}) exists and is unique.
\end{Prop}
\begin{Proof}
\changes{\margnote{IBP removed}
Using Young's inequality we have for $\Phi_h^n\in\Scn{n}$
\beq
\lv \bhbil{\Phi_h^n}{\Phi_h^n}{(\vec T^a_h)^n}\rv\leq c(\epsilon)\mhbil{\Phi_h^n}{\Phi_h^n}+\epsilon\ahbil{\Phi_h^n}{\Phi_h^n}.
\eeq
Hence for the scheme (\ref{eqn:BDF2_fd_scheme}) we have for all $\epsilon>0$
\begin{align}
\frac{3}{2}\mhbil{\Phi_h^n}{\Phi_h^n}&+\tau\left(\ahbil{\Phi_h^n}{\Phi_h^n}+\bhbil{\Phi_h^n}{\Phi_h^n}{(\vec T^a_h)^n}\right)\\
\notag
&\geq(\frac{3}{2}-c(\epsilon)\tau)\mhbil{\Phi_h^n}{\Phi_h^n}+\tau(1-\epsilon)\ahbil{\Phi_h^n}{\Phi_h^n},\notag
\end{align}
hence for $\tau\leq\tau_0$, the  system matrix $\vec A^n=\left(\frac{3}{2}\vec{M}^{n}+\tau\left(\vec{S}^{n}+\vec{B}^{n}\right)\right),n=2,\dots,N$ is  positive definite. 
}
\end{Proof}

We now prove the fully discrete analogues to the stability bounds of Lemma \ref{Lem:sd_stability}. 
\changes{
We make use of the  following result from \cite[Lemma 4.1]{dziuk2011runge} that provides basic estimates. There is a constant $\mu$ (independent of the discretisation parameters $\tau, h$ and the length of the time interval $T$) such that for all $\vec \alpha,\vec\beta\in\Reals^J$, for $\tau\leq \tau_0$, for $k,j=-1,0,1,j\geq k$ and for $n\in\{1,\dots,N-1\}$ we have\margnote{ref 2. pt. 24.}
\beq
\label{eqn:DLM_mass}
 \left( \vec M^{n+j}-\vec M^{n+k}\right) \vec \alpha  \cdot  \vec \beta \leq \mu(j-k)\tau\left(\vec M^{n+k}\vec \alpha  \cdot  \vec \alpha\right)^{\frac{1}{2}}\left(\vec M^{n+k}\vec \beta  \cdot  \vec \beta\right)^{\frac{1}{2}}.
 \eeq
 }
 \begin{Lem}[Stability of the fully discrete scheme (\ref{eqn:BDF2_fd_scheme})]\label{Lem:BDF2_fd_stability}
Assume the starting value for the scheme satisfies the bound
\begin{align}
\label{eqn:BDF2_starting_stability_L2}
\ltwon{U_h^1}{\Gctn{1}}^2&\leq c\ltwon{U_h^0}{\Gctn{0}}^2,
\end{align}
then the fully discrete solution $U_h^n,n=2,\dots,N$ of the BDF2 scheme (\ref{eqn:BDF2_fd_scheme})  satisfies the following bounds for $\tau\leq\tau_0$, where $\tau_0$ depends on the data of the problem and the arbitrary tangential velocity $\vec a_\tangent$,
\begin{align}
\label{eqn:BDF2_fd_stability_L2}
\ltwon{U_h^n}{\Gctn{n}}^2+\tau\sum_{i=2}^n\ltwon{\nabla_{\Gctn{i}}U_h^i}{\Gctn{i}}^2&\leq c\ltwon{U_h^0}{\Gctn{0}}^2,\\
\label{eqn:BDF2_fd_stability_L2_l}
\ltwon{u_h^n}{\Gtn{n}}^2+\tau\sum_{i=2}^n\ltwon{\nabla_{\Gtn{i}}u_h^i}{\Gtn{i}}^2&\leq c\ltwon{u_h^0}{\Gtn{0}}^2.
\end{align}
Furthermore if, along with (\ref{eqn:BDF2_starting_stability_L2}), we assume the starting values satisfy the bound
\beq\label{eqn:BDF_starting_MD_stability}
\tau\ltwon{\mdth{\vec V^a_h}U_h^L(\cdot,t^{1}-0)}{\Gctn{2}}^2+\ltwon{\nabla_{\Gctn{1}}U_h^{1}}{\Gctn{1}}^2\leq c\Bigg(\ltwon{U_h^0}{\Gctn{0}}^2+\ltwon{\nabla_{\Gctn{0}}U_h^0}{\Gctn{0}}^2\Bigg),
\eeq
then for $n\in\{2,\dots,N\}$, we have the stability bounds
\begin{align}
\label{eqn:BDF2_md_stability_fd}
\tau\sum_{i=1}^{n-1}\ltwon{\mdth{\vec V^a_h}U_h^L(\cdot,t^{i+1}-0)}{\Gctn{i}}^2+\ltwon{\nabla_{\Gctn{n}}U_h^{n}}{\Gctn{n}}^2&\leq c\Hiln{U_h^0}{1}{\Gctn{0}}^2,\\
\label{eqn:BDF2_md_stability_fd_l}
\tau\sum_{i=1}^{n-1}\ltwon{\mdth{\vec v^a_h}u_h^L(\cdot,t^{i+1}-0)}{\Gtn{i}}^2+\ltwon{\nabla_{\Gtn{n}}u_h^{n}}{\Gtn{n}}^2\leq& c\Hiln{u_h^0}{1}{\Gtn{0}}^2.
\end{align}
\end{Lem}
\begin{Proof}
We begin with the proof of (\ref{eqn:BDF2_fd_stability_L2}). We work with the matrix vector form of the scheme (\ref{eqn:BDF2_mv_fd_scheme}) and we multiply by a vector $\vec U^{n+1}$ which gives
\begin{align}\label{eqn:BDF_l2_stab_pf_1}
&\frac{3}{2\tau}\vec M^{n+1}\vec U^{n+1}\cdot\vec U^{n+1}-\frac{2}{\tau}\vec M^{n}\vec U^{n}\cdot\vec U^{n+1}\\
\notag
&+\frac{1}{2\tau}\vec M^{n-1}\vec U^{n-1}\cdot\vec U^{n+1}+\left(\vec S^{n+1}+\vec B^{n+1}\right)\vec U^{n+1}\cdot\vec U^{n+1}=0.
\end{align}
We first note that a calculation yields for $\vec \alpha,\vec \beta,\vec \kappa\in\Reals^J$
\begin{align}
\left(\frac{3}{2}\vec \alpha-2\vec \beta+\frac{1}{2}\vec \kappa\right)\cdot\vec \alpha=\frac{1}{4}\left(\lv\vec\alpha\rv^2-\lv\vec\beta\rv^2+\lv2\vec \alpha-\vec\beta\rv^2-\lv2\vec\beta-\vec\kappa\rv^2\right)+\frac{1}{4}\lv\vec\alpha-2\vec\beta+\vec\kappa\rv^2.
\end{align}
Using this result we see that
\begin{align}\label{eqn:BDF_l2_stab_pf_2}
\frac{3}{2}\vec M^{n+1}&\vec U^{n+1}\cdot\vec U^{n+1}-2\vec M^{n}\vec U^{n}\cdot\vec U^{n+1}+\frac{1}{2}\vec M^{n-1}\vec U^{n-1}\cdot\vec U^{n+1}\\
\notag
=&
\frac{3}{2}\left(\vec M^{n+1}-\vec M^n\right)\vec U^{n+1}\cdot\vec U^{n+1}+\frac{1}{2}\left(\vec M^{n-1}-\vec M^n\right)\vec U^{n-1}\cdot\vec U^{n+1}\\
\notag
&
+\frac{1}{4}\Bigg(\vec M^n\vec U^{n+1}\cdot\vec U^{n+1}-\vec M^n\vec U^{n}\cdot\vec U^{n}\\
&
\notag
+\vec M^n\left(2\vec  U^{n+1}-\vec U^n\right)\cdot\left(2\vec U^{n+1}-\vec U^n\right)
-\vec M^n\left(2\vec U^{n}-\vec U^{n-1}\right)\cdot\left(2\vec U^{n}-\vec U^{n-1}\right)\\
\notag
&+\vec M^n\left(\vec U^{n+1}-2\vec U^{n}+\vec U^{n-1}\right)\cdot\left(\vec U^{n+1}-2\vec U^{n}+\vec U^{n-1}\right)\Bigg)
\\
=&
\notag
\frac{1}{4}\vec M^{n+1}\vec U^{n+1}\cdot\vec U^{n+1}-\frac{1}{4}\vec M^{n}\vec U^{n}\cdot\vec U^{n}\\
\notag
&+\frac{1}{4}\vec M^{n}\left(2\vec  U^{n+1}-\vec U^n\right)\cdot\left(2\vec U^{n+1}-\vec U^n\right)\\
&
\notag
-\frac{1}{4}\vec M^{n-1}\left(2\vec U^{n}-\vec U^{n-1}\right)\cdot\left(2\vec U^{n}-\vec U^{n-1}\right)\\
\notag
&+\frac{1}{4}\vec M^{n}\left(\vec U^{n+1}-2\vec U^{n}+\vec U^{n-1}\right)\cdot\left(\vec U^{n+1}-2\vec U^{n}+\vec U^{n-1}\right)\\
\notag
&
+\frac{5}{4}\left(\vec M^{n+1}-\vec M^{n}\right)\vec U^{n+1}\cdot\vec U^{n+1}+\frac{1}{2}\left(\vec M^{n-1}-\vec M^{n}\right)\vec U^{n-1}\cdot\vec U^{n+1}\\
&
\notag
+\frac{1}{4}\left(\vec M^{n-1}-\vec M^{n}\right)\left(2\vec U^n-\vec U^{n-1}\right)\cdot\left(2\vec U^n-\vec U^{n-1}\right).
\end{align}
\changes{
The last three terms on the right hand side may be estimated as follows. Using (\ref{eqn:mh_pb_td}) \margnote{ref 2. pt. 267}}
\begin{align}\label{eqn:BDF_l2_stab_pf_3}
\frac{5}{4}\left(\vec M^{n+1}-\vec M^{n}\right)\vec U^{n+1}\cdot\vec U^{n+1}&=\frac{5}{4}\left(\mhbil{U_h^{n+1}}{U_h^{n+1}}-\mhbil{\bhn{U}{n+1}(\cdot,t^{n})}{\bhn{U}{n+1}(\cdot,t^{n})}\right)\\
\notag
&
\geq -c\tau\ltwon{U_h^{n+1}}{\Gctn{n+1}}^2.
\end{align}
Using (\ref{eqn:DLM_mass}), Young's inequality and (\ref{eqn:pbn_pn_l2}) we have
\begin{align}\label{eqn:BDF_l2_stab_pf_4}
\frac{1}{2}\left(\vec M^{n-1}-\vec M^{n}\right)&\vec U^{n-1}\cdot\vec U^{n+1}\\
\notag
&\geq -\frac{\mu}{2}\tau\left(\mhbil{\bhn{U}{n+1}(\cdot,t^{n-1})}{\bhn{U}{n+1}(\cdot,t^{n-1})}+\ltwon{U_h^{n-1}}{\Gctn{n-1}}^2\right)\\
&\notag
\geq
 -c\tau\left(\ltwon{U_h^{n+1}}{\Gctn{n+1}}^2+\ltwon{U_h^{n-1}}{\Gctn{n-1}}^2\right).
\end{align}
\changes{
For the third term we use (\ref{eqn:DLM_mass}) to conclude \margnote{ref 2. pt. 28.,pt. 29.}
\begin{align}\label{eqn:BDF_l2_stab_pf_5}
\frac{1}{4}&\left(\vec M^{n-1}-\vec M^{n}\right)\left(2\vec U^n-\vec U^{n-1}\right)\cdot\left(2\vec U^n-\vec U^{n-1}\right)
\geq
\\
\notag
&
-c\tau\mhbil{2\bhn{U}{n}(\cdot,t^{n-1})-U_h^{n-1}}{2\bhn{U}{n}(\cdot,t^{n-1})-U_h^{n-1}}.
\end{align}
Applying (\ref{eqn:BDF_l2_stab_pf_2})---(\ref{eqn:BDF_l2_stab_pf_5}) in (\ref{eqn:BDF_l2_stab_pf_1}) 
and reverting to the bilinear forms,  we arrive at
\begin{align}
\frac{1}{4}\pdtau&\Big(\mhbil{U_h^n}{U_h^n}+\mhbil{2\bhn{U}{n}(\cdot,t^{n-1})-U_h^{n-1}}{2\bhn{U}{n}(\cdot,t^{n-1})-U_h^{n-1}}\Big)\\
\notag
+&(1-\epsilon)\ltwon{\nabla_{\Gctn{n+1}}U_h^{n+1}}{\Gctn{n+1}}^2\leq  c\Big(c(\epsilon)\ltwon{U_h^{n+1}}{\Gctn{n+1}}^2+\ltwon{U_h^{n-1}}{\Gctn{n-1}}^2\\
\notag
&+\mhbil{2\bhn{U}{n}(\cdot,t^{n-1})-U_h^{n-1}}{2\bhn{U}{n}(\cdot,t^{n-1})-U_h^{n-1}}\Big),
\end{align}
where we have used Young's inequality to bound the non-symmetric term and $\epsilon>0$ is a positive constant of our choice.
Summing over $n$ and multiplying by $\tau$ gives (where we have suppressed the dependence of the constants on $\epsilon$)
\begin{align}
\frac{1}{4}&\Big(\ltwon{U_h^{k}}{\Gctn{k}}^2+\mhbil{2\bhn{U}{k}(\cdot,t^{k-1})-U_h^{k-1}}{2\bhn{U}{k}(\cdot,t^{k-1})-U_h^{k-1}}\Big)\\
\notag
+&\tau\sum_{i=2}^k\ltwon{\nabla_{\Gctn{i}}U_h^{i}}{\Gctn{i}}^2\leq  c\tau\sum_{i=0}^k\ltwon{U_h^{i}}{\Gctn{i}}^2\\
\notag
&+c\tau\sum_{i=1}^k\mhbil{2\bhn{U}{i}(\cdot,t^{i-1})-U_h^{i-1}}{2\bhn{U}{i}(\cdot,t^{i-1})-U_h^{i-1}}\\
\notag
&+\frac{1}{4}\Big(\ltwon{U_h^{1}}{\Gctn{1}}^2+\mhbil{2\bhn{U}{1}(\cdot,t^{0})-U_h^{0}}{2\bhn{U}{1}(\cdot,t^{0})-U_h^{0}}\Big).
\end{align}
With the assumptions on the starting values, a discrete Gronwall argument\margnote{ref 2. pt. 30.} completes the proof.
}
The estimate (\ref{eqn:BDF2_fd_stability_L2_l}) follows by the usual norm equivalence.

\changes{
In order to show the bound (\ref{eqn:BDF2_md_stability_fd}), we recall the following basic identity  given in \cite[pg. 1653]{elliott1993global}, given  vectors $\vec \alpha,\vec \beta, \vec \kappa\in\Reals^J$, \margnote{ref 2. pt 31.}
\begin{align}\label{eqn:Ell_Stu_basic}
\frac{3}{2}\vec \alpha\cdot(\vec \alpha-\vec \beta) -2\vec \beta\cdot(\vec \alpha-\vec \beta)& +\frac{1}{2}\vec \kappa\cdot(\vec \alpha-\vec \beta)
=
\\
\notag
&
\lv\vec \alpha-\vec \beta\rv^2+\frac{1}{4}\Big(\lv\vec \alpha-\vec \beta\rv^2-\lv \vec \beta-\vec \kappa\rv^2+\lv\vec \alpha-2\vec \beta+\vec \kappa\rv^2\Big).
\end{align}
}
\changes{
\margnote{Proof changed to circumvent IBP}
We work with the matrix vector form of the scheme (\ref{eqn:BDF2_mv_fd_scheme}), multiplying with $\vec U^{n+1}-\vec U^{n}$ and using (\ref{eqn:Ell_Stu_basic}) we have\margnote{ref 2. pt 32}
\begin{align}
\frac{1}{\tau}&\Bigg(\vec M^n\left(\vec U^{n+1}-\vec U^{n}\right)\cdot\left(\vec U^{n+1}-\vec U^{n}\right)
+\frac{1}{4}\Big(\vec M^n\left(\vec U^{n+1}-\vec U^{n}\right)\cdot\left(\vec U^{n+1}-\vec U^{n}\right)\\
\notag
&-\vec M^n\left(\vec U^{n}-\vec U^{n-1}\right)\cdot\left(\vec U^{n}-\vec U^{n-1}\right)\\
\notag
&+\vec M^n\left(\vec U^{n+1}-2\vec U^n+\vec U^{n-1}\right)\cdot\left(\vec U^{n+1}-2\vec U^n+\vec U^{n-1}\right)\Big)\Bigg)
\\
\notag
&+\left(\vec S^{n+1}+\vec B^{n+1}\right)\vec U^{n+1}\cdot\left(\vec U^{n+1}-\vec U^{n}\right)
+\frac{1}{2\tau}\left(\vec M^{n-1}-\vec M^{n}\right)\vec U^{n-1}\cdot\left(\vec U^{n+1}-\vec U^{n}\right)
\\
\notag
&+\frac{3}{2\tau}\left(\vec M^{n+1}-\vec M^{n}\right)\vec U^{n+1}\cdot\left(\vec U^{n+1}-\vec U^{n}\right)=0.
\end{align}
Dropping a positive term and rearranging gives\margnote{ref 2. pt 33}
\begin{align}\label{BDF2_stability_md_pf_1}
&\vec M^{n+1}\left(\vec U^{n+1}-\vec U^{n}\right)\cdot\left(\vec U^{n+1}-\vec U^{n}\right)
+\frac{\tau}{4}\pdtau\left(\vec M^n\left(\vec U^{n}-\vec U^{n-1}\right)\cdot\left(\vec U^{n}-\vec U^{n-1}\right)\right)
\\
\notag
&+\frac{\tau}{2}\left(\vec S^{n+1}\vec U^{n+1}\cdot\vec U^{n+1}-\vec S^n\vec U^{n}\cdot\vec U^{n}\right)+\tau\left(\vec B^{n+1}\vec U^{n+1}\cdot\vec U^{n+1}-\vec B^n\vec U^{n}\cdot\vec U^{n}\right)
\\
\notag
&\leq
-\frac{\tau}{2}\vec S^{n+1}(\vec U^{n+1}-\vec U^{n})\cdot(\vec U^{n+1}-\vec U^{n})+\frac{\tau}{2}\left(\vec S^{n+1}-\vec S^{n}\right)\vec U^{n}\cdot\vec U^{n}
\\
\notag
&
+\tau\vec B^{n+1}\left(\vec U^{n+1}-\vec U^n\right)\cdot\vec U^{n}+\tau\left(\vec B^{n+1}-\vec B^{n}\right) \vec U^n\cdot\vec U^n
\\
\notag
&+\frac{1}{2}\left(\vec M^{n}-\vec M^{n-1}\right)\vec U^{n-1}\cdot\left(\vec U^{n+1}-\vec U^{n}\right)
-\frac{3}{2}\left(\vec M^{n+1}-\vec M^{n}\right)\vec U^{n+1}\cdot\left(\vec U^{n+1}-\vec U^{n}\right)
\\
\notag
&+\frac{5}{4}\left(\vec M^{n+1}-\vec M^{n}\right)\left(\vec U^{n+1}-\vec U^{n}\right)\cdot\left(\vec U^{n+1}-\vec U^{n}\right).\\
\notag
&:= I+II+III+IV+V+VI+VII.
\end{align}
For the first two terms on the right hand side of (\ref{BDF2_stability_md_pf_1})
\changes{ we proceed as in \cite[Proof of Lemma 4.1]{doi:10.1137/110828642} using (\ref{eqn:taupdtau_ah}) and (\ref{eqn:pbn_pn_grad}) we get the following bound,\margnote{ref 2. pt. 34.} }  
 \begin{align}\label{eqn:md_fd_stab_pf_1_1_2}
I+II=-\frac{\tau^3}{2}&\ahbil{\mdth{\vec V^a_h}U_h^L(\cdot,t^{n+1}-0)}{\mdth{\vec V^a_h}U_h^L(\cdot,t^{n+1}-0)}
\\
+\frac{\tau}{2}&\Big(\ahbil{\bhn{U}{n}(\cdot,t^{n+1})}{\bhn{U}{n}(\cdot,t^{n+1})}-\ahbil{{U_h^n}}{{U_h^n}}\Big)\notag\\
&\leq c\tau^2\ltwon{\nabla_{\Gctn{n}}U^{n}_h}{\Gctn{n}}^2\notag.
\end{align}
For the third term on the right hand side of (\ref{BDF2_stability_md_pf_1}), we have 
\begin{align}\label{eqn:md_fd_stab_pf_1_3}
III\leq&\big\vert \tau^2\bhbil{\mdth{\vec V^a_h}U_h^L(\cdot,t^{n+1}-0)}{\bhn{U}{n}(\cdot,t^{n+1})}{(\vec T^a_h)^{n+1}}\big\vert\\
\leq& c\tau^2\ltwon{\mdth{\vec V^a_h}U_h^L(\cdot,t^{n+1}-0)}{\Gctn{n+1}}\ltwon{\nabla_{\Gctn{n+1}}\bhn{U}{n}(\cdot,t^{n+1})}{\Gctn{n+1}}\notag\\
\leq&\epsilon\tau^2\ltwon{\mdth{\vec V^a_h}U_h^L(\cdot,t^{n+1}-0)}{\Gctn{n+1}}^2+c(\epsilon)\tau^2\ltwon{\nabla_{\Gctn{n}}U_h^{n}}{\Gctn{n}}^2\notag.
\end{align}
where $\epsilon$ is a positive constant of our choice and we have used Young's inequality and (\ref{eqn:pbn_pn_grad}) in the last step. For the fourth term on the right hand side of (\ref{BDF2_stability_md_pf_1}), we have using (\ref{eqn:taupdtau_bh}),   (\ref{eqn:pbn_pn_l2}) and (\ref{eqn:pbn_pn_grad})
\begin{align}\label{eqn:md_fd_stab_pf_1_4}
IV\leq&\lv \tau\left(\bhbil{\bhn{U}{n}(\cdot,t^{n+1})}{\bhn{U}{n}(\cdot,t^{n+1})}{(\vec T^a_h)^{n+1}}\right)-\left(\bhbil{{U_h}^{n}}{{U_h}^{n}}{(\vec T^a_h)^{n}}\right)\rv\\
\leq& c\tau^2\ltwon{U_h^n}{\Gctn{n}}\ltwon{\nabla_{\Gctn{n}}U_h^n}{\Gctn{n}}\notag.
\end{align}
For the fifth term using (\ref{eqn:DLM_mass}) we 
\beq 
V\leq \mu\tau\left(\vec M^{n-1}\left(\vec U^{n+1}-\vec U^{n}\right)\cdot\left(\vec U^{n+1}-\vec U^{n}\right)\right)^{1/2}\ltwon{U_h^{n-1}}{\Gctn{n-1}}.
\eeq
For the sixth term we use (\ref{eqn:mh_pb_td}) and (\ref{eqn:taupdtau_mh}) to give for all $\epsilon>0$,  \margnote{ref 2. pt. 35}
\begin{align}\label{BDF2_stability_md_pf_3}
 VI&=\frac{3\tau}{2}\left(\mhbil{U_h^{n+1}}{\mdth{\vec V^a_h}U_h^L(\cdot,t^{n+1}-0)}-\mhbil{\bhn{U}{n+1}(\cdot,t^n)}{\mdth{\vec V^a_h}U_h^L(\cdot,t^{n}+0)}\right)
\\
\notag
&\leq c(\epsilon)\tau^2\ltwon{U_h^{n+1}}{\Gctn{n+1}}+\epsilon\tau^2\ltwon{\mdth{\vec V^a_h}U_h^L(\cdot,t^{n+1}-0)}{\Gctn{n+1}}^2.
\end{align}
For the seventh term we apply (\ref{eqn:mh_pb_td}) to obtain \margnote{ref 2. pt. 35}
\beq\label{BDF2_stability_md_pf_4}
 VII\leq c\tau\vec M^{n+1}\left(\vec U^{n+1}-\vec U^{n}\right)\cdot\left(\vec U^{n+1}-\vec U^{n}\right)= c\tau^3\ltwon{\mdth{\vec V^a_h}U_h^L(\cdot,t^{n+1}-0)}{\Gctn{n+1}}^2.
\eeq
Writing (\ref{BDF2_stability_md_pf_1}) in terms of the bilinear forms, applying the estimates (\ref{eqn:md_fd_stab_pf_1_1_2})---(\ref{BDF2_stability_md_pf_4}) and summing gives, 
\begin{align}
\sum_{i=2}^{n}{\tau^2}&\ltwon{\mdth{\vec V^a_h}U_h^L(\cdot,t^{i}-0)}{\Gctn{i}}^2+{\tau}\ltwon{\nabla_\Gctn{n}U_h^n}{\Gctn{n}}
\leq  c\tau^2\ltwon{\mdth{\vec V^a_h}U_h^L(\cdot,t^{1}-0)}{\Gctn{1}}^2
\\
\notag
&+c{\tau}\ltwon{\nabla_\Gctn{1}U_h^1}{\Gctn{1}}+c\tau^2\sum_{i=0}^n\ltwon{U_h^i}{\Gctn{i}}^2+c\tau^2\sum_{i=2}^n\ltwon{\nabla_{\Gctn{i}}U_h^i}{\Gctn{i}}^2.
\end{align}
Dividing by $\tau$, applying the stability bound (\ref{eqn:BDF2_fd_stability_L2}) and the assumptions on the starting data (\ref{eqn:BDF2_starting_stability_L2}) and (\ref{eqn:BDF_starting_MD_stability}) completes the proof of (\ref{eqn:BDF2_md_stability_fd}). As usual the equivalence of norms yields (\ref{eqn:BDF2_md_stability_fd_l}).
}
\end{Proof}
\begin{The}[Error bound for the fully discrete scheme (\ref{eqn:BDF2_fd_scheme})]
\label{the:BDF2_fd_convergence}
\changes{
Let $u$ be a sufficiently smooth solution of (\ref{eqn:pde}),  let the geometry be sufficiently regular
and let $u_h^i, (i=0,\dots,N)$ denote the lift of the solution of the BDF2 fully discrete scheme (\ref{eqn:BDF2_fd_scheme}). Furthermore, assume that initial data is sufficiently smooth and the initial approximations for the scheme are such that \margnote{ref 2. pt 36.}
\begin{equation}\label{eqn:IC_approx}
\ltwon{u(\cdot,0)-\Ritz u(\cdot,0)}{\G^0}+\ltwon{\Ritz u(\cdot,0)-u_h^0}{\G^0}\leq ch^2,
\end{equation}
and
\begin{equation}
\ltwon{u(\cdot,t^1)-\Ritz u(\cdot,t^1)}{\G(t^1)}+\ltwon{\Ritz u(\cdot,t^1)-u_h^1}{\G(t^1)}\leq c(h^2+\tau^2),
\end{equation}
hold. Furthermore, assume the starting values satisfy the stability assumptions (\ref{eqn:BDF2_starting_stability_L2}) and (\ref{eqn:BDF_starting_MD_stability}).
Then for $0<h\leq h_0,0<\tau\leq \tau_0$, with $h_0$ dependent on the data of the problem and $\tau_0$ dependent on the data of the problem and the arbitrary tangential velocity $\vec a_\tangent$, the following error bound holds. For $n\in\{2,\dots,N\}$ the solution of the fully discrete BDF2 scheme satisfies 
\begin{align}
\ltwon{u(\cdot,t^n)-u_h^n}{\Gtn{n}}^2+&c_1h^2\tau\sum_{i=2}^n\ltwon{\nabla_{\G_h^i}\left(u(\cdot,t^i)-u_h^i\right)}{\Gtn{i}}^2 \\ 
\notag
\leq& c\left(h^4+\tau^4\right)\Bigg(\sup_{s\in[0,T]}\Hiln{u}{2}{\G(s)}^2\\
\notag
&+\int_0^T\Hiln{u}{2}{\Gt}^2+\Hiln{\mdt{\vec v_a}u}{2}{\Gt}^2+\Hiln{\mdt{\vec v_a}(\mdt{\vec v_a}u)}{1}{\Gt}^2\diff t\Bigg).
\end{align}
}
\end{The}
We follow a similar strategy  to that employed in the semidiscrete case to prove the theorem.
We decompose the error as in \S \ref{subsec:semidisc_error_decomp} setting
\begin{equation}\label{eqn:fullydisc_error_decomp}
u(\cdot,t^n)-u^n_h=\rho^n+\theta^n, \quad \rho^n=\rho(\cdot,t^n)=u(\cdot,t^n)-\Ritz u(\cdot,t^n),\quad \theta^n=\Ritz u(\cdot,t^n)-u^n_h\in\Scl,
\end{equation}
with $\Ritz$ the Ritz projection defined in (\ref{eqn:RP_definition}) and $u_h^n$ the lift of the solution to the fully discrete scheme at time $t^n$.

From the scheme (\ref{eqn:BDF2_fd_scheme}) on the interval $[t^{n-1},t^{n+1}]$ we have
\begin{align}\label{eqn:BDF2_fd_err_1} 
\Ltbil{{u}_h^L}{\varphi_h^{n+1}}&
+\abil{u_h^{n+1}}{\varphi_h^{n+1}}+\bbil{u_h^{n+1}}{\varphi_h^{n+1}}{(\vec t^a_h)^{n+1}}\\
\notag
=&\Ltbil{{u}_h^L}{\varphi_h^{n+1}}-\Lthbil{U_h^L}{\Phi_h^{n+1}}
+\abil{u_h^{n+1}}{\varphi_h^{n+1}}-\ahbil{U_h^{n+1}}{\Phi_h^{n+1}}\\
\notag
&+\bbil{u_h^{n+1}}{\varphi_h^{n+1}}{(\vec t^a_h)^{n+1}}-\bhbil{U_h^{n+1}}{\Phi_h^{n+1}}{(\vec T^a_h)^{n+1}}
\\
\notag:=&H_1(\varphi_h^{n+1}).
\end{align}
From the definition of the Ritz projection (\ref{eqn:RP_definition}) we have 
\begin{align}\label{eqn:BDF2_fd_err_2}
&\Ltbil{\Ritz u}{\varphi_h^{n+1}}+\abil{\Ritz u^{n+1}}{\varphi_h^{n+1}}+\bbil{\Ritz u^{n+1}}{\varphi_h^{n+1}}{(\vec t^a_h)^{n+1}}
\\
\notag
&=-\Ltbil{\rho}{\varphi_h^{n+1}}+\Ltbil{u}{\varphi_h^{n+1}}+\abil{u^{n+1}}{\varphi_h^{n+1}}+\bbil{\Ritz u^{n+1}}{\varphi_h^{n+1}}{(\vec t^a_h)^{n+1}}
\\
\notag
&:=H_2(\varphi_h^{n+1}).
\end{align}
Taking the difference of (\ref{eqn:BDF2_fd_err_2}) and (\ref{eqn:BDF2_fd_err_1}) we arrive at the error relation between the fully discrete  solution and the Ritz projection, for $\varphi^{n+1}_h=(\Phi^{n+1})^l\in\Scln{n+1}$
\begin{equation}\label{eqn:BDF2_fd_error_relation}
\Ltbil{{\theta}^L}{\varphi_h^{n+1}}+\abil{\theta^{n+1}}{\varphi_h^{n+1}}+\bbil{\theta^{n+1}}{\varphi_h^{n+1}}{(\vec t^a_h)^{n+1}}
=H_2(\varphi_h^{n+1})-H_1(\varphi_h^{n+1}).
\end{equation}

\begin{Lem}\label{Lem:H1}
For $H_1$ defined in (\ref{eqn:BDF2_fd_err_1}) and for all $\epsilon>0$, we have the estimate
\begin{align}
\lv H_1(\varphi_h^{n+1})\rv\leq& \frac{c(\epsilon)}{\tau}h^4\int_{t^{n-1}}^{t^{n+1}}\Hiln{u_h^L}{1}{\Gt}^2+\ltwon{\mdt{\vec v^a_h}u_h^L}{\Gt}^2\diff t\\
\notag
&+c(\epsilon)h^4\Hiln{u_h^{n+1}}{1}{\Gtn{n+1}}^2+c\ltwon{\varphi_h^{n+1}}{\Gtn{n+1}}^2+\epsilon\ltwon{\nabla_{\Gtn{n+1}}\varphi_h^{n+1}}{\Gtn{n+1}}^2\notag
\end{align}
\end{Lem}
\begin{Proof}
From the definition of $H_1$ (\ref{eqn:BDF2_fd_err_1}) we have
\begin{align}\label{eqn:BDF2_fd_convergence_pf_H1}
\lv H_1(\varphi_h^{n+1})\rv&\leq
\lv\Ltbil{U_h^L}{\varphi_h^{n+1}}-\Lthbil{U_h^L}{\Phi_h^{n+1}}\rv+\lv\abil{u_h^{n+1}}{\varphi_h^{n+1}}-\ahbil{U_h^{n+1}}{\Phi_h^{n+1}}\rv\\
&+\lv\bbil{u_h^{n+1}}{\varphi_h^{n+1}}{(\vec t^a_h)^{n+1}}-\bhbil{U_h^{n+1}}{\Phi_h^{n+1}}{(\vec T^a_h)^{n+1}}\rv\notag\\
&:=I+II+III\notag.
\end{align}
For the first term,  we follow \citep[Proof of Lemma 4.3]{doi:10.1137/110828642}, using the transport  formulas (\ref{eqn:transp_L2_h}) and (\ref{eqn:transp_L2}) together with (\ref{eqn:pert_m}) and (\ref{eqn:pert_g}) we have
\begin{align}\label{eqn:H1_1}
I=\leq&
\frac{c}{\tau}\Bigg\vert\int_{t^{n-1}}^{t^{n+1}}
\mhbil{\mdth{\vec V^a_h}U_h^L(\cdot,t)}{\bhn{\Phi}{n+1}(\cdot,t)}+\ghbil{U_h^L(\cdot,t)}{\bhn{\Phi}{n+1}(\cdot,t)}{\vec V^a_h(\cdot,t)}
\\
\notag
&-\mbil{\mdth{\vec v^a_h}u_h^L(\cdot,t)}{\bhn{\varphi}{n+1}(\cdot,t)}+\gbil{u_h^L(\cdot,t)}{\bhn{\varphi}{n+1}(\cdot,t)}{\vec v^a_h(\cdot,t)}
\diff t\Bigg\vert.\\
\notag
\leq& \frac{ch^2}{\tau}\int_{t^{n-1}}^{t^{n+1}}\left(\ltwon{\mdth{\vec v^a_h}u_h^L}{\Gt}\ltwon{\bhn{\varphi}{n+1}}{\Gt}+
\Hiln{u_h^L}{1}{\Gt}\Hiln{\bhn{\varphi}{n+1}}{1}{\Gt}\right)\diff t\
\\
\notag
\leq&\epsilon\ltwon{\nabla_\Gtn{n+1}\varphi_h^{n+1}}{\Gtn{n+1}}^2+c\ltwon{\varphi_h^{n+1}}{\Gtn{n+1}}^2
\\
\notag
&+\frac{c(\epsilon)}{\tau}h^4\int_{t^{n-1}}^{t^{n+1}}\ltwon{\mdth{\vec v^a_h}u_h^L}{\Gt}^2+\Hiln{u_h^L}{1}{\Gt}^2\diff t
\end{align}
 where $\epsilon$ is a positive constant of our choice. Using (\ref{eqn:pert_a}) we conclude that for all $\epsilon>0$
\begin{align} \label{eqn:H1_2}
II&\leq ch^2\ltwon{\nabla_\Gtn{n+1}\varphi_h^{n+1}}{\Gtn{n+1}}\ltwon{\nabla_\Gtn{n+1}u_h^{n+1}}{\Gtn{n+1}}\\
&\leq c(\epsilon)h^4\ltwon{\nabla_\Gtn{n+1}u_h^{n+1}}{\Gtn{n+1}}^2+\epsilon\ltwon{\nabla_\Gtn{n+1}\varphi_h^{n+1}}{\Gtn{n+1}}^2\notag
\end{align}
Using (\ref{eqn:pert_b}) we have for all $\epsilon>0$
\begin{align}\label{eqn:H1_3}
 III&\leq ch^2\ltwon{u_h^{n+1}}{\Gtn{n+1}}\ltwon{\nabla_\Gtn{n+1}\varphi_h^{n+1}}{\Gtn{n+1}}\\
\notag &\leq c(\epsilon)h^4\ltwon{u_h^{n+1}}{\Gtn{n+1}}^2+\epsilon\ltwon{\nabla_\Gtn{n+1}\varphi_h^{n+1}}{\Gtn{n+1}}^2.
\end{align}
Applying the estimates (\ref{eqn:H1_1})---(\ref{eqn:H1_3}) in (\ref{eqn:BDF2_fd_convergence_pf_H1}) completes the proof of the Lemma.
\end{Proof}
\begin{Lem}\label{Lem:H2}
For $H_2$ defined in (\ref{eqn:BDF2_fd_err_2}) and for all $\epsilon>0$, we have the estimate
\begin{align}
\lv H_2(\varphi_h^{n+1})\rv\leq&
 \frac{c}{\tau}h^4\int_{t^{n-1}}^{t^{n+1}}\Hiln{u}{2}{\Gt}^2+\Hiln{\mdt{\vec v_a}u}{2}{\Gt}^2\diff t\\
&+c\tau^3\int_{t^{n-1}}^{t^{n+1}}\Hiln{u}{2}{\Gt}^2+\Hiln{\mdt{\vec v_a}u}{2}{\Gt}^2+\Hiln{\mdt{\vec v_a}(\mdt{\vec v_a}u)}{1}{\Gt}^2\diff t\notag\\
&+ch^4\Hiln{u}{2}{\Gtn{n+1}}^2+c\ltwon{\varphi_h^{n+1}}{\Gtn{n+1}}^2+\epsilon\ltwon{\nabla_{\Gtn{n+1}}\varphi_h^{n+1}}{\Gtn{n+1}}^2\notag
\end{align}
\end{Lem}
\begin{Proof}
We set
\beq
\sigma(t)=
\begin{cases}
\frac{3}{2\tau}\quad &t\in[t^n,t^{n+1}]\\
-\frac{1}{2\tau}\quad &t\in[t^{n-1},t^{n}].
\end{cases}
\label{eqn:sigma_def}
\eeq
We start by noting that using the transport formula (\ref{eqn:transp_L2}), 
\begin{align}
\lv\Ltbil{\rho}{\varphi_h^{n+1}}\rv&=\lv\int_{t^{n-1}}^{t^{n+1}}\sigma(t)\left(\mbil{\mdth{\vec v^a_h}\rho(\cdot,t)}{\bhn{\varphi}{n+1}(\cdot,t)}+\gbil{\rho(\cdot,t)}{\bhn{\varphi}{n+1}(\cdot,t)}{\vec v^a_h}\diff t\right)\rv\\
&\leq \frac{c}{\tau}\int_{t^{n-1}}^{t^{n+1}}\left(\ltwon{\mdth{\vec v^a_h}\rho(\cdot,t)}{\Gt}+\ltwon{\rho(\cdot,t)}{\Gt}
\right)\ltwon{\bhn{\varphi}{n+1}(\cdot,t)}{\Gt}\diff t\notag.
\end{align}
Young's inequality, (\ref{eqn:pbn_pn_l2}), (\ref{eqn:Ritz_bound}) and (\ref{eqn:MD_Ritz_bound}), yield the estimate
\begin{align}\label{eqn:H2_1}
\lv\Ltbil{\rho}{\varphi_h^{n+1}}\rv&\\
\notag
\leq& \frac{ch^4}{\tau}\int_{t^{n-1}}^{t^{n+1}}\Big(\Hiln{\mdt{\vec v_a}u(\cdot,t)}{2}{\Gt}^2+\Hiln{u(\cdot,t)}{2}{\Gt}^2
\Big)\diff t+\ltwon{{\varphi_h}^{n+1}}{\Gtn{n+1}}^2.
\end{align}
Integrating in time the variational form (\ref{eqn:ale_vf}) over the interval $[t^{n},t^{n+1}]$ with $\varphi=\bhn{\varphi}{n+1}$ we have
\begin{align}
\mbil{u^{n+1}}{\bhn{\varphi}{n+1}(\cdot,t^{n+1})}&-\mbil{u^{n}}{\bhn{\varphi}{n+1}(\cdot,t^{n})}+\int_{t^n}^{t^{n+1}}\abil{u(\cdot,t)}{\bhn{\varphi}{n+1}(\cdot,t)}\diff t\\
\notag
=\int_{t^n}^{t^{n+1}}&-\bbil{u(\cdot,t)}{\bhn{\varphi}{n+1}(\cdot,t)}{\vec a_\tangent}+\mbil{u(\cdot,t)}{\mdt{\vec v_a}\bhn{\varphi}{n+1}(\cdot,t)}\diff t.
\end{align}
Similarly integrating in time the variational form (\ref{eqn:ale_vf}) over the interval $[t^{n-1},t^{n+1}]$ with $\varphi=\bhn{\varphi}{n+1}$ we have
\begin{align}
\mbil{u^{n+1}}{\bhn{\varphi}{n+1}(\cdot,t^{n+1})}&-\mbil{u^{n-1}}{\bhn{\varphi}{n+1}(\cdot,t^{n-1})}+\int_{t^{n-1}}^{t^{n+1}}\abil{u(\cdot,t)}{\bhn{\varphi}{n+1}(\cdot,t)}\diff t\\
\notag
=\int_{t^{n-1}}^{t^{n+1}}&-\bbil{u(\cdot,t)}{\bhn{\varphi}{n+1}(\cdot,t)}{\vec a_\tangent}+\mbil{u(\cdot,t)}{\mdt{\vec v_a}\bhn{\varphi}{n+1}(\cdot,t)}\diff t.
\end{align}
From the definition (\ref{eqn:L2}), we observe that
\begin{align}
\Ltbil{u}{\varphi^{n+1}}=&
\frac{2}{\tau}\left(\mbil{u^{n+1}}{\bhn{\varphi}{n+1}(\cdot,t^{n+1})}-\mbil{u^{n}}{\bhn{\varphi}{n+1}(\cdot,t^{n})}\right)\\
\notag
&-\frac{1}{2\tau}\left(\mbil{u^{n+1}}{\bhn{\varphi}{n+1}(\cdot,t^{n+1})}-\mbil{u^{n-1}}{\bhn{\varphi}{n+1}(\cdot,t^{n-1})}\right)\\
\notag
=&\int_{t^{n-1}}^{t^{n+1}}\sigma(t)\Bigg(\mbil{u(\cdot,t)}{\mdt{\vec v_a}\bhn{\varphi}{n+1}(\cdot,t)}-\abil{u(\cdot,t)}{\bhn{\varphi}{n+1}(\cdot,t)}\\
\notag
&-\bbil{u(\cdot,t)}{\bhn{\varphi}{n+1}(\cdot,t)}{\vec a_\tangent(\cdot,t)}\Bigg)\diff t,
\end{align}
with $\sigma$ as defined in (\ref{eqn:sigma_def}).
Thus we have
\begin{align}\label{eqn:H2_1234}
\Ltbil{u}{\varphi_h^{n+1}}&+\abil{u^{n+1}}{\varphi_h^{n+1}}+\bbil{\Ritz u^{n+1}}{\varphi_h^{n+1}}{(\vec t^a_h)^{n+1}}=\\
&\notag
\left(\bbil{\Ritz u^{n+1}}{\varphi_h^{n+1}}{(\vec t^a_h)^{n+1}}-\bbil{u^{n+1}}{\varphi_h^{n+1}}{\vec a_\tangent^{n+1}}\right)\\
&\notag
+\left(\bbil{u^{n+1}}{\varphi_h^{n+1}}{\vec a_\tangent^{n+1}}-\int_{t^{n-1}}^{t^{n+1}}\sigma(t)\bbil{u(\cdot,t)}{\bhn{\varphi}{n+1}(\cdot,t)}{\vec a_\tangent(\cdot,t)}\diff t\right)\\
&\notag
+\left(\abil{u^{n+1}}{\varphi_h^{n+1}}-\int_{t^{n-1}}^{t^{n+1}}\sigma(t)\abil{u(\cdot,t)}{\bhn{\varphi}{n+1}(\cdot,t)}\diff t\right)\\
&\notag
+\int_{t^{n-1}}^{t^{n+1}}\sigma(t)\mbil{u(\cdot,t)}{\mdt{\vec v_a}\bhn{\varphi}{n+1}(\cdot,t)}\\
\notag
&:= I + II +III +IV
\end{align}
The first term on the right of  (\ref{eqn:H2_1234})
is estimated as follows, we have
\begin{align}\label{eqn:fd_convergence_pf_H2_I_1}
\lv I\rv\leq&\lv-\bbil{\rho^{n+1}}{\varphi_h^{n+1}}{(\vec t^a_h)^{n+1}}\rv+\lv\bbil{u^{n+1}}{\varphi_h^{n+1}}{(\vec t^a_h)^{n+1}-\vec a_\tangent^{n+1}}\rv
\end{align}
For the first term on the right hand side of (\ref{eqn:fd_convergence_pf_H2_I_1}) we use (\ref{eqn:Ritz_bound}) to see that for all $\epsilon>0$
\beq\label{eqn:fd_convergence_pf_H2_I_1_1}
\lv-\bbil{\rho^{n+1}}{\varphi_h^{n+1}}{(\vec t^a_h)^{n+1}}\rv\leq c(\epsilon)h^4\Hiln{u}{2}{\Gtn{n+1}}^2+\epsilon\ltwon{\nabla_{\Gtn{n+1}}\varphi_h^{n+1}}{\Gtn{n+1}}^2.
\eeq
For the next term on the right hand side of (\ref{eqn:fd_convergence_pf_H2_I_1}) we apply (\ref{tang_velocity_bound}) and observe that for all $\epsilon>0$
\beq\label{eqn:fd_convergence_pf_H2_1_1_2}
\lv\bbil{u^{n+1}}{\varphi_h^{n+1}}{(\vec t^a_h)^{n+1}-\vec a_\tangent^{n+1}}\rv\leq c(\epsilon)h^4\ltwon{u}{\Gtn{n+1}}^2+\epsilon\ltwon{\nabla_{\Gtn{n+1}}\varphi_h^{n+1}}{\Gtn{n+1}}^2.
\eeq
Thus we  have 
\beq\label{eqn:fd_convergence_pf_H2_1}
\lv I\rv\leq c(\epsilon)h^4\Hiln{u}{2}{\Gtn{n+1}}^2+\epsilon\ltwon{\nabla_{\Gtn{n+1}}\varphi_h^{n+1}}{\Gtn{n+1}}^2,
\eeq
for all $\epsilon>0$.
For the second term on the right of  (\ref{eqn:H2_1234}) we have
\begin{align}
\lv II\rv\leq&\frac{1}{\tau}\Bigg(\int_{t^n}^{t^{n+1}}(t^{n+1}-t)(t^{n+1}-3t-4t^n)\lv\frac{\diff ^2}{\diff t^2}\bbil{u(\cdot,t)}{\bhn{\varphi}{n+1}(\cdot,t)}{\vec a_\tangent(\cdot,t)}\rv\diff t\\
\notag
&+\int_{t^{n-1}}^{t^{n}}(t-t^{n-1})^2\lv\frac{\diff ^2}{\diff t^2}\bbil{u(\cdot,t)}{\bhn{\varphi}{n+1}(\cdot,t)}{\vec a_\tangent(\cdot,t)}\rv\diff t\Bigg).
\end{align}
The estimate (\ref{eqn:d2_b}) and the fact that $\mdth{\vec v^a_h}\bhn{\varphi}{n+1}=0$ yield
\begin{align}
\lv II\rv\leq\tau&\int_{t^{n-1}}^{t^{n+1}}\Big(\ltwon{u}{\Gt}\\
\notag
&+\ltwon{\mdth{\vec v^a_h}u}{\Gt}+\ltwon{\mdth{\vec v^a_h}(\mdth{\vec v^a_h}u)}{\Gt}\Big)\ltwon{\nabla_\Gt\bhn{\varphi}{n+1}}{\Gt}\diff t.
\end{align}
Young's inequality and (\ref{eqn:pbn_pn_grad}) give for all $\epsilon>0$,
\begin{align}\label{eqn:H2_4}
\lv II\rv\leq& c(\epsilon)\tau^3\int_{t^{n-1}}^{t^{n+1}}\left(\ltwon{u}{\Gt}^2+\ltwon{\mdth{\vec v^a_h}u}{\Gt}^2+\ltwon{\mdth{\vec v^a_h}(\mdth{\vec v^a_h}u)}{\Gt}^2\right)\diff t
\\
\notag
&+\epsilon\ltwon{\nabla_\Gtn{n+1}{\varphi_h}^{n+1}}{\Gtn{n+1}}^2.
\end{align}
The third term on the right of (\ref{eqn:H2_1234}) is estimated in the same way using (\ref{eqn:d2_a}) and (\ref{eqn:pbn_pn_grad}) to give for all $\epsilon>0$,
\begin{align}\label{eqn:H2_5}
\lv III\rv\leq c(\epsilon)\tau^3\int_{t^{n-1}}^{t^{n+1}}&\Big(\ltwon{\nabla_\Gt\mdth{\vec v^a_h}(\mdth{\vec v^a_h}u)}{\Gt}^2+\ltwon{\nabla_\Gt\mdth{\vec v^a_h}u}{\Gt}^2\\
\notag
&
+\ltwon{\nabla_\Gt u}{\Gct}^2\Big)\diff t
+\epsilon\ltwon{\nabla_\Gtn{n+1}{\varphi_h}^{n+1}}{\Gtn{n+1}}^2.
\end{align}
The fourth term on the right of  (\ref{eqn:H2_1234}) may be estimated  using (\ref{md_l2_bound}) together with the fact that $\mdth{\vec v^a_h}\bhn{\varphi}{n+1}=0$ which gives for all $\epsilon>0$,
\beq\label{eqn:H2_6}
\lv IV\rv\leq \frac{c(\epsilon)}{\tau}h^4\int_{t^n}^{t^{n+1}}\ltwon{u}{\Gt}^2\diff t+\epsilon\ltwon{\nabla_\Gtn{n+1}{\varphi_h}^{n+1}}{\Gtn{n+1}}^2.
\eeq
The estimates (\ref{eqn:H2_1}), (\ref{eqn:fd_convergence_pf_H2_1}), (\ref{eqn:H2_4}), (\ref{eqn:H2_5}) and (\ref{eqn:H2_6}) together with the estimates (\ref{mdmd_l2_bound}) and (\ref{mdmd_grad_bound}) completes 
 the proof of the Lemma.
\end{Proof}

We may now finally complete the proof of Theorem \ref{the:BDF2_fd_convergence}.

\begin{Proof}[of Theorem \ref{the:BDF2_fd_convergence}]
\changes{
With the error decomposition of (\ref{eqn:fullydisc_error_decomp})  and the estimates on the Ritz projection error \ref{eqn:Ritz_bound} it remains to bound $\theta$. With the same argument as used in the proof of Lemma \ref{Lem:BDF2_fd_stability}, i.e., (\ref{eqn:BDF_l2_stab_pf_1})---(\ref{eqn:BDF_l2_stab_pf_5}) and the usual estimation of the non-symmetric term using Young's inequality, we have
 \margnote{ref 2. pt 37.}
\begin{align}
\frac{1}{4}\pdtau&\Big(\mbil{\theta^n}{\theta^n}+\mbil{2\bn{\theta}{n}(\cdot,t^{n-1})-\theta^{n-1}}{2\bn{\theta}{n}(\cdot,t^{n-1})-\theta^{n-1}}\Big)\\
\notag
+&(1-\epsilon)\ltwon{\nabla_{\Gtn{n+1}}\theta^{n+1}}{\Gtn{n+1}}^2\leq \lv H_1(\theta^{n+1})\rv +\lv H_2(\theta^{n+1})\rv +c\Big(c(\epsilon)\ltwon{\theta^{n+1}}{\Gtn{n+1}}^2\\
\notag
&+\ltwon{\theta^{n}}{\Gtn{n}}^2+\mbil{2\bn{\theta}{n}(\cdot,t^{n-1})-\theta^{n-1}}{2\bn{\theta}{n}(\cdot,t^{n-1})-\theta^{n-1}}\Big),
\end{align}
for $\epsilon$ a positive constant of our choice.
Inserting the bounds from Lemmas \ref{Lem:H1} and \ref{Lem:H2} we obtain
\begin{align}
\frac{1}{4}\pdtau&\Big(\mbil{\theta^n}{\theta^n}+\mbil{2\bn{\theta}{n}(\cdot,t^{n-1})-\theta^{n-1}}{2\bn{\theta}{n}(\cdot,t^{n-1})-\theta^{n-1}}\Big)\\
\notag
+&(1-\epsilon)\ltwon{\nabla_{\Gtn{n+1}}\theta^{n+1}}{\Gtn{n+1}}^2 \leq c\Big(\ltwon{\theta^{n+1}}{\Gtn{n+1}}^2+\ltwon{\theta^{n}}{\Gtn{n}}^2\\
\notag
&+\mbil{2\bn{\theta}{n}(\cdot,t^{n-1})-\theta^{n-1}}{2\bn{\theta}{n}(\cdot,t^{n-1})-\theta^{n-1}}\Big)
\\ 
&
\notag
+\frac{c}{\tau}h^4\int_{t^{n-1}}^{t^{n+1}}\Hiln{u}{2}{\Gt}^2+\Hiln{\mdt{\vec v_a}u}{2}{\Gt}^2+\Hiln{u_h^L}{1}{\Gt}^2+\ltwon{\mdt{\vec v^a_h}u_h^L}{\Gt}^2\diff t
\\
&
\notag
+c\tau^3\int_{t^{n-1}}^{t^{n+1}}\Hiln{u}{2}{\Gt}^2+\Hiln{\mdt{\vec v_a}u}{2}{\Gt}^2+\Hiln{\mdt{\vec v_a}(\mdt{\vec v_a}u)}{1}{\Gt}^2\diff t
\\
&
\notag
+ch^4\Hiln{u}{2}{\Gtn{n+1}}^2+ch^4\Hiln{u_h^{n+1}}{1}{\Gtn{n+1}}^2,
\end{align}
where we have suppressed the dependence of the constants on $\epsilon$.
Summing over time, multiplying by $\tau$ and choosing $\epsilon>0$ suitably yields (where we have dropped a positive term), for $n\in\{2,\dots,N\}$
\begin{align}
\ltwon{\theta^n}{\Gtn{n}}^2&+c_1\tau\sum_{k=2}^n\ltwon{\nabla_\Gtn{k}\theta^{k}}{\Gtn{k}}^2
\leq \ltwon{\theta^1}{\Gtn{1}}^2+c\tau\sum_{i=1}^n\ltwon{\theta^i}{\Gtn{i}}^2\\
\notag
&
+\mbil{2\bn{\theta}{1}(\cdot,t^{0})-\theta^{0}}{2\bn{\theta}{1}(\cdot,t^{0})-\theta^{0}}
\\
\notag
&+ch^4\int_{0}^{t^{n}}\Hiln{u}{2}{\Gt}^2+\Hiln{\mdt{\vec v_a}u}{2}{\Gt}^2+\Hiln{u_h^L}{1}{\Gt}^2+\ltwon{\mdt{\vec v^a_h}u_h^L}{\Gt}^2\diff t
\\
\notag
&+c\tau^4\int_{0}^{t^{n}}\Hiln{u}{2}{\Gt}^2+\Hiln{\mdt{\vec v_a}u}{2}{\Gt}^2+\Hiln{\mdt{\vec v_a}(\mdt{\vec v_a}u)}{1}{\Gt}^2\diff t
\\
\notag
&
+c\tau h^4\sum_{i=2}^n\left(\Hiln{u}{2}{\Gtn{i}}^2+ch^4\Hiln{u_h^{i}}{1}{\Gtn{i}}^2\right).
\end{align}
A discrete Gronwall \margnote{ref 2. pt. 38.} argument together with the stability bounds of Lemmas \ref{Lem:sd_stability} \ref{Lem:BDF2_fd_stability} and the assumptions on the approximation of the initial data and starting values complete the proof.
}
\end{Proof}

\subsection{Fully discrete BDF1 ALE-ESFEM scheme}
\changes{\margnote{ref 2. pt. 39.}
We could also have considered an implicit Euler time discretisation of the semidiscrete scheme (\ref{eqn:sd_scheme}) as follows.
Given $U^0_h\in\Scn{0}$ find $U_h^{n+1}\in\Scn{n+1},n\in\{0,\dots,N-1\rbrace$ such that for all $\Phi_h^{n+1}\in\Scn{n+1}$
and for $n\in\lbrace0,\dots,N-1\rbrace$
\beq\label{eqn:fd_scheme}
\frac{1}{\tau}\left(\mhbil{U_h^{n+1}}{{\Phi_h}^{n+1}}-\mhbil{U_h^{n}}{\bhn{\Phi}{n+1}(\cdot,t^{n})}\right)+\ahbil{U^{n+1}_h}{\Phi^{n+1}_h}=-\bhbil{U_h^{n+1}}{\Phi_h^{n+1}}{(\vec T^a_h)^{n+1}}.
\eeq
}
Using the ideas in the analysis presented above  it is a relatively straight forward extension of \citep{doi:10.1137/110828642}
to show the following error bound.
\begin{Cor}[Error bound for an implicit Euler time discretisation]\label{Cor:fd_convergence}
Let $u$ be a sufficiently smooth solution of (\ref{eqn:pde}) and  let the geometry be sufficiently regular.
Furthermore let $u_h^i, (i=0,\dots,N)$ denote the lift of the solution of the implicit Euler fully discrete scheme (\ref{eqn:fd_scheme}). Furthermore, assume that initial data is sufficiently smooth and that the approximation of the initial data is such that 
\begin{equation}
\ltwon{u(\cdot,0)-\Ritz u(\cdot,0)}{\G^0}+\ltwon{\Ritz u(\cdot,0)-u_h^0}{\G^0}\leq ch^2,
\end{equation}
holds.
Then for $0<h\leq h_0,0<\tau\leq \tau_0$ (with $h_0$ dependent on the data of the problem and $\tau_0$ dependent on the data of the problem and the arbitrary tangential velocity $\vec a_\tangent$) the following error bound holds.\changes{ For $n\in\{0,\dots,N\}$\margnote{ref 2. pt 40.}
\begin{align}
\ltwon{u(\cdot,t^n)-u_h^n}{\Gtn{n}}^2+c_1h^2\tau\sum_{i=1}^n\ltwon{\nabla_{\G^i_h}\left(u(\cdot,t^i)-u_h^i\right)}{\Gtn{i}}^2 \\ \leq c\left(h^4+\tau^2\right)\sup_{t\in[0,T]}\left(\Hiln{u}{2}{\Gt}^2+\Hiln{\mdt{\vec v_a}u}{2}{\Gt}^2\right).\notag
\end{align}
}
\end{Cor}
\section{Numerical experiments}\label{Sec:examples}
We report on numerical simulations that support our theoretical results and illustrate that, for certain material velocities, the arbitrary tangential velocity may be chosen such that the meshes generated during the evolution are more suitable for computation than in the Lagrangian case.  We also report on an experiment in which we investigate numerically the long time behaviour of solutions to (\ref{eqn:pde}) with different initial data when the evolution of the surface is a periodic function of time.  The  code for the simulations made use of  the finite element library ALBERTA 
\cite{schmidt2005design} and for the visualisation we used PARAVIEW \cite{henderson2004paraview}.\changes{In many of the examples the velocity fields and the suitable right hand sides (in the case of benchmark examples) were computed using Maple\textsuperscript{TM}.

For each of the simulations, an initial triangulation $\G_h^0$ is obtained by first defining a coarse macro triangulation that interpolates at the vertices the continuous surface and subsequently  refining and projecting the new nodes onto the continuous surface. The vertices  are then advected with the velocity $\vec v_a$ (c.f. \S 1). In practice it is often the case that this velocity must be determined by solving an ODE, throughout the above analysis we have assumed this ODE is solved exactly and hence that the vertices lie on the continuous surface at all times. 
\margnote{ref 2. pt. 43.}
}

\begin{Example}[Benchmarking experiments]\label{eg:benchmark}
We define the level set function
\beq\label{eqn:benchmark_LS}
d(\vec x,t)=\frac{x_1^2}{a(t)}+x_2^2+x_3^2-1,
\eeq
and consider the surface
\beq\label{eqn:benchmark_surface}
\Gt=\left\{\vec x\in\Reals^3\lv d(\vec x,t)=0, x_3\geq 0\right.\right\}.
\eeq
The surface is  the surface of a  hemiellipsoid with time dependent axis. We set $a(t)=1+0.25\sin(t)$ and
we assume that the material velocity of the surface $\vec v$ has zero tangential component. Therefore the material velocity of the surface is given by \citep{dziuk2007finite}
\beq\label{eqn:benchmark_vel}
\vec v=\frac{-\pdt d}{\lv \nabla d\rv}\frac{\nabla d}{\lv \nabla d\rv}=v\normal,
\eeq
\changes{
with
\margnote{ref 2. pt 44.}
\beq\label{eqn:benchmark_vel_explicit}
v(\vec x,t)=\frac{-\pdt d(\vec x,t)}{\lv \nabla d(\vec x,t)\rv}\quad\text{ and }\quad\normal(\vec x,t)=\frac{\nabla d(\vec x,t)}{\lv \nabla d(\vec x,t)\rv}\quad\text{ for $x\in\Gt$, $t\in[0,T]$},
\eeq
where $d$ is given by (\ref{eqn:benchmark_LS}).
}

We consider a time interval $[0,2]$ and insert a suitable right hand side in (\ref{eqn:pde}) such that the exact solution is $u(\vec x,t)=\sin(t)x_1x_2$,  i.e.,\changes{ we compute a right hand side for (\ref{eqn:pde}) from the equation
\beq\label{eqn:RHS}
f=\pdt u +\vec v\nabla u+u\nabla_\Gt \cdot\vec v-\lap_\Gt u.
\eeq
}

To investigate the performance of the proposed BDF2-ALE  ESFEM scheme we report on two numerical experiments. First we consider the Lagrangian scheme i.e., $\vec a_\tangent=\vec 0$. Secondly we consider an evolution in which the arbitrary tangential velocity is nonzero. The velocity is defined as follows;
\beq\label{eqn:benchmark_ALE_velocity}
v^a_1(\vec x,t)=\frac{0.25\cos(t)}{2(1+0.25\sin(t))^{1/2}}x_0,\quad v^a_2(\vec x,t)=v^a_3(\vec x,t)=0,\quad\vec x_0\in\G^0.
\eeq
The arbitrary tangential velocity is then determined by $\vec a_\tangent=\vec v^a-\vec v$ where $\vec v^a$ and $\vec v$ are defined by (\ref{eqn:benchmark_ALE_velocity}) and (\ref{eqn:benchmark_vel}) respectively. We note that 
$\vec v^a\cdot \conormal =0$ as the conormal to the boundary of $\Gt$ is given by $(0,0,-1)^T$.

\changes{
We remark that for this example, the continuous surface and the choice of the arbitrary velocity $\vec v^a$ are such that the lift (c.f., \eqref{eqn:lift}) of the  triangulated surface (with straight boundary faces) is the continuous surface in both the Lagrangian and the ALE case. This holds as  the normal to the continuous surface $\normal(\vec x,t)$ is a vector in the plane $x_3=0$ and the boundary curves $\partial\Gt$ and $\partial\Gct$ (in both the Lagrangian and ALE case) are curves in the plane $x_3=0$. Thus the assumptions of Remark \ref{rem:bdry} are satisfied and the preceding analysis is applicable.
\margnote{ref 2. pt. 44.}
}
\end{Example}

\begin{Defn}{Experimental order of convergence (EOC)}
For a series of triangulations $\left\{\mathcal{T}_i\right\}_{i=0,\dots,N}$ we denote by $\{e_i\}_{i=0,\dots,N}$ the error and by
$h_i$ the mesh size of $\mathcal{T}_i$. The  EOC is given by
\beq\label{eqn:EOC_def}
EOC(e_{i,i+1},h_{i,i+1})=\ln(e_{i+1}/e_i)/\ln(h_{i+1}/h_i).
\eeq
\end{Defn}
\changes{
In Tables \ref{tab:Lag_EOC} and  \ref{tab:ALE_EOC}  we report on the mesh size at the final time together with the errors and EOCs in equivalent norms to the norms appearing in Theorem \ref{the:BDF2_fd_convergence} for the two numerical simulations considered in Example \ref{eg:benchmark}. Specifically we lift the continuous solution onto the discrete surface (the inverse of the lift defined in \eqref{eqn:lift}) and measure the errors in the following norm and seminorm
\begin{align*}
\Lp{\infty}(\Lp{2})&:= \sup_{n\in[2,\dotsc,N]}\ltwon{u(\cdot,t^n)^{-l}-U_h^n}{\Gctn{n}}\\
\Lp{2}\left(\Hil{1}\right)&:= \sum_{i=2}^N\left(\tau\ltwon{\nabla_{\Gctn{i}}\left(u(\cdot,t^i)^{-l}-U_h^i\right)}{\Gctn{i}}^2\right)^{1/2}
\end{align*}
The EOCs were computed using the mesh size at the final time and the timestep was coupled to the initial mesh size. The starting values for the scheme were taken to be the interpolant of the exact solution. We observe that the EOCs support the error bounds of Theorem \ref{the:BDF2_fd_convergence} and that for this example the errors with the Lagrangian and ALE schemes are similar in magnitude.
\margnote{ref 2. pt. 45. (Precise definition of the errors)}
We remark that in all the computations the integrals have been evaluated using numerical quadrature of a sufficiently high order such that the effects of quadrature are negligible in the evaluation of the convergence rates.
}
\begin{table}[h!]
\centering
\begin{tabular}{ccccc}
\toprule
$h$&$\Lp{\infty}(\Lp{2})$&$EOC$&$\Lp{2}\left(\Hil{1}\right)$&$EOC$\\
\midrule
0.88146&   0.07772&   -   &0.63634 &  -\\
   0.47668   &0.02087 &  2.13842   &0.36133  & 0.92064\\
   0.24445  & 0.00546  & 2.00845   &0.18755  & 0.98184\\
   0.12307   &0.00140   &1.97958  & 0.09480   &0.99420\\
   0.06165  & 0.00036   &1.96828 &  0.04754   &0.99823\\
   \bottomrule
\end{tabular}\\
\caption[]{Errors and EOC in the $\Lp{\infty}{\left(0,T;\Lp{2}\right)}$  seminorm  and the $\Lp{2}{\left(0,T;\Hil{1}\right)}$ norm
for Example \ref{eg:benchmark} with the Lagrangian scheme ($\vec a_\tangent=\vec 0$).}\label{tab:Lag_EOC}
\end{table}
\begin{table}[h!]
\centering
\begin{tabular}{ccccc}
\toprule
$h$&$\Lp{\infty}(\Lp{2})$&$EOC$&$\Lp{2}\left(\Hil{1}\right)$&$EOC$\\
\midrule
   0.85679 &  0.07876   &-&   0.63090&   -\\
   0.44695   &0.02134  & 2.00683 &  0.35151&   0.89884\\
   0.22693   &0.00560  & 1.97379  & 0.18173&   0.97332\\
   0.11415   &0.00143   &1.98248   &0.09177&   0.99437\\
   0.05722   &0.00036&   1.98228   &0.04601&   0.99973\\
      \bottomrule
\end{tabular}\\
\caption[]{Errors and EOC in the $\Lp{\infty}{\left(0,T;\Lp{2}\right)}$  norm  and the $\Lp{2}{\left(0,T;\Hil{1}\right)}$ seminorm
for Example \ref{eg:benchmark} with the velocity defined by (\ref{eqn:benchmark_ALE_velocity}) which includes a nonzero arbitrary tangential component $\vec a_\tangent$ .}\label{tab:ALE_EOC}
\end{table}

\begin{Example}[Comparison of the Lagrangian and ALE schemes]\label{eg:comparison}
We define the level set function 
\beq\label{eqn:complex_LS}
d(\vec x,t)=\frac{x_1^2}{a(t)^2}+G(x_2^2)+G\left(\frac{x_3^2}{L(t)^2}\right)-1,
\eeq
\changes{
\margnote{ref 2. pt. 46} where $a(t)=0.1+0.01\sin(2\pi t)$, $L(t)=1+0.3\sin(4\pi t)$ and
$G(s)=31.25s(s-0.36)(s-0.95)$.}
We consider the surface
\beq\label{eqn:complex_surface}
\Gt=\left\{\vec x\in\Reals^3\lv d(\vec x,t)=0\right.\right\}.
\eeq
 To compare the Lagrangian and the ALE numerical schemes we first consider a numerical scheme where the nodes are moved with the material velocity, which we  assume is the normal velocity.  For this Lagrangian scheme we approximate the nodal velocity by solving the ODE (\ref{eqn:benchmark_vel}) at each node numerically with $d$ as in (\ref{eqn:complex_LS}). Secondly we consider an evolution of the form proposed in \cite{EllSty12} where the arbitrary tangential velocity is nonzero. The evolution is defined as follows; for each node $(X_j(t),Y_j(t),Z_j(t))^T:=\vec X_j,j=1,\dotsc,J$,  given nodes
 $\vec X_j(0),j=1,\dotsc,J$ on $\G^0$, we set
\beq\label{eqn:complex_ALE_velocity}
X_j(t)=X_j(0)\frac{a(t)}{a(0)},\ Y_j(t)=Y_j(0)\text{ and }Z_j(t)=Z_j(0)\frac{L(t)}{L(0)},\quad t\in[0,T].
\eeq
Thus $d(\vec X_j(t),t)=0,j=1,\dotsc,J,t\in[0,T].$ \changes{In this case at a vertex $\vec X_j,j=1,\dotsc,J$, the arbitrary tangential velocity $\vec a_\tangent(\vec X_j,t)$ is  given by $\vec a_\tangent(\vec X_j,t)=\frac{\diff}{\diff t}\vec X_j(t)-\vec v(\vec X_j(t),t)$. We note that the value of $\vec a_\tangent$ at the vertices is sufficient to define the tangential velocity that enters the scheme $\vec T^a_h$ (c.f., \eqref{eqn:disc_surface_tang_velocity}).\margnote{ref 2. pt. 46}

 We insert a suitable right hand side for (\ref{eqn:pde}) by computing an $f$ (as in Example \ref{eg:benchmark}) such that the exact solution is $u(\vec x,t)=\cos(\pi t)x_1x_2x_3$ and consider a time interval $[0,1]$.\margnote{ref 2. pt. 46}
}

We used  CGAL \cite{cgal:ry-smg-13b} to generate an initial triangulation $\Gctn{0}$ of $\G^0$. The mesh had 15991 vertices  (the righthand mesh at $t=1$ in Figure \ref{fig:complex_mesh} is identical to the initial mesh). 
 We used the same initial triangulation for both schemes. We considered a time interval corresponding  to a single period of the evolution, i.e., $[0,1]$ and selected a timestep of $10^{-3}$ and used a BDF2 time discretisation, i.e., the scheme (\ref{eqn:BDF2_fd_scheme}). The starting values for the scheme were taken to be the interpolant of the exact solution.

Figure \ref{fig:complex_mesh} shows snapshots of the meshes obtained with the two different velocities. We clearly observe that moving  the vertices of the mesh with the velocity with a nonzero tangential component generates meshes that appear much more suitable for computation than the meshes obtained when the vertices are moved with the material velocity. Figure \ref{fig:complex_errors} shows the interpolant of the error, i.e.,
the Figure shows plots of the function $e_h^n\in\Sc(t^n)$ with nodal values  given by $e_h^n(\vec X_j)=\lv (U_h(\vec X_j))^n-u(\vec X_j,t^n)\rv, j=1,\dotsc,J$. We observe that the ALE scheme has a significantly smaller error than the Lagrangian scheme.
 \begin{figure}[ht]
 \centering 
    \subfigure[][{$t=0.2$}]{
  \includegraphics[ trim = 0mm 0mm 0mm 0mm,  clip,    width=0.33\linewidth
 ]{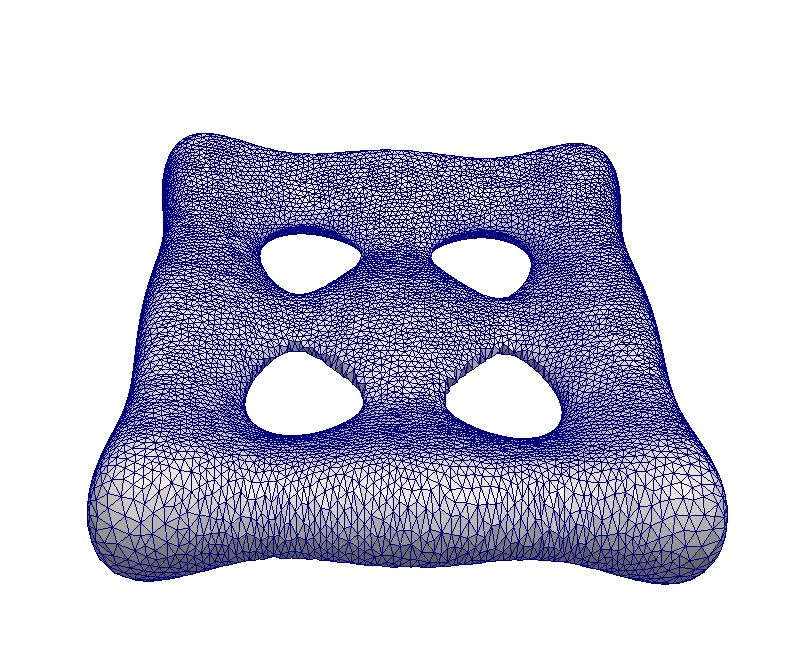}
 \includegraphics[ trim = 0mm 0mm 0mm 0mm,  clip,    width=0.33\linewidth
 ]{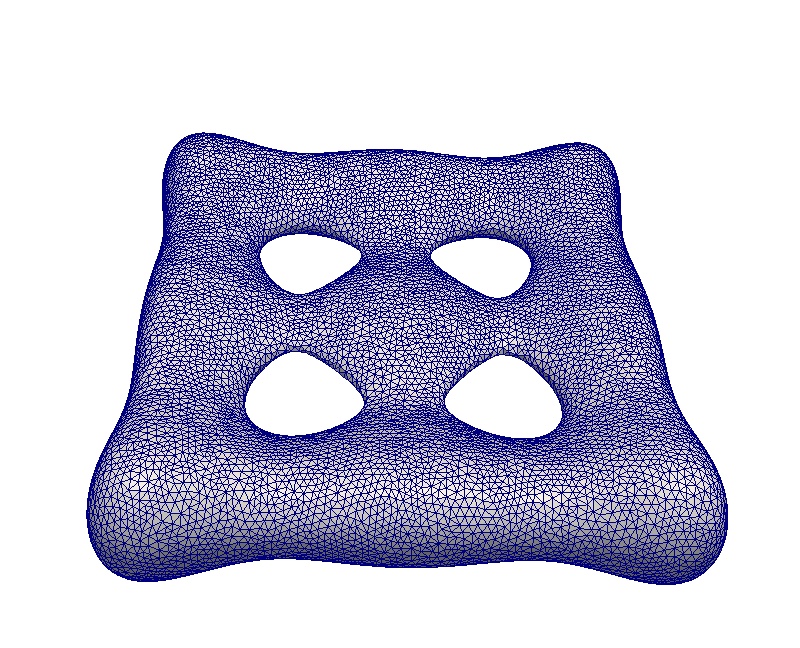}
}
    \subfigure[][{$t=0.4$}]{
  \includegraphics[ trim = 0mm 0mm 0mm 0mm,  clip,    width=0.33\linewidth
 ]{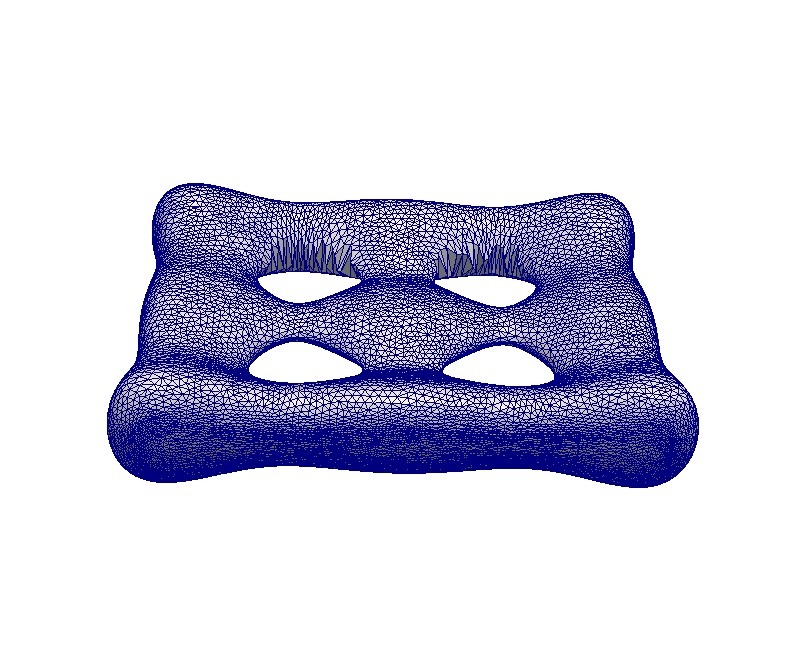}
 \includegraphics[ trim = 0mm 0mm 0mm 0mm,  clip,    width=0.33\linewidth
 ]{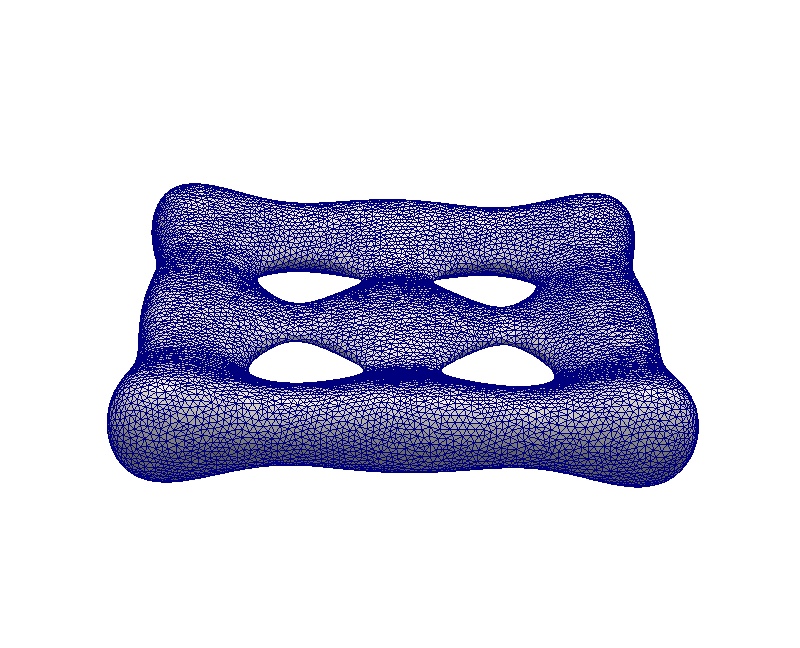}
 }
    \subfigure[][{$t=0.7$}]{
  \includegraphics[ trim = 0mm 0mm 0mm 0mm,  clip,    width=0.33\linewidth
 ]{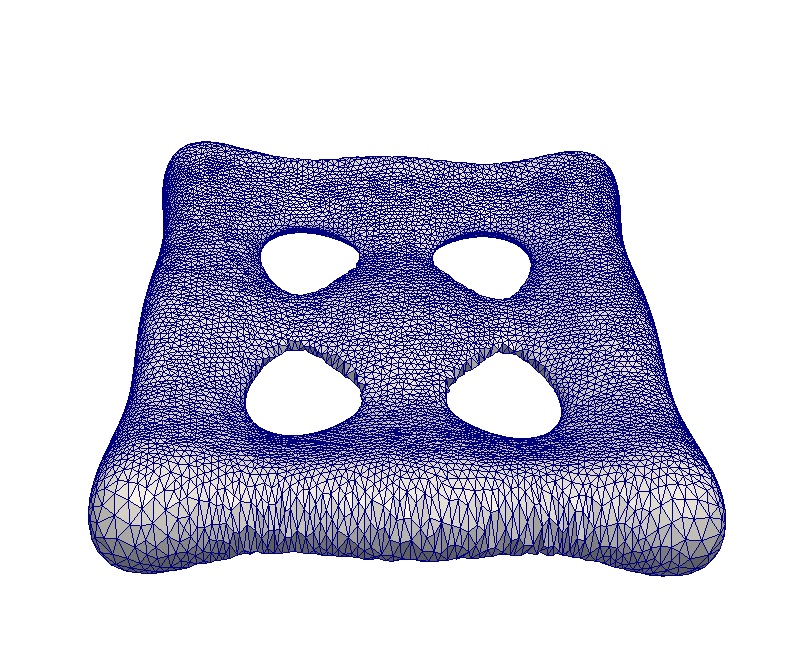}
 \includegraphics[ trim = 0mm 0mm 0mm 0mm,  clip,    width=0.33\linewidth
 ]{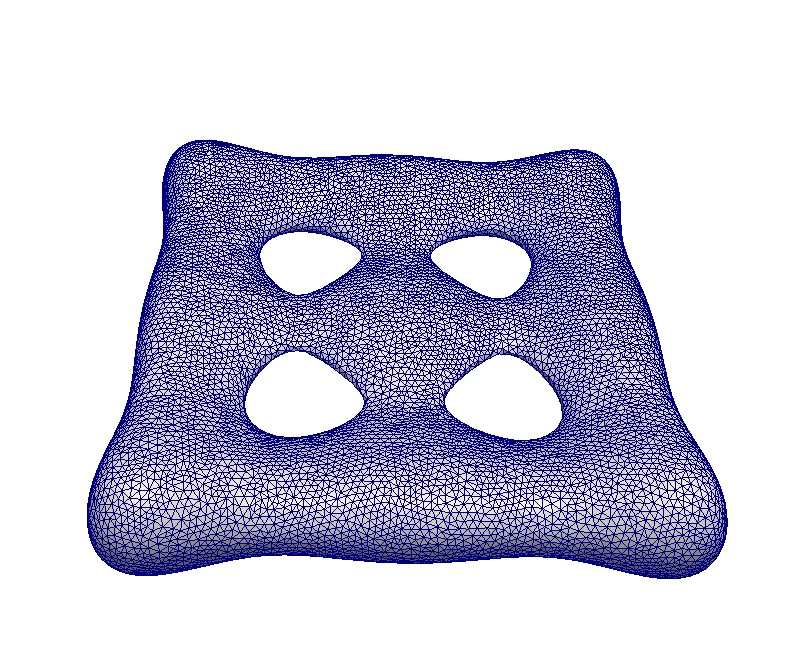}
 }
   \subfigure[][{$t=1$}]{
  \includegraphics[ trim = 0mm 0mm 0mm 0mm,  clip,    width=0.33\linewidth
 ]{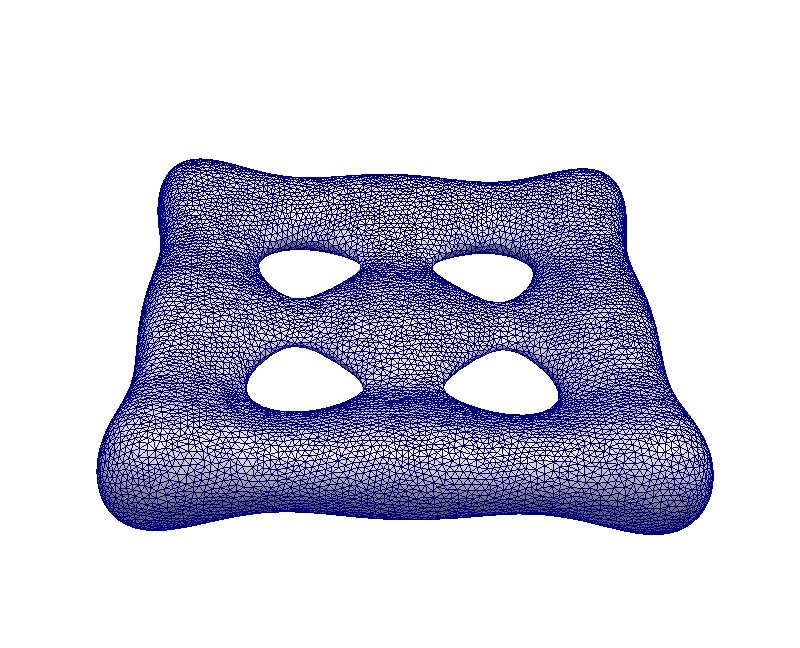}
 \includegraphics[ trim = 0mm 0mm 0mm 0mm,  clip,    width=0.33\linewidth
 ]{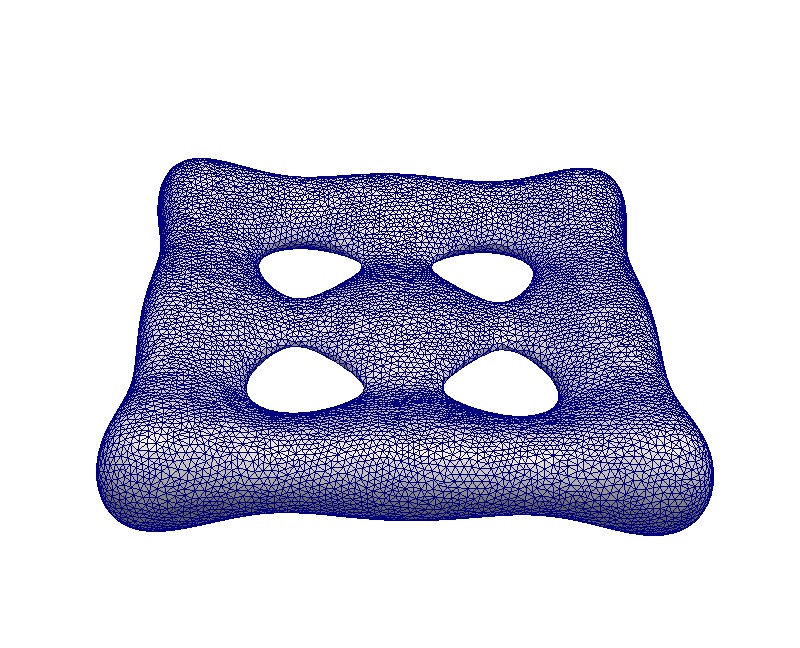}
 }
 \caption{Meshes obtained for Example \ref{eg:comparison} with an approximation of  the  Lagrangian (zero tangential) velocity   (lefthand column) and with the ALE velocity (\ref{eqn:complex_ALE_velocity}) (righthand column).}\label{fig:complex_mesh}
 \end{figure}
  \begin{figure}[ht]
 \centering 
    \subfigure[][{$t=0.2$}]{
  \includegraphics[ trim = 20mm 140mm 20mm 20mm,  clip,    width=0.33\linewidth
 ]{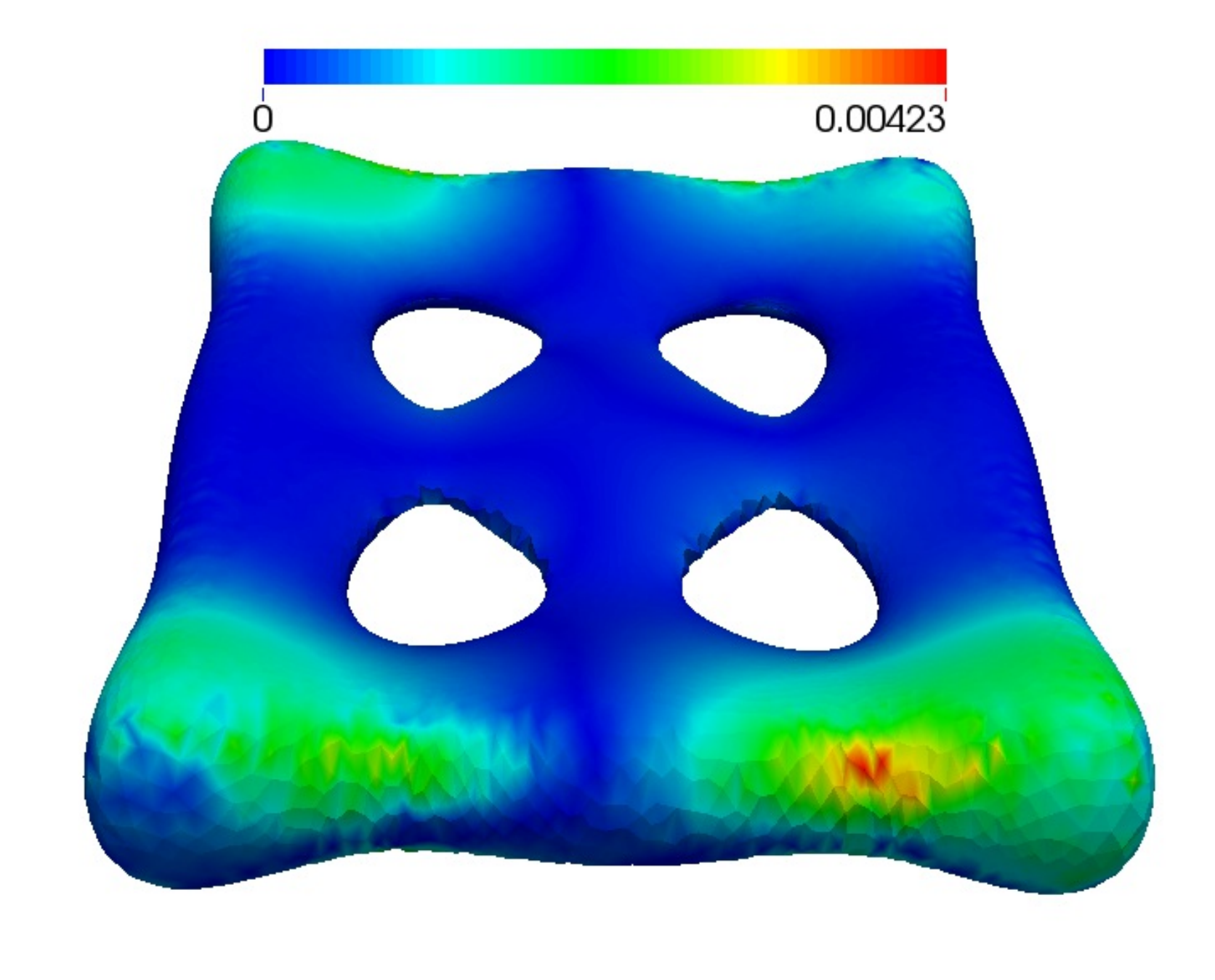}
 \includegraphics[ trim = 20mm 140mm 20mm 20mm,  clip,    width=0.33\linewidth
 ]{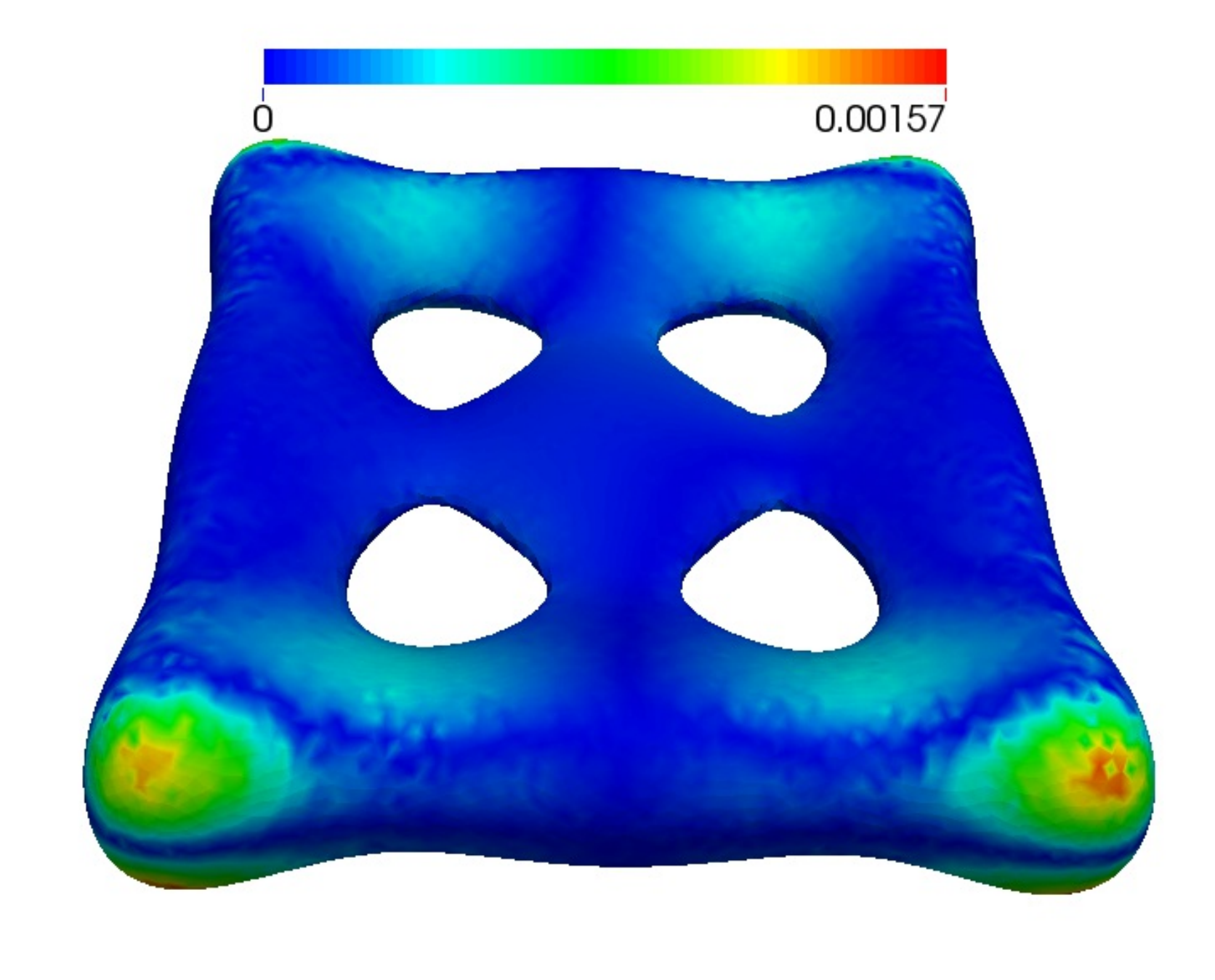}
}
    \subfigure[][{$t=0.4$}]{
  \includegraphics[ trim = 20mm 140mm 20mm 20mm,  clip,    width=0.33\linewidth
 ]{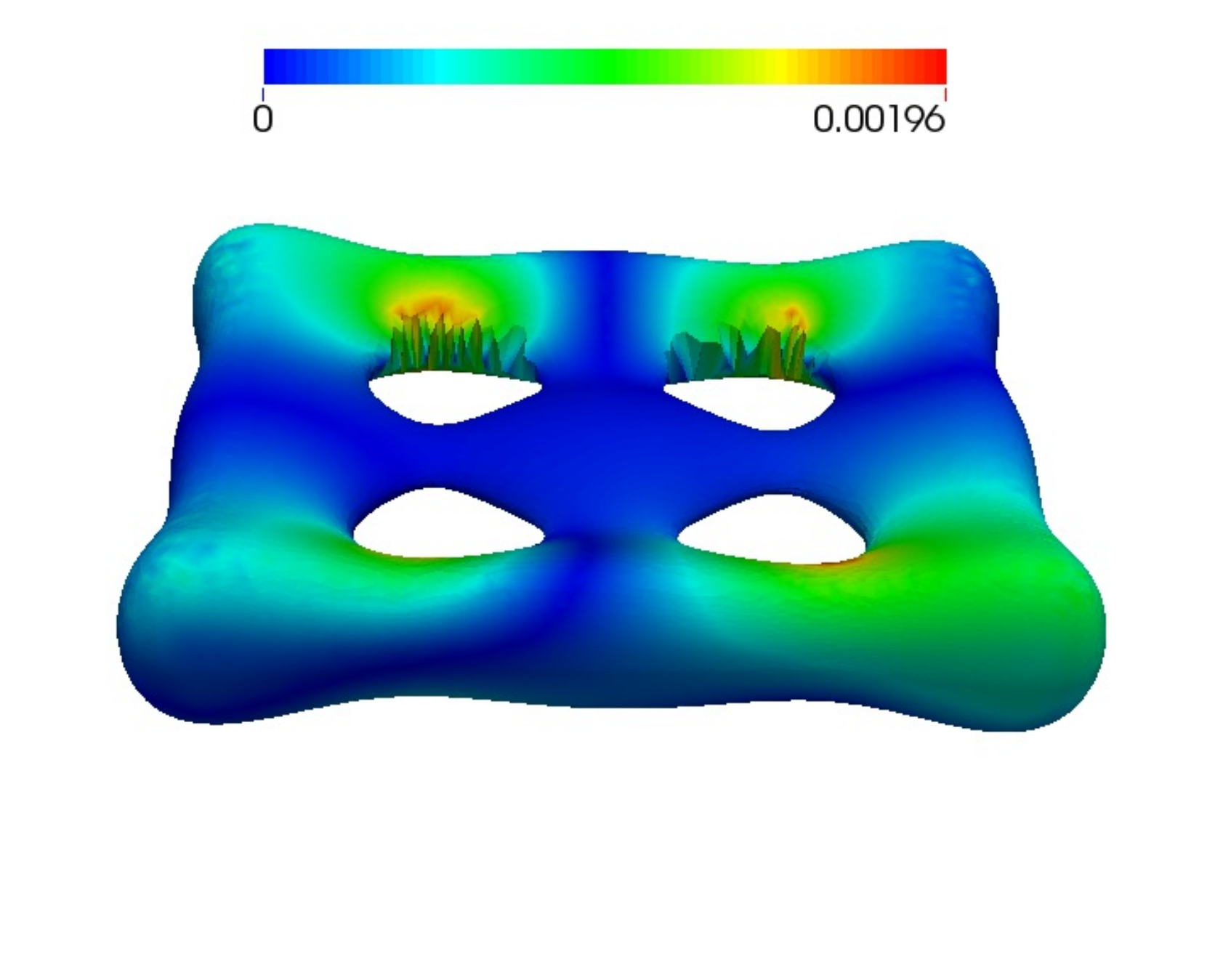}
 \includegraphics[ trim = 20mm 140mm 20mm 20mm,  clip,    width=0.33\linewidth
 ]{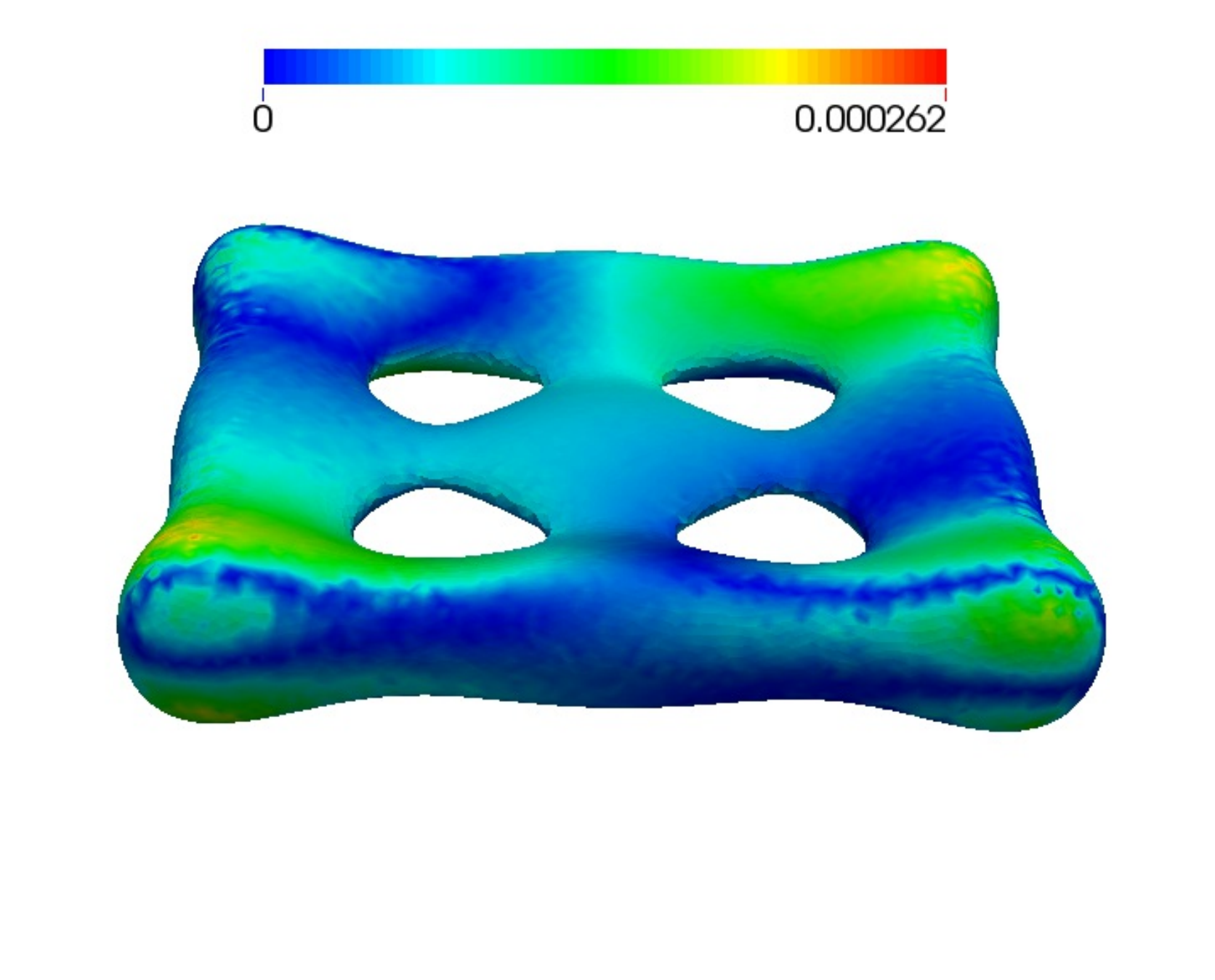}
 }
    \subfigure[][{$t=0.7$}]{
  \includegraphics[ trim = 20mm 140mm 20mm 20mm,  clip,    width=0.33\linewidth
 ]{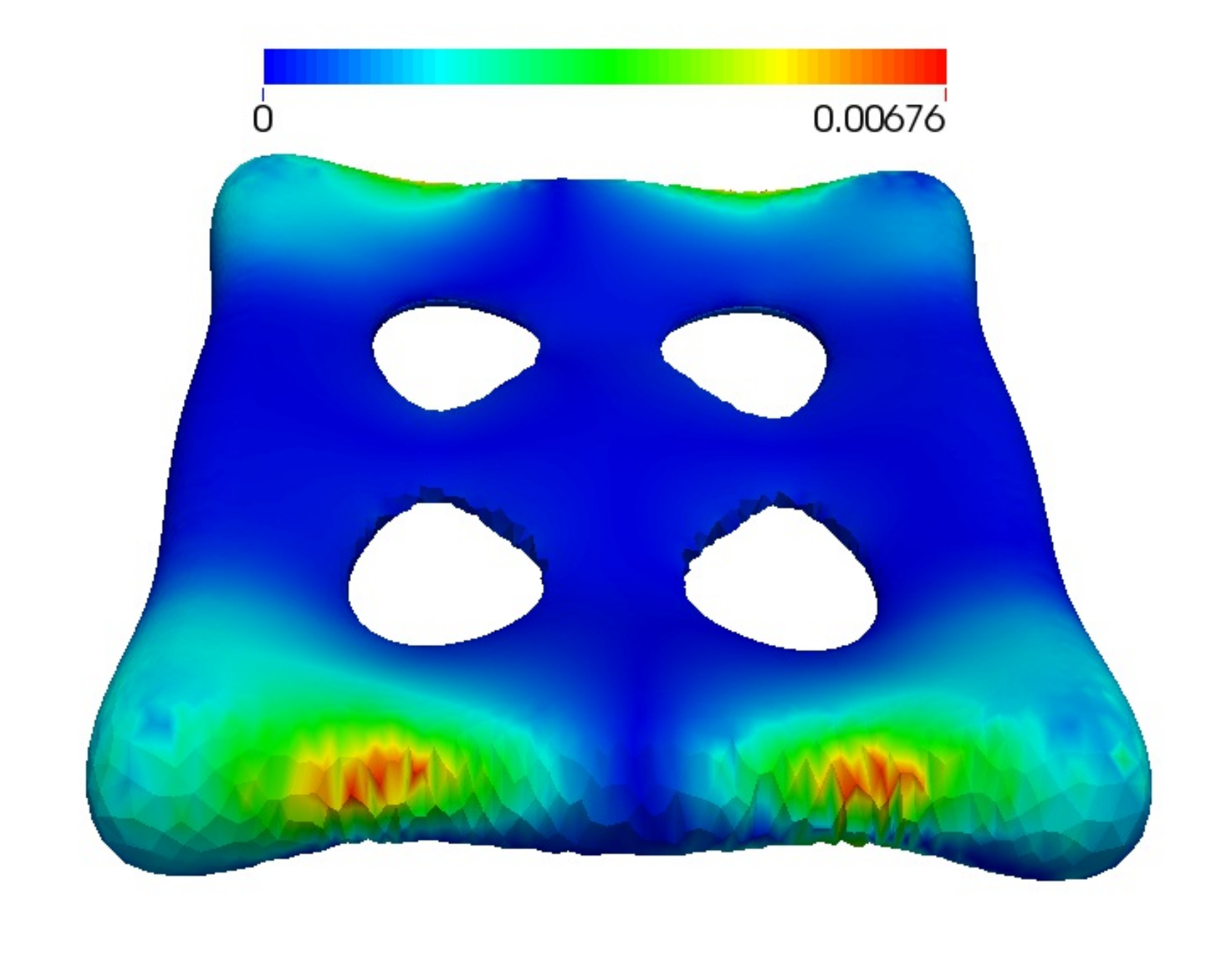}
 \includegraphics[ trim = 20mm 140mm 20mm 20mm,  clip,    width=0.33\linewidth
 ]{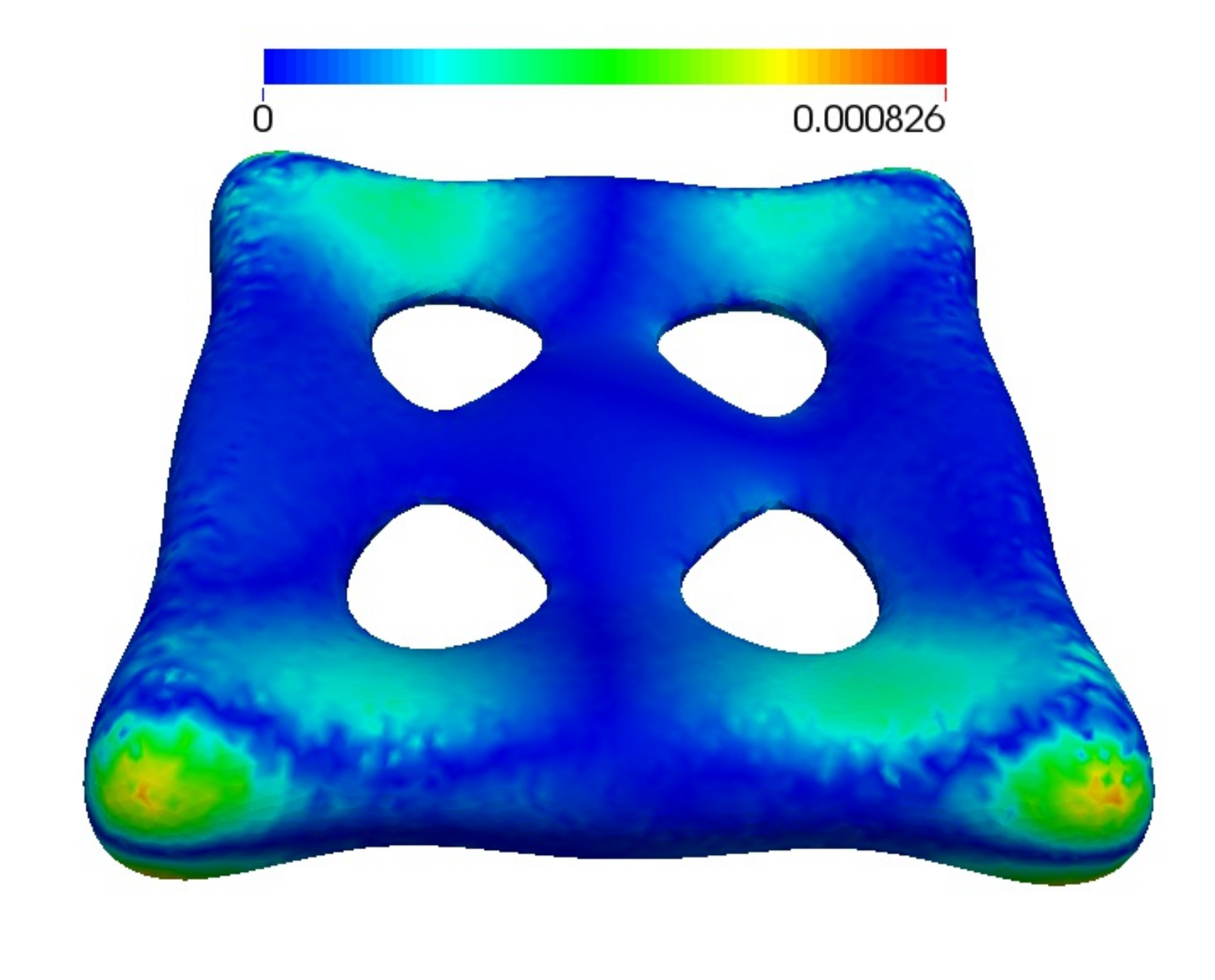}
 }
   \subfigure[][{$t=1$}]{
  \includegraphics[ trim = 20mm 140mm 20mm 20mm,  clip,    width=0.33\linewidth
 ]{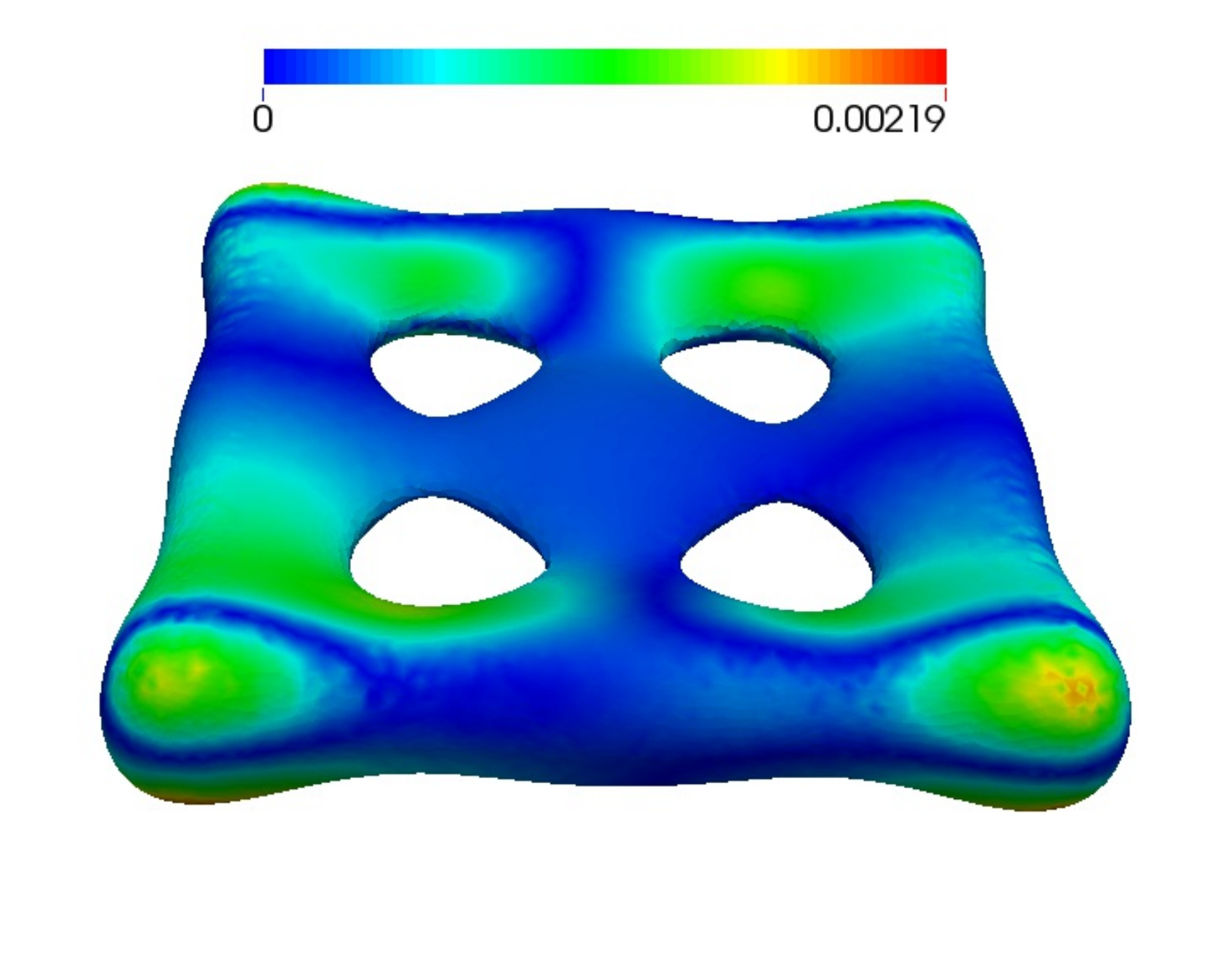}
 \includegraphics[ trim = 20mm 140mm 20mm 20mm,  clip,    width=0.33\linewidth
 ]{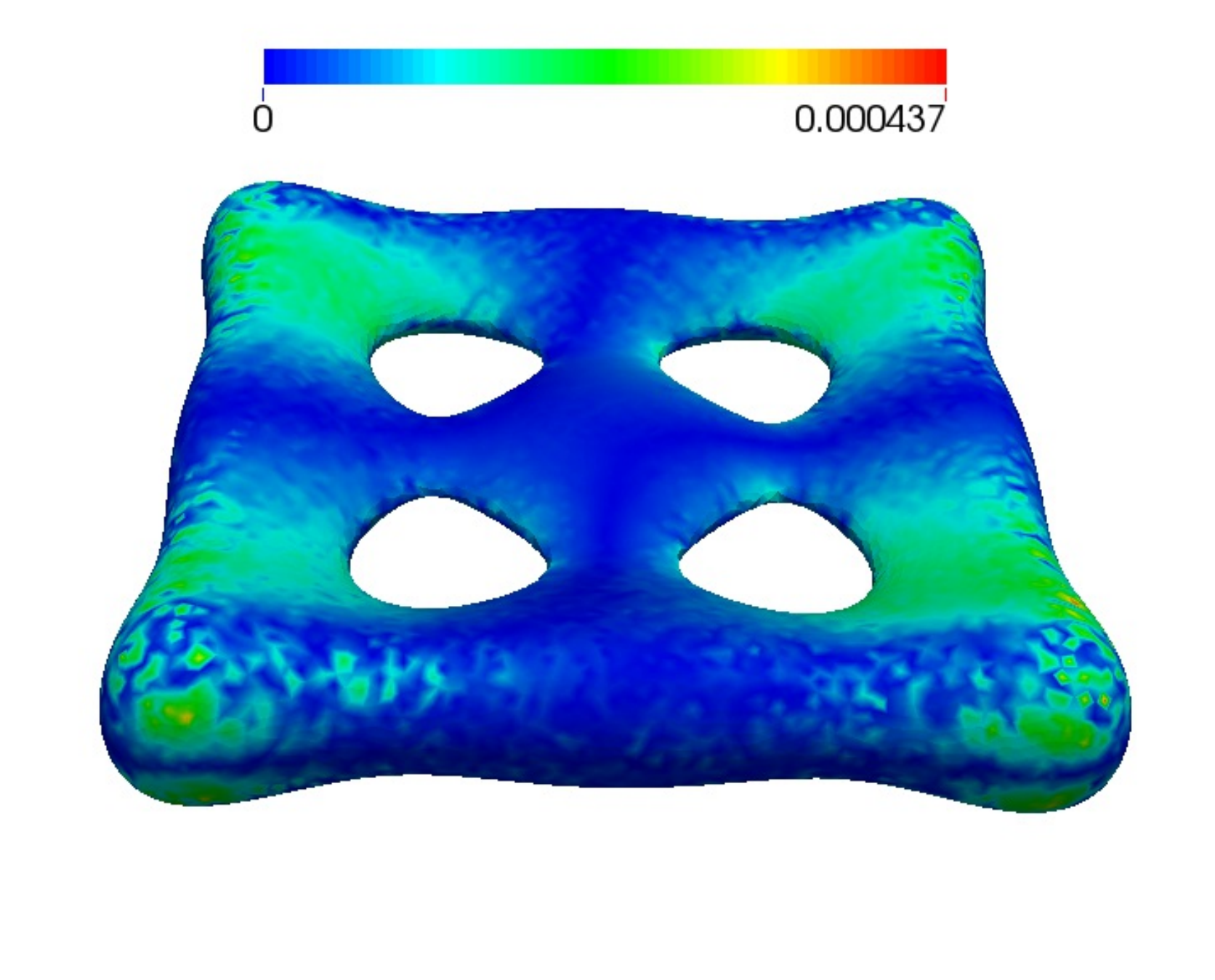}
 }
 \caption{Snapshots of the interpolant of the error using the two different schemes for Example \ref{eg:comparison}, the left hand column corresponds to the Lagrangian scheme and the righthand column corresponds to the ALE scheme.}\label{fig:complex_errors}
 \end{figure}
\end{Example}

\changes{
\begin{Example}[Simulation on a surface with changing conormal]\label{eg:graph}
We compute on a graph $\Gt$ above the unit disc which is given by the following parameterisation
\beq\label{eqn:gamma_disc}
\vec x(\vec \theta,t)=\left(\theta_1,\theta_2,2\sin(2\pi t)(1-\theta_1^2-\theta_2^2)\right),\qquad\vec \theta=(\theta_1,\theta_2)^T\in B_1(0),t\in[0,0.25].
\eeq
Defining the height of the graph
\beq
z(\vec \theta,t)=2\sin(2\pi t)(1-\theta_1^2-\theta_2^2),\qquad\vec \theta=(\theta_1,\theta_2)^T\in B_1(0),t\in[0,0.25],
\eeq
we set the material velocity to be the normal velocity of the graph which is given by 
\beq\label{eqn:vel_Lag_graph}
\vec v(\vec \theta,t)=\frac{-\pdt z(\vec \theta,t)\left((\nabla z(\vec \theta,t))^T,-1)^T\right)}{1+\lv\nabla z(\vec \theta,t)\rv^2},\qquad \vec \theta\in B_1(0),t\in[0,0.25].
\eeq
We will again compare a Lagrangian and ALE scheme. For the ALE scheme we define the arbitrary velocity 
\beq\label{eqn:vel_ALE_graph}
\vec v^a_1(\vec \theta,t)=0,\vec v^a_2(\vec \theta,t)=0,\vec v^a_3(\vec \theta,t)=\pdt z(\vec \theta,t)\qquad \vec \theta\in B_1(0),t\in[0,0.25].
\eeq
The arbitrary tangential velocity is then determined by $\vec a_\tangent=\vec v^a-\vec v$.

For this example we define the initial triangulation $\Gctn{0}$ (which is used for both schemes) with curved boundary faces in such a way that the initial triangulation is an exact triangulation of the unit disc, i.e., $\Gctn{0}=\G^0$. We also note that as the the  velocity fields $\vec v$ and $\vec v^a$, defined by (\ref{eqn:vel_Lag_graph}) and (\ref{eqn:vel_ALE_graph}) respectively, are zero on the boundary, the triangulation of the boundary remains exact for all times.

    \begin{figure}[ht]
 \centering 
    \subfigure[][{$t=0.05$}]{
  \includegraphics[ trim =10mm 200mm 10mm 10mm,  clip,    width=0.31\linewidth
 ]{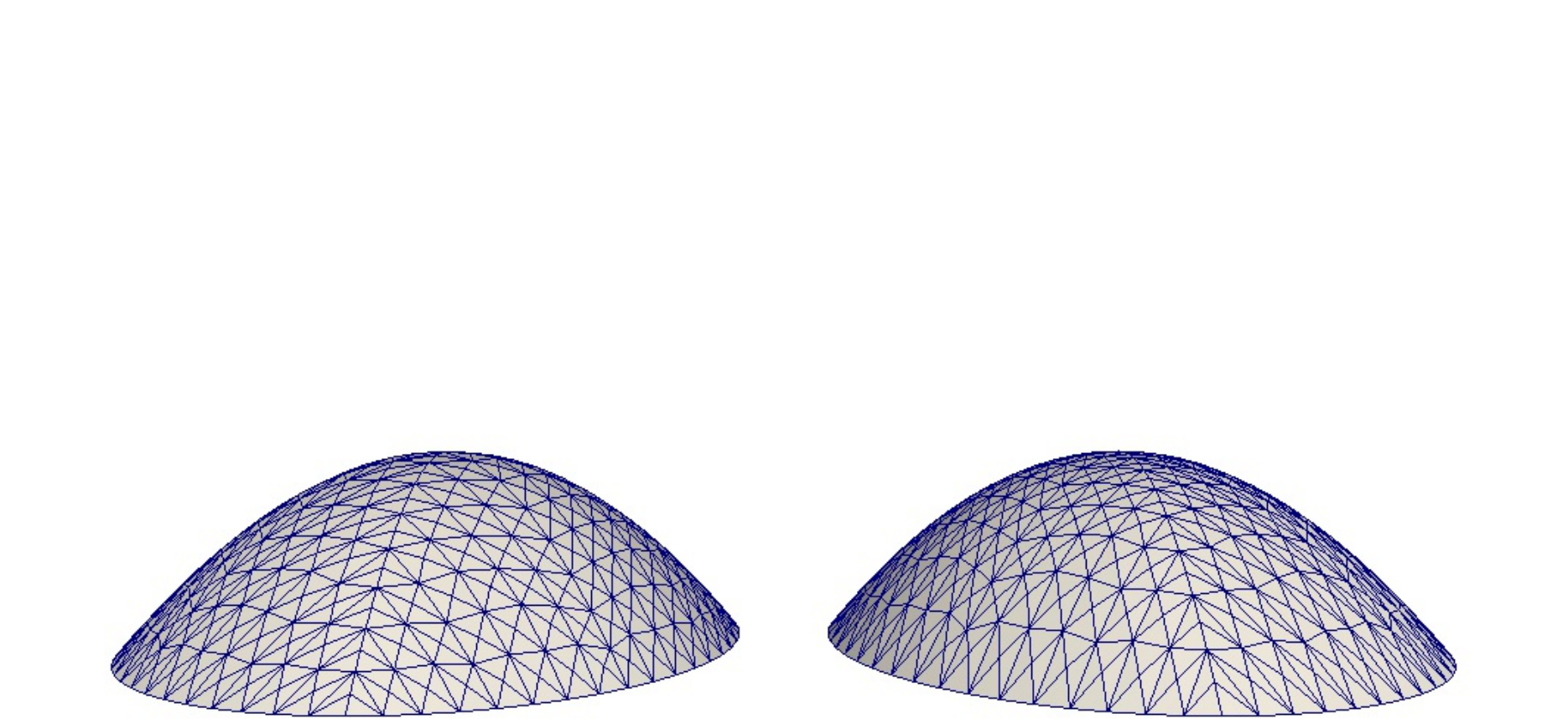}
}   
\subfigure[][{$t=0.1$}]{
  \includegraphics[ trim =10mm 200mm 10mm 10mm,  clip,    width=0.31\linewidth
 ]{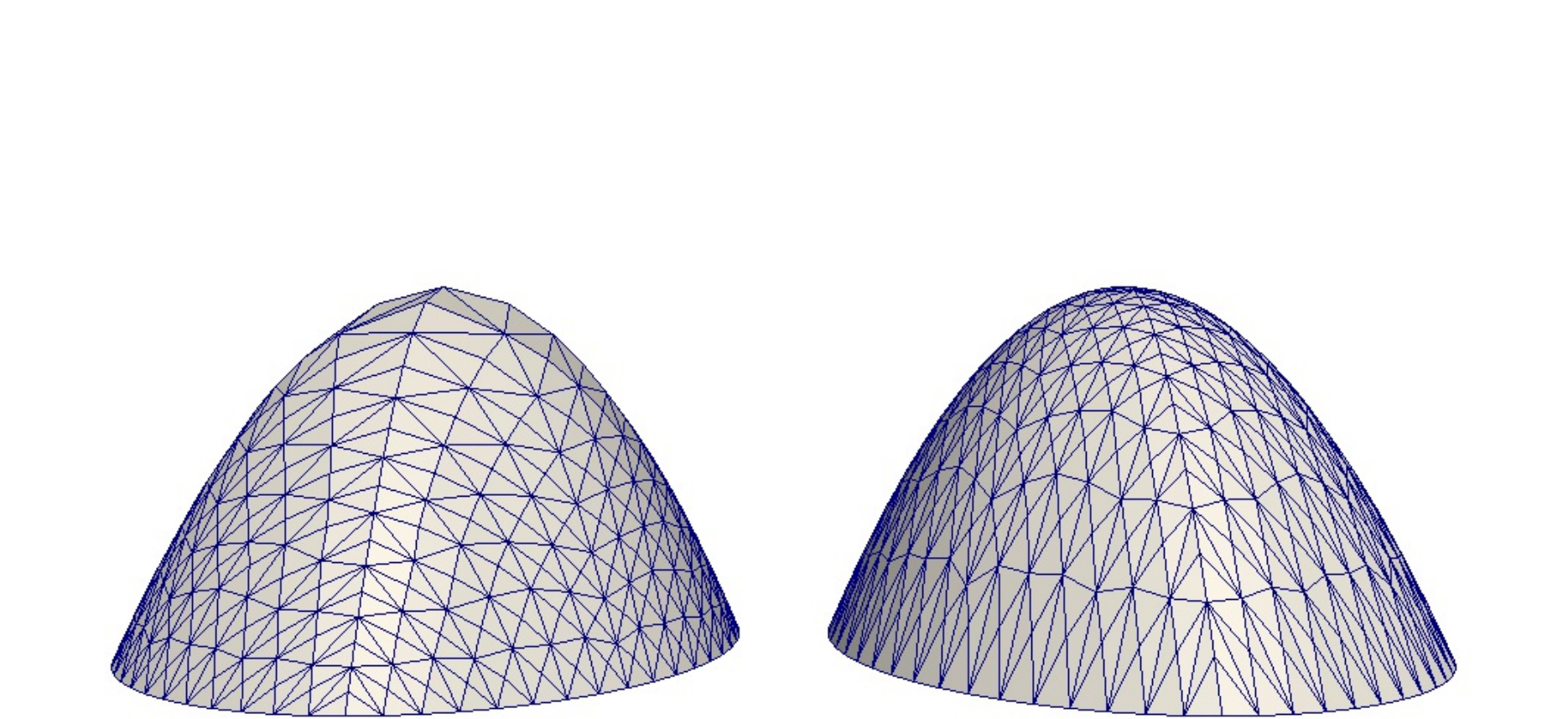}
}
\subfigure[][{$t=0.25$}]{
  \includegraphics[ trim =10mm 200mm 10mm 10mm,  clip,    width=0.31\linewidth
 ]{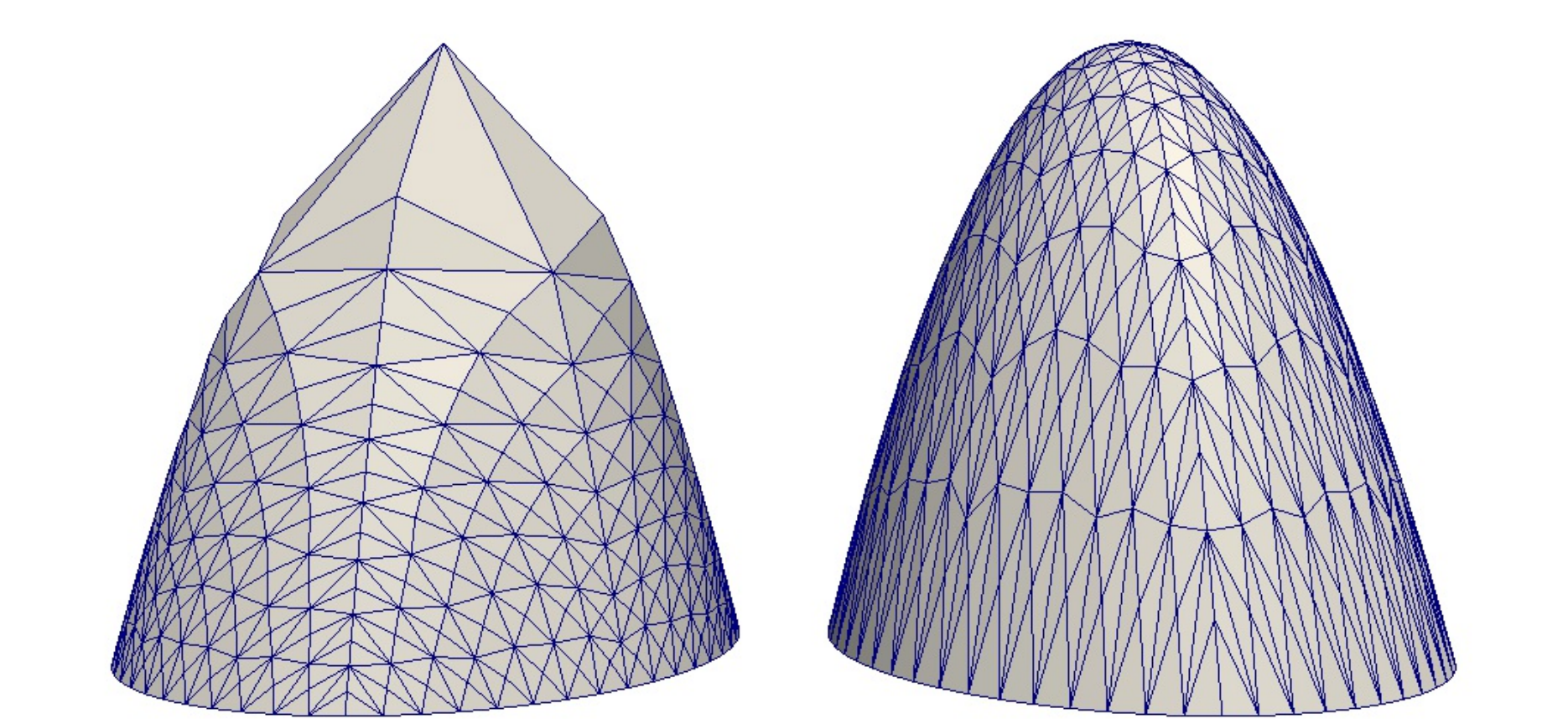}
}
 \caption{Snapshots of the meshes obtained in Example \ref{eg:graph} on moving the vertices with the Lagrangian (\ref{eqn:vel_Lag_graph}) and ALE velocities (\ref{eqn:vel_ALE_graph}), left and right hand meshes in each subfigure respectively, starting with the same initial triangulation with 545 vertices.
 We observe that moving the vertices with the Lagrangian velocity  leads to a mesh that resolves the surface poorly at the final time.}\label{fig:graph_meshes}
 \end{figure}

In  Figure \ref{fig:graph_meshes} we show some snapshots of the evolution of the same initial triangulation using the Lagrangian and ALE velocities, we have used a coarse initial triangulation so that the individual elements are clearly visible. For this example the ALE velocity clearly yields a mesh more suitable for computation.

 We consider the following equation on the surface (\ref{eqn:gamma_disc});
 \beq\label{eqn:conormal_pde}
\mdt{\vec v} u+u\nabla_{\Gt}\cdot\vec v-\lap_{\Gt}u=10\sin\left(2\pi x_3^2\right)\quad\mbox{ on }\Gt, t\in[0,0.25],
\eeq
with natural boundary conditions of the form (\ref{eqn:BCs}). We take the initial data 
 $u(\vec x,0)=0$. We selected a timestep of $10^{-5}$. We employed the BDF1 (implicit Euler) scheme (\ref{eqn:fd_scheme}) to compute the discrete solutions.

    \begin{figure}[ht]
 \centering 
    \subfigure[][{Mesh with 8321 vertices}]{
  \includegraphics[ trim =40mm 200mm 40mm 10mm,  clip,    width=0.31\linewidth
 ]{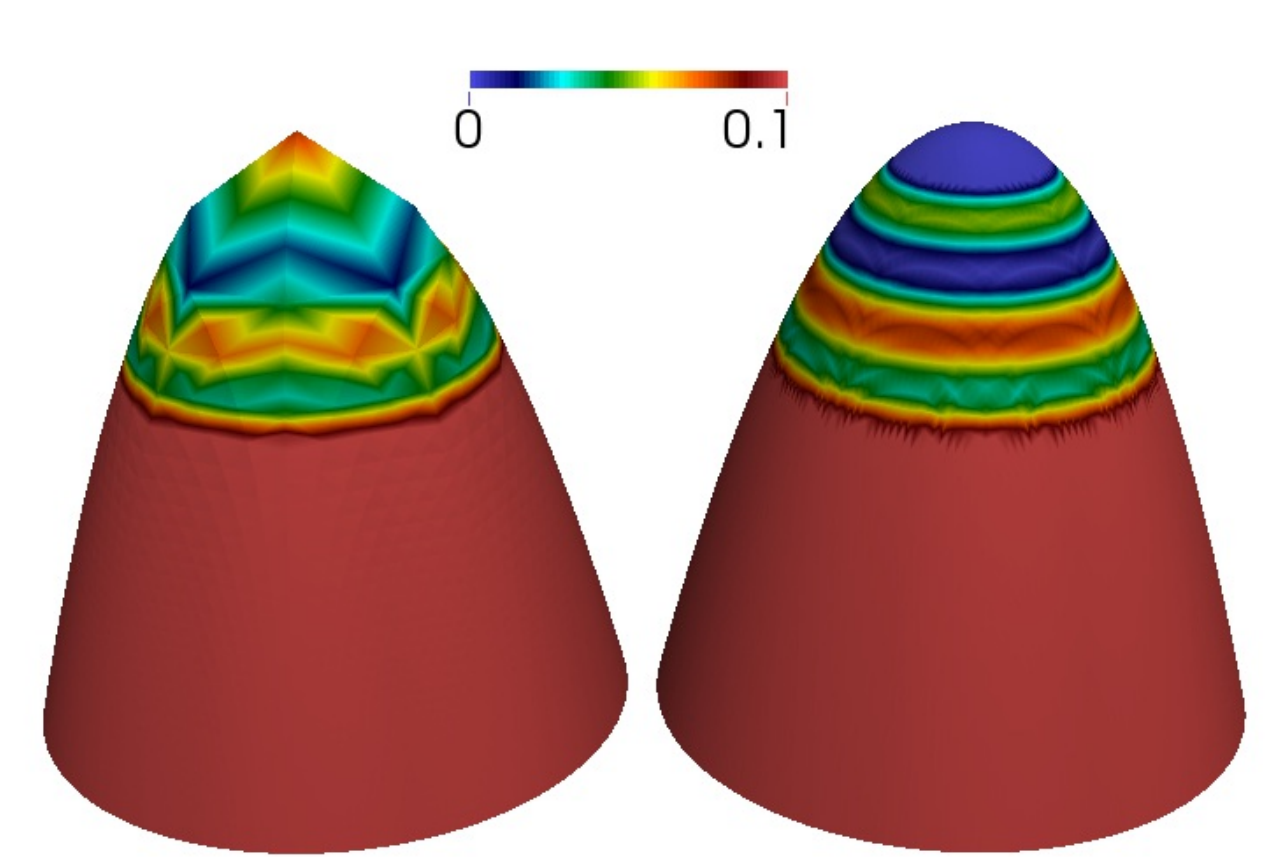}
}   
\subfigure[][{Mesh with 33025 vertices}]{
  \includegraphics[ trim =40mm 200mm 40mm 10mm,  clip,    width=0.31\linewidth
 ]{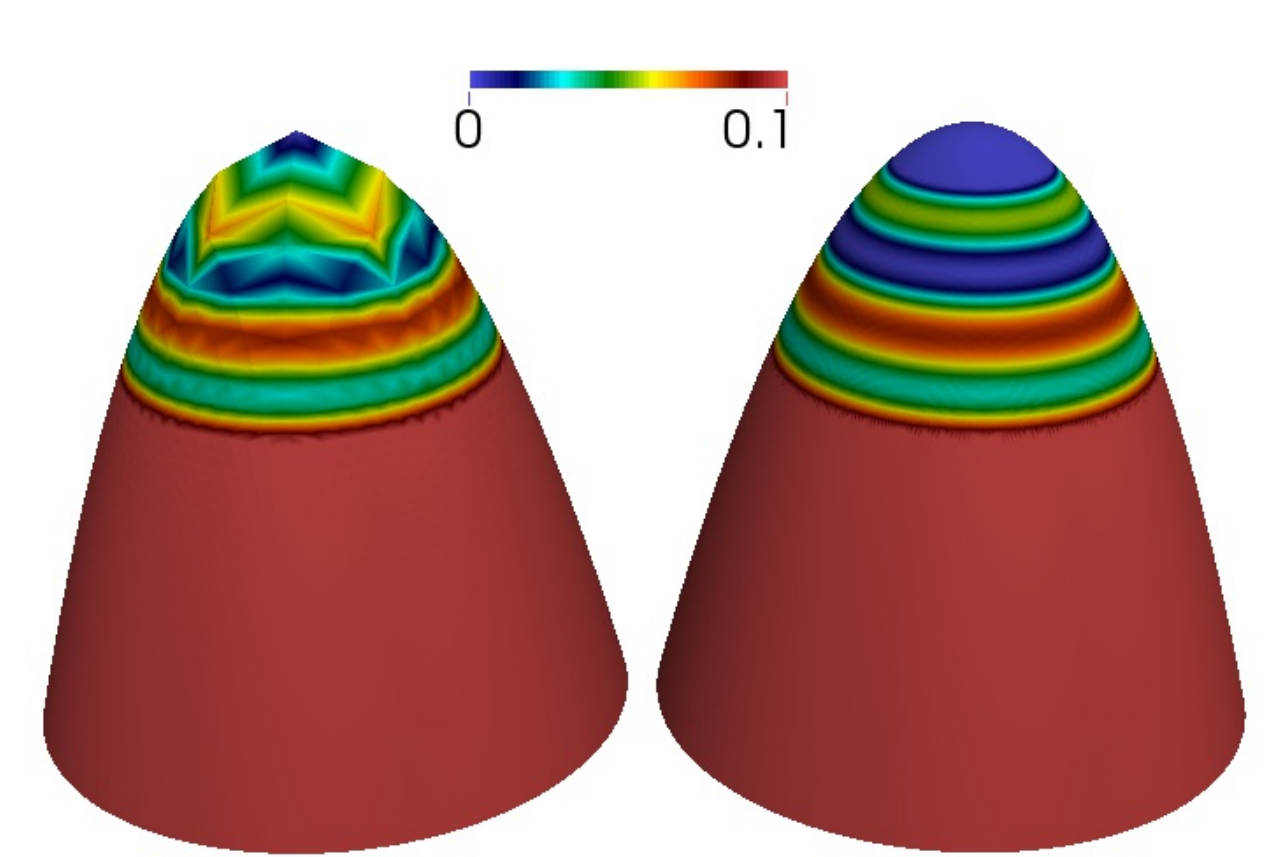}
}
\subfigure[][{Mesh with 131585 vertices}]{
  \includegraphics[ trim =40mm 200mm 40mm 10mm,  clip,    width=0.31\linewidth
 ]{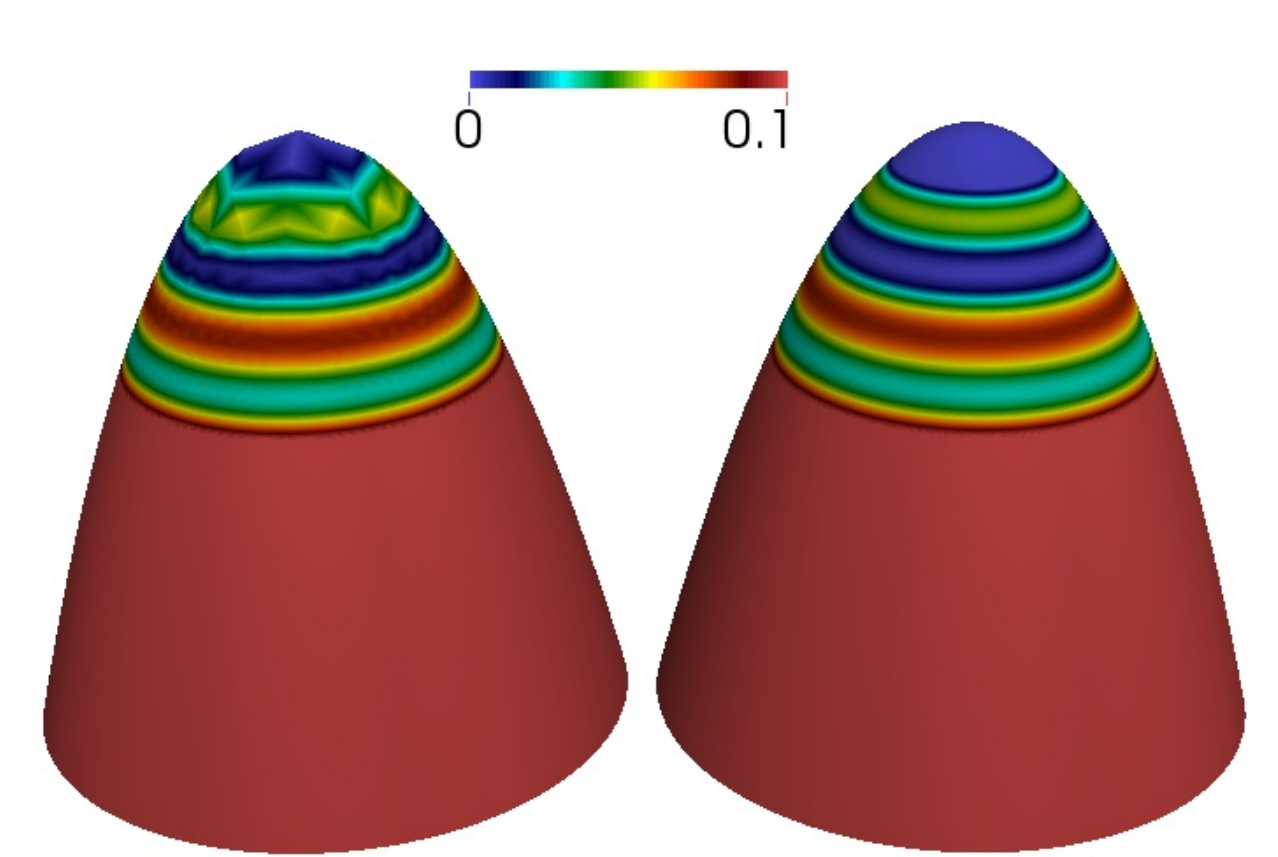}
}
 \caption{Snapshots of the computed solution for Example \ref{eg:graph} at the end time $t=0.25$ (slightly tilted). Each subfigure shows  the results of the Lagrangian  (left) and ALE schemes (right) starting with the same initial triangulations for both schemes.}\label{fig:graph_solutions}
 \end{figure}
 
 Figure \ref{fig:graph_solutions} shows the computed to solution to (\ref{eqn:conormal_pde}) at the final time on successive refinements of the mesh with the ALE and Lagrangian schemes.
 We observe good agreement between the solutions with the coarser and finer meshes in the ALE case and qualitative agreement between these solutions and the solution with Lagrangian scheme on the finest mesh (although even with the finest mesh the resolution of the surface is poor in the Lagrangian case). On the coarser meshes the Lagrangian scheme does not adequately resolve the surface and the source term and hence generates qualitatively different solutions to the fine mesh Lagrangian and (all three) ALE simulations.
\end{Example}
}

\begin{Example}[\changes{Long time Lagrangian simulations on a surface with periodic evolution}]\label{eg:per}\margnote{ref 2. pt 50.}
We consider a surface  
\beq\label{eqn:periodic_LS}
\Gt=\left\{\vec x\in\Reals^3\lv\frac{x_1^2}{a(t)^2}+\frac{x_2^2}{b(t)^2}+\frac{x_3^2}{c(t)^2}-1=0\right.\right\},
\eeq
with $a(t)= 1-0.1\sin(\pi t), b(t)=1-0.2\sin(\pi t)$ and $c(t)=1+0.1\sin(\pi t)$. The surface is therefore an ellipsoid with time dependent axes and the initial surface at $t=0$ is the surface of the unit sphere. We assume the material velocity of the surface is the normal velocity. 
We consider (\ref{eqn:pde}) posed on the surface with four different initial conditions
\begin{align}
\label{eqn:IC1}
u_1(\vec x,0)&=1\quad\vec x\in\G^0,
\\
\label{eqn:IC2}
u_2(\vec x,0)&=1+\sin(2 \pi x_1)\quad\vec x\in\G^0,
\\
\label{eqn:IC3}
u_3(\vec x,0)&=1+4\sin(8 \pi x_1)+3\cos(6\pi x_2)+2\sin(8 \pi x_3)\quad\vec x\in\G^0,
\\
\label{eqn:IC4}
u_4(\vec x,0)&=1+8\sin(16 \pi x_1)+7\cos(14\pi x_2)+6\sin(24 \pi x_3)\quad\vec x\in\G^0.
\end{align}

We used the Lagrangian BDF1 scheme (\ref{eqn:fd_scheme}) to simulate the equation on a triangulation of the sphere with 
16386 vertices and selected a timestep of $10^{-4}$. We approximated the initial data for the numerical method as follows
\begin{align}
\label{eqn:IC_h1}
U_{h,1}(\vec x,0)&=1\quad\vec x\in\Gctn{0},
\\
\label{eqn:IC_h2}
U_{h,2}(\vec x,0)&=\tilde{I}_h u_2(\vec x,0)+\int_{\Gctn{0}}\left(1-\tilde{I}_h u_2(\cdot,0)\right)\quad\vec x\in\Gctn{0},
\\
\label{eqn:IC_h3}
U_{h,3}(\vec x,0)&=\tilde{I}_h u_3(\vec x,0)+\int_{\Gctn{0}}\left(1-\tilde{I}_h u_3(\cdot,0)\right)\quad\vec x\in\Gctn{0},
\\
\label{eqn:IC_h4}
U_{h,4}(\vec x,0)&=\tilde{I}_h u_4(\vec x,0)+\int_{\Gctn{0}}\left(1-\tilde{I}_h u_4(\cdot,0)\right)\quad\vec x\in\Gctn{0},
\end{align}
where $\tilde{I}_h:\Cont{}(\G^0)\to\Sc(0)$ denotes the linear Lagrange interpolation operator. The approximations of the initial conditions for the numerical scheme were chosen such that  the initial approximations have the same total mass. We note that the approximations of the the initial conditions satisfy  (\ref{eqn:IC_approx}). 

Figure \ref{fig:IC_h} shows plots of the initial conditions (\ref{eqn:IC_h2}), (\ref{eqn:IC_h3}) and (\ref{eqn:IC_h4}) on the discrete surface. Figure \ref{fig:per_solution} shows snapshots of the discrete solution for the case of constant initial data (\ref{eqn:IC_h1}) we observe that the numerical solution  appears to converge rapidly to a periodic  function. We wish to investigate numerically the effect of the initial data on  this periodic solution, to this end we compute the numerical solution
with the initial conditions (\ref{eqn:IC_h2}), (\ref{eqn:IC_h3}) and (\ref{eqn:IC_h4}) and compare these numerical solutions to that obtained with constant initial data. Figure \ref{fig:per_comparison} shows the $\Lp{2}(\Gc(t))$ norm of the difference between the numerical solutions with the non-constant initial data and the numerical solution with the constant initial data  versus time. It appears that the numerical solutions converge to the same periodic solution for all four different initial conditions.
    \begin{figure}[ht]
 \centering 
    \subfigure[][{$U_{h,2}(\cdot,0)$}]{
  \includegraphics[ trim = 20mm 130mm 20mm 0mm,  clip,    width=0.31\linewidth
 ]{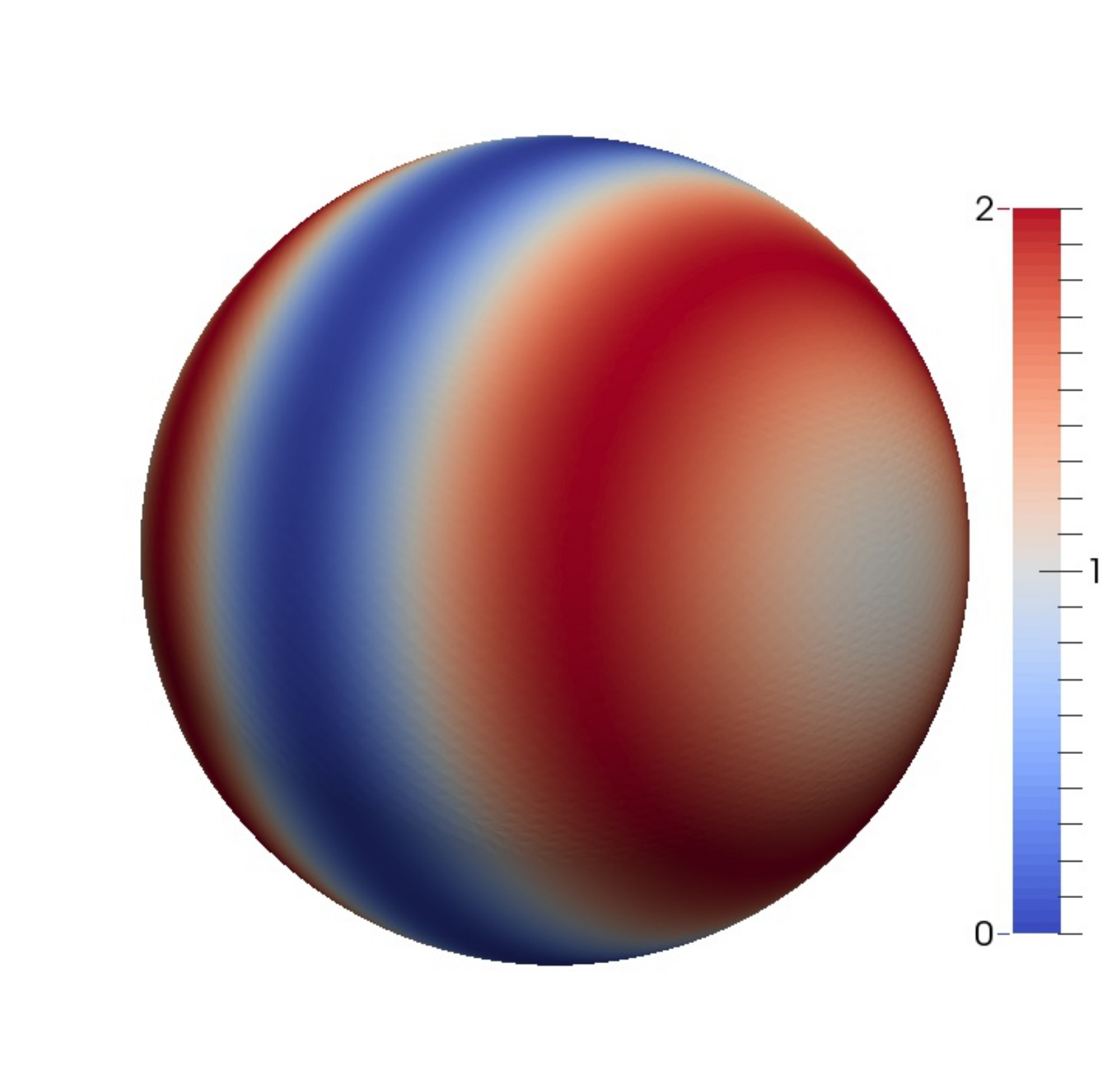}
}
    \subfigure[][{$U_{h,3}(\cdot,0)$}]{
  \includegraphics[ trim = 20mm 130mm 20mm 0mm,  clip,    width=0.31\linewidth
 ]{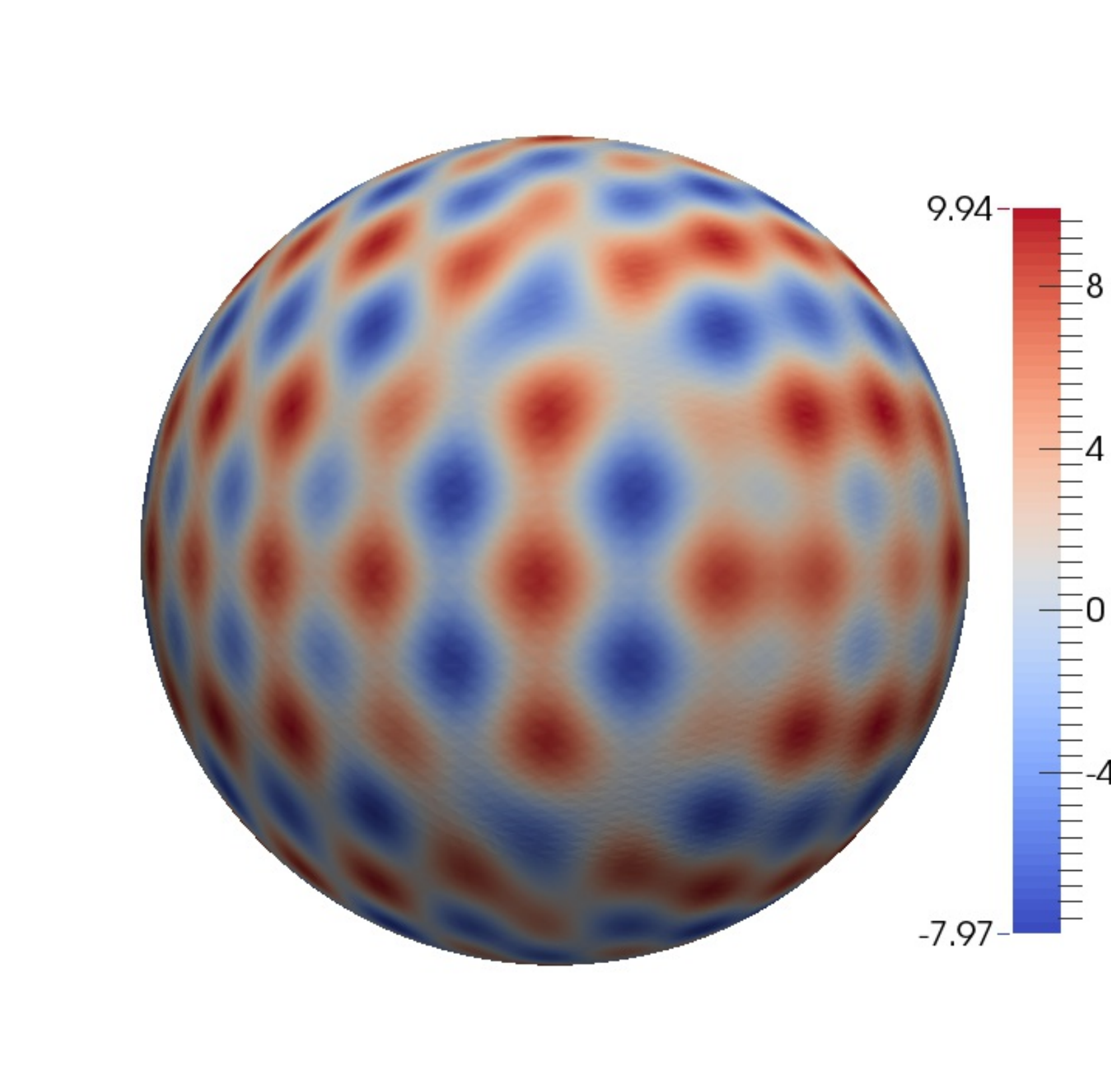}
}
    \subfigure[][{$U_{h,4}(\cdot,0)$}]{
  \includegraphics[ trim = 20mm 130mm 20mm 0mm,  clip,    width=0.31\linewidth
 ]{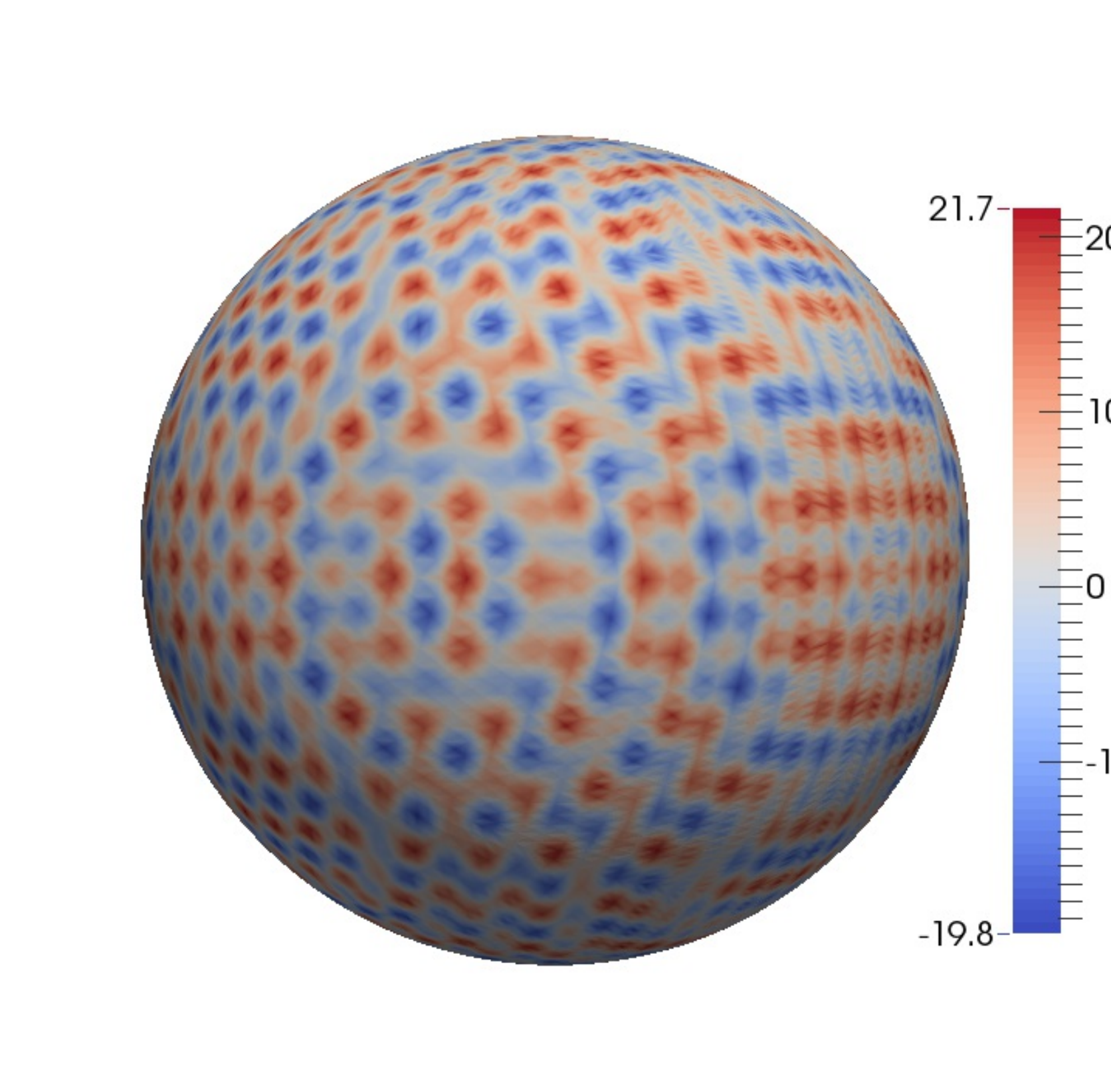}
}
 \caption{Initial conditions (\ref{eqn:IC_h2}), (\ref{eqn:IC_h3}) and (\ref{eqn:IC_h4})  for the simulations of Example \ref{eg:per}.}\label{fig:IC_h}
 \end{figure}

    \begin{figure}[ht]
 \centering 
     \subfigure[][{$t=0$}]{
  \includegraphics[ trim = 20mm 0mm 10mm 0mm,  clip,    width=0.23\linewidth
 ]{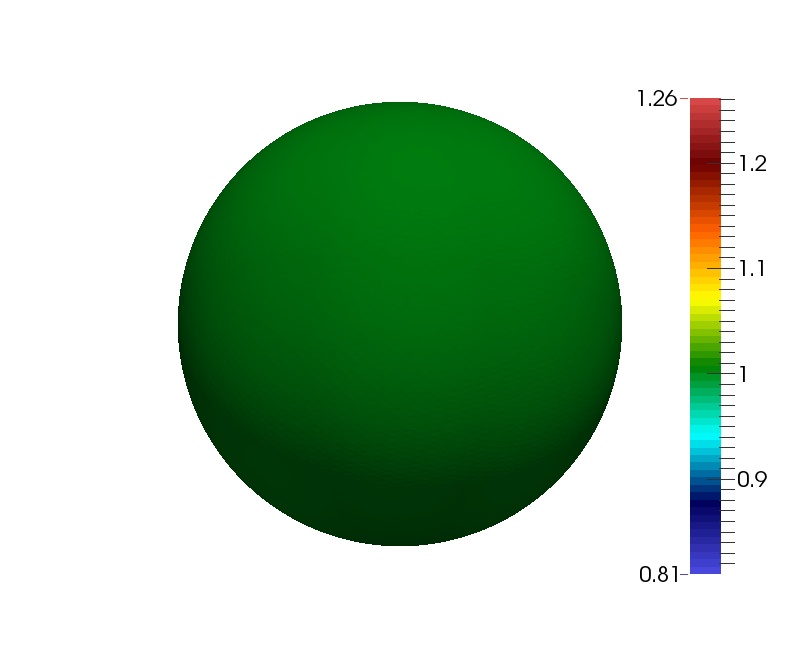}
}
    \subfigure[][{$t=0.2$}]{
  \includegraphics[ trim = 20mm 0mm 10mm 0mm,  clip,    width=0.23\linewidth
 ]{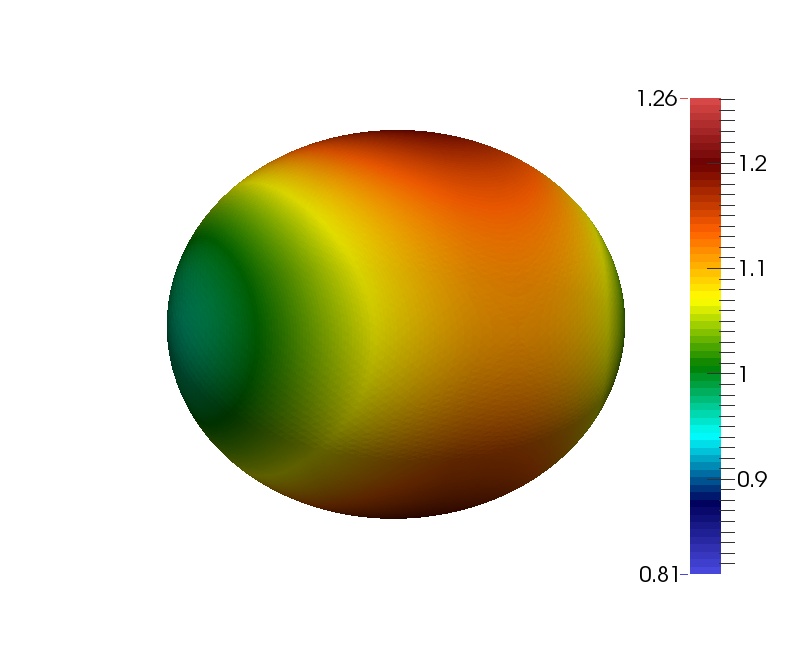}
}
    \subfigure[][{$t=1$}]{
  \includegraphics[ trim = 20mm 0mm 10mm 0mm,  clip,    width=0.23\linewidth
 ]{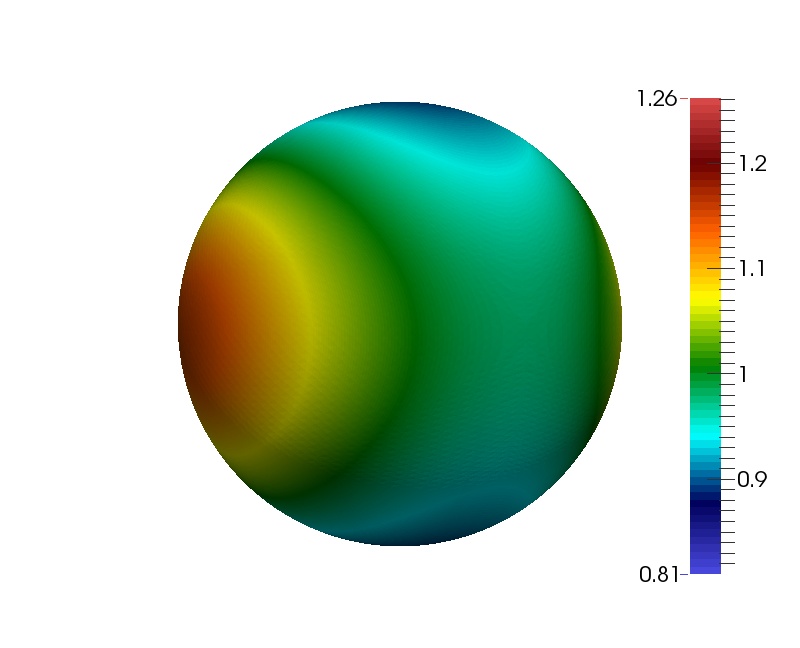}
}
    \subfigure[][{$t=1.5$}]{
  \includegraphics[ trim = 20mm 0mm 10mm 0mm,  clip,    width=0.23\linewidth
 ]{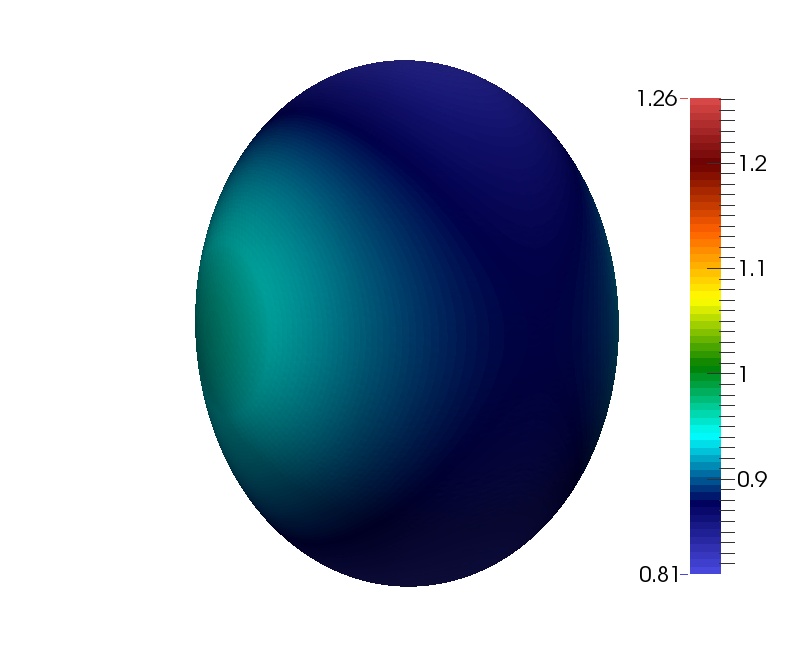}
}
    \subfigure[][{$t=2$}]{
  \includegraphics[ trim = 20mm 0mm 10mm 0mm,  clip,    width=0.23\linewidth
 ]{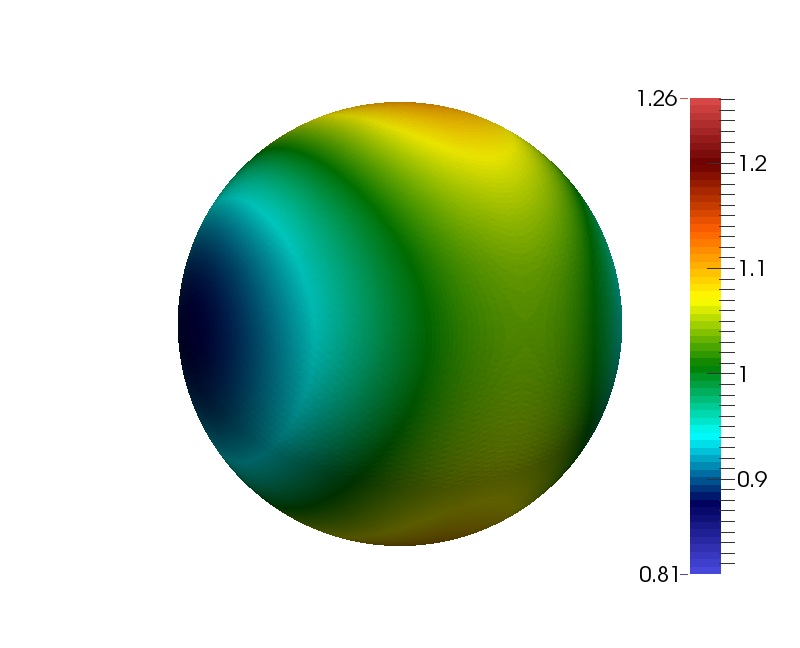}
}
    \subfigure[][{$t=2.2$}]{
  \includegraphics[ trim = 20mm 0mm 10mm 0mm,  clip,    width=0.23\linewidth
 ]{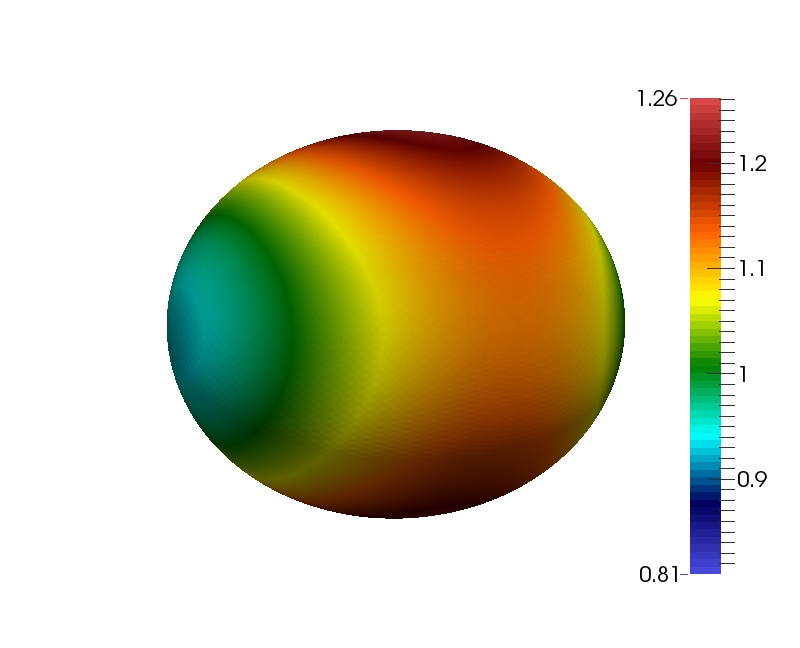}
}
    \subfigure[][{$t=3$}]{
  \includegraphics[ trim = 20mm 0mm 10mm 0mm,  clip,    width=0.23\linewidth
 ]{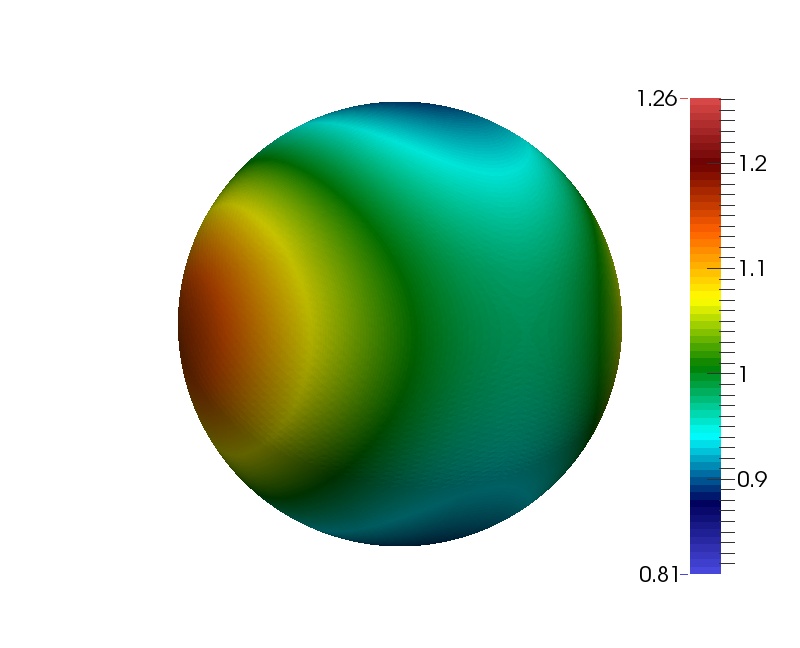}
}
    \subfigure[][{$t=3.5$}]{
  \includegraphics[ trim = 20mm 0mm 10mm 0mm,  clip,    width=0.23\linewidth
 ]{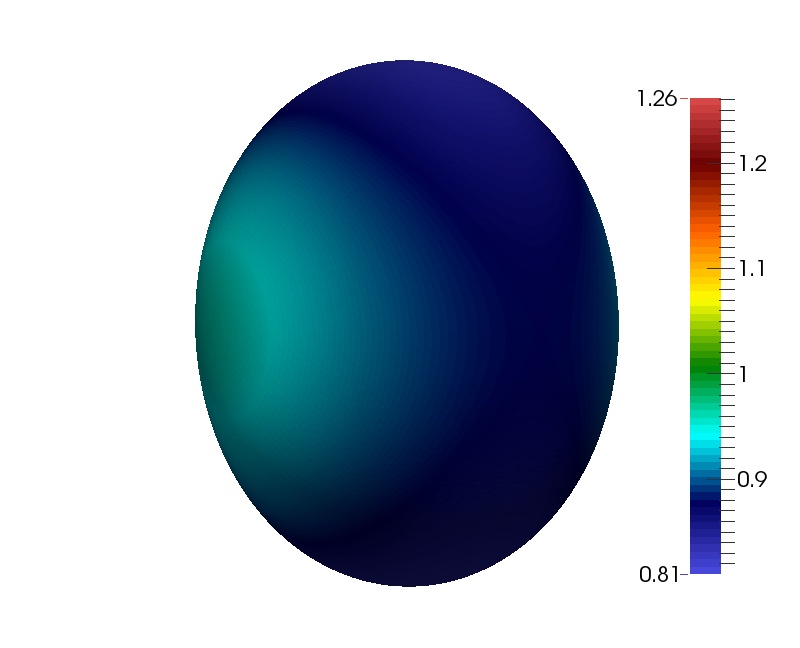}
}
    \subfigure[][{$t=4$}]{
  \includegraphics[ trim = 20mm 0mm 10mm 0mm,  clip,    width=0.23\linewidth
 ]{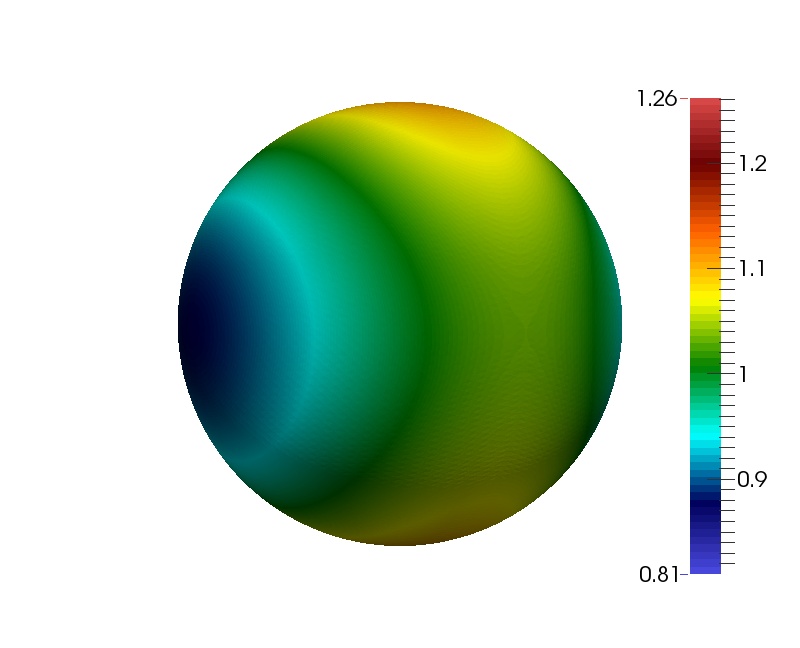}
}
    \subfigure[][{$t=4.2$}]{
  \includegraphics[ trim = 20mm 0mm 10mm 0mm,  clip,    width=0.23\linewidth
 ]{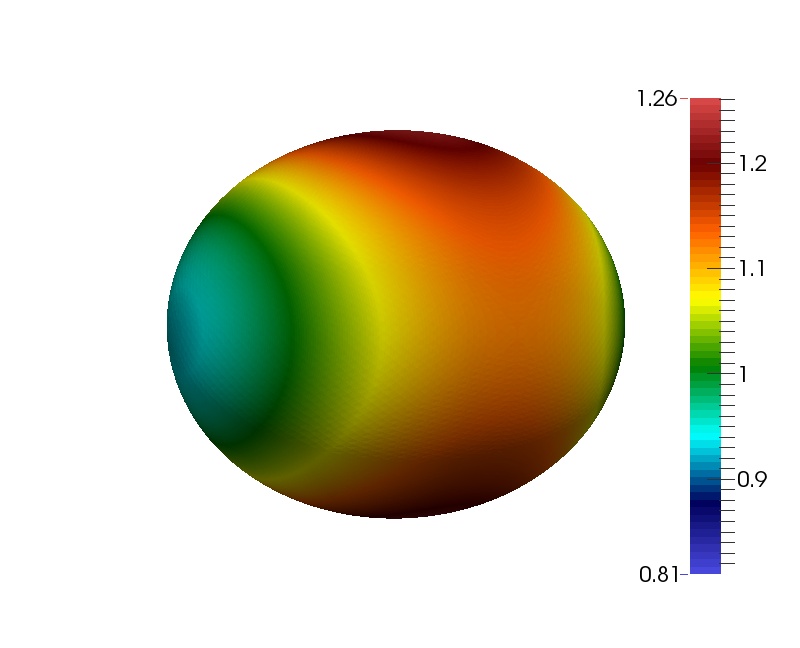}
}
    \subfigure[][{$t=5$}]{
  \includegraphics[ trim = 20mm 0mm 10mm 0mm,  clip,    width=0.23\linewidth
 ]{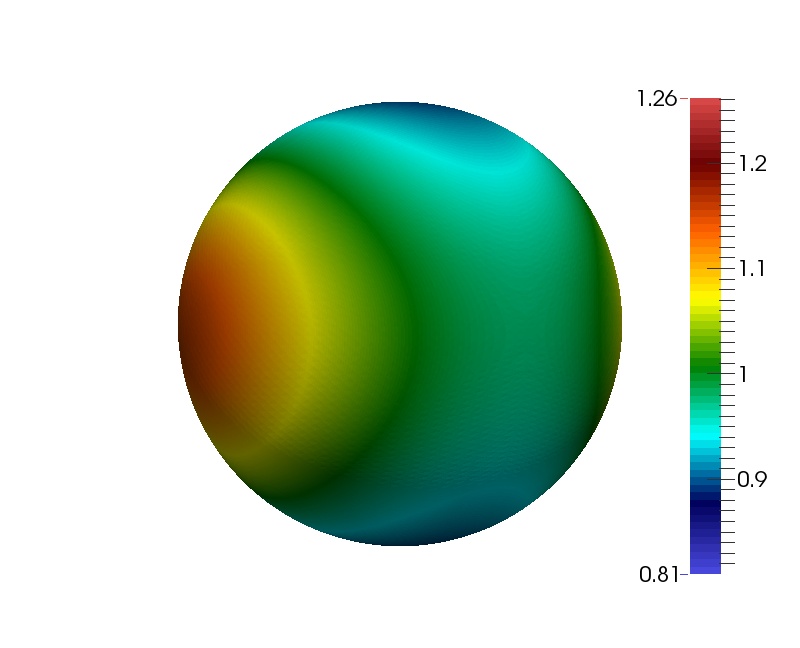}
}
    \subfigure[][{$t=5.5$}]{
  \includegraphics[ trim = 20mm 0mm 10mm 0mm,  clip,    width=0.23\linewidth
 ]{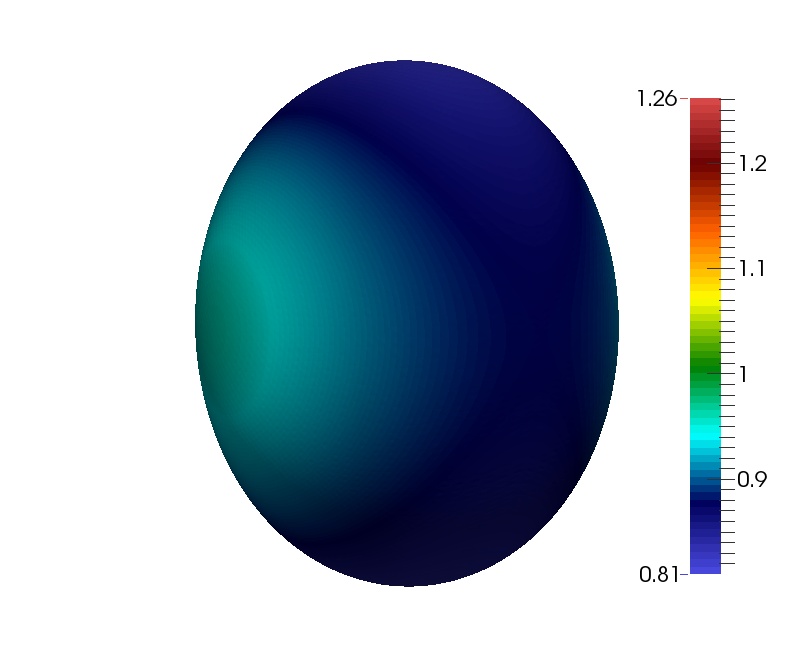}
}
 \caption{ Snapshots of the numerical solution of Example \ref{eg:per} with constant initial data (\ref{eqn:IC_h1}).}\label{fig:per_solution}
 \end{figure}
 
     \begin{figure}[ht]
 \centering 
  \includegraphics[ trim = 0mm 0mm 0mm 0mm,  clip,    width=0.75\linewidth
 ]{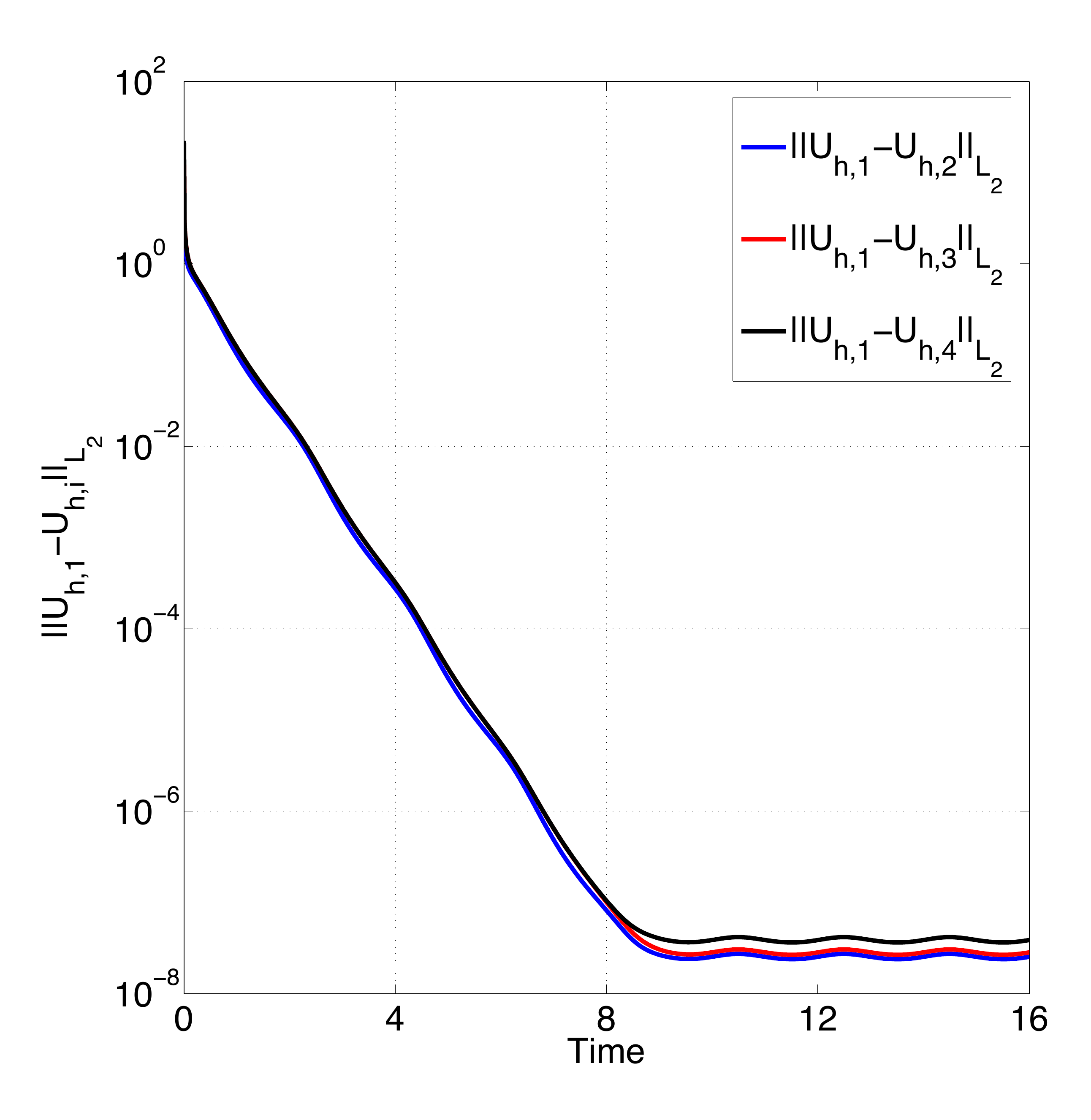}
\caption{ The $\Lp{2}(\Gc(t))$ norm of the difference between numerical solution with constant initial data (\ref{eqn:IC_h1})  and the numerical solutions corresponding to the initial conditions (\ref{eqn:IC_h2}), (\ref{eqn:IC_h3}) and (\ref{eqn:IC_h4}) (blue, red and black lines respectively).}\label{fig:per_comparison}
 \end{figure}
\end{Example}

\appendix
\section{Transport formula}
The following transport formula play a fundamental role in the formulation and analysis of the numerical method.
\begin{Lem}[Transport formula]\label{lem:transport}
Let $\surface(t)$ be a smoothly evolving surface with material velocity $\vec v$, let $f$ and $g$ be sufficiently smooth functions and $\vec w$ a sufficiently smooth vector
field such that all the following quantities exist. Then
\begin{align}
\frac{\diff }{\diff t}\int_{\surface(t)}f=&
\int_{\surface(t)}\mdt{\vec v}f+f\nabla_\G\cdot \vec v\label{eqn:transport_scalar},\\
\frac{\diff }{\diff t}\int_{\surface(t)} f\vec w\cdot\nabla_\G g=&
\int_{\surface(t)}\left( \mdt{\vec v}f\vec w\cdot\nabla_\G g+ f\mdt{\vec v}\vec w\cdot\nabla_\G g+ f\vec w\cdot\nabla_\G \mdt{\vec v}g\right)\label{eqn:transport_bbil}\\
\notag&+\int_{\surface(t)}\nabla_\G\cdot\vec v\left(f\vec w\cdot\nabla_\G g\right)-\int_{\surface(t)}f\vec w\cdot \defB{\vec v}{\G}\nabla_\G g\\
\frac{\diff }{\diff t}\int_{\surface(t)}\nabla_\G f\cdot\nabla_\G g=&
\int_{\surface(t)}\left(\nabla_\G \mdt{\vec v}f\cdot\nabla_\G g+\nabla_\G \mdt{\vec v}g\cdot\nabla_\G f\right)\label{eqn:transport_dirichlet_ip}\\
\notag&+\int_{\surface(t)}\left(\nabla_\G\cdot\vec v-2\defD{\vec v}{\G}\right)\nabla_\G f\cdot\nabla_\G g,
\end{align}
with the deformation tensors defined by
$$
\defB{\vec v_{ij}}{\G}=(\nabla_\G)_i v_j-\sum_{l=1}^{m+1}\nu_l\nu_i(\nabla_\G)_jv_l\quad\text{ and }\quad
\defD{\vec v_{ij}}{\G}=\frac{1}{2}\left((\nabla_\G)_iv_j+(\nabla_\G)_jv_i\right),
$$ 
respectively.
\end{Lem}
\begin{Proof}
Proofs of (\ref{eqn:transport_scalar}) and (\ref{eqn:transport_dirichlet_ip}) are given in \citep{dziuk2007finite}. 
The proof of (\ref{eqn:transport_bbil}) is as follows, (for further details see the proof of (\ref{eqn:transport_dirichlet_ip})
in \citep[Appendix]{dziuk2007finite}) we have
\begin{align}\label{eqn:A4}
&\frac{\diff }{\diff t}\int_{\surface(t)} f\vec w\cdot\nabla_\G g=\int_{\surface(t)}\mdt{\vec v}\left(f\vec w\cdot\nabla_\G g\right)+\nabla_\G\cdot\vec v\left(f\vec w\cdot\nabla_\G g\right)\\
&=\int_{\surface(t)}\bigg(\mdt{\vec v}f\left(\vec w\cdot\nabla_\G g\right)+f\left(\mdt{\vec v}\vec w\right)\cdot\nabla_\G g+f\vec w\cdot\left(\mdt{\vec v}\nabla_\G g\right)+\nabla_\G\cdot\vec v\left(f\vec w\cdot\nabla_\G g\right)\bigg)\notag.
\end{align}
Finally, application of the following result from \citep[Lemma 2.6]{2013arXiv1307.1056D}
\[
\mdt{\vec v}\nabla_{\Gt}g=\nabla_{\Gt}\mdt{\vec v}g-\defB{\vec v}{\G}\nabla_{\Gt}g,
\]
in (\ref{eqn:A4}) completes the proof of the Lemma.
\end{Proof}

For the analysis of the second order scheme we note that repeated application of the transport formula together with the smoothness of the velocity yields the following bounds,\changes{\margnote{ ref 2. pt. 54.} see \citep[Lemma 9.1]{dziuk2011runge} for a similar discussion.} Let $\G$ be a smoothly evolving surface with material velocity $\vec v$, let $f$ and $g$ be sufficiently smooth functions and $\vec w$ a sufficiently smooth vector and further assume $\mdt{\vec v}g=0$ then 
\begin{align}
\label{eqn:d2_a}
\lv\frac{\diff^2 }{\diff t^2}\int_{\G}\nabla_\G f\cdot\nabla_\G g\rv
&\leq
\lv\int_{\G}\nabla_\G\mdt{\vec v}(\mdt{\vec v} f)\cdot\nabla_\G g\rv\\
\notag
&+c\left(\lv\int_{\G}\nabla_\G\mdt{\vec v}f\cdot\nabla_\G g\rv+\lv\int_{\G}\nabla_\G f\cdot\nabla_\G g\rv\right),\\
\label{eqn:d2_b}
\lv\frac{\diff^2 }{\diff t^2}\int_{\G}f\vec w\cdot\nabla_\G g\rv
&\leq
\lv\int_{\G}\mdt{\vec v}(\mdt{\vec v} f)\vec w\cdot\nabla_\G g\rv\\
\notag
&+c\left(\lv\int_{\G}\mdt{\vec v}f\vec w\cdot\nabla_\G g\rv+\lv\int_{\G}f\vec w\cdot\nabla_\G g\rv\right).
\end{align}
 \section{Approximation results}
  For a function $\eta\in C^0(\Gt)$ we denote by
 $I_h\eta\in\Scl$  the lift of the linear Lagrange interpolant of $\tilde{I}_h\eta\in\Sc$, i.e., $I_h\eta=(\tilde{I}_h\eta)^l$. The following Lemma was shown in \cite{dziuk1988finite}.
 \begin{Lem}[Interpolation bounds]\label{Lem:interp}
 For an $\eta\in\Hil{2}{\Gt}$ there exists a unique $I_h\eta\in\Sclt$ such that
 \beq
 \ltwon{\eta-I_h\eta}{\Gt}+h\ltwon{\nabla_\Gt(\eta-I_h\eta)}{\Gt}\leq ch^2\Hiln{\eta}{2}{\Gt}.
 \eeq
 \end{Lem}
The following results provide estimates for the difference between the continuous velocity (here we mean the velocity that includes the arbitrary tangential motion and {\it not} the material velocity) and the discrete velocity of the smooth surface together with an estimate on the material derivative. 
\begin{Lem}{Velocity and material derivative estimates}
\begin{align}
\lv \mdth{\vec v^a_h}\left(\vec v_a - \vec v^a_h\right)\rv +h\lv\nabla_\Gt\mdth{\vec v^a_h} \left(\vec v_a-\vec v^a_h\right)\rv
&
\leq ch^2\quad \text{on }\G\label{mdt_velocity_bound}
\\
\ltwon{\vec a_\tangent - \vec t^a_h}{\Gt} +h\ltwon{\nabla_\Gt \left(\vec a_\tangent - \vec t^a_h\right)}{\Gt}&\leq ch^2\Hiln{\vec a_\tangent}{2}{\Gt}.\label{tang_velocity_bound}\\
\ltwon{\mdt{\vec v_a}z-\mdth{\vec v^a_h}z}{\Gt}&\leq ch^2\Hiln{z}{1}{\Gt}\label{md_l2_bound}\\
\ltwon{\nabla_\Gt\left(\mdt{\vec v_a}z-\mdth{\vec v^a_h}z\right)}{\Gt}&\leq ch\Hiln{z}{2}{\Gt}\label{md_grad_bound}\\
\ltwon{\mdt{\vec v_a}\mdt{\vec v_a}z-\mdt{\vec v^a_h}\mdth{\vec v^a_h}z}{\Gt}&\leq ch^2\Hiln{\mdt{\vec v}z}{1}{\Gt}\label{mdmd_l2_bound}\\
\ltwon{\nabla_\Gt\left(\mdt{\vec v_a}\mdt{\vec v_a}z-\mdt{\vec v^a_h}\mdth{\vec v^a_h}z\right)}{\Gt}&\leq ch\Hiln{\mdt{\vec v}z}{2}{\Gt}\label{mdmd_grad_bound}.
\end{align}
\end{Lem}
\begin{Proof}
The estimate  (\ref{mdt_velocity_bound}) is shown in \citep[Lemma 7.3]{lubich2012variational}, (\ref{tang_velocity_bound}) follows from Lemma \ref{Lem:interp} and the fact that $\vec T^a_h$ is the interpolant of the arbitrary tangential velocity and $\vec t^a_h$ is its lift.  Estimates (\ref{md_l2_bound}) and (\ref{md_grad_bound}) are shown in \citep[Cor. 5.7]{dziuk2010l2}. The estimates (\ref{mdmd_l2_bound}) and (\ref{mdmd_grad_bound}) follow easily from (\ref{mdt_velocity_bound}), (\ref{md_l2_bound}) and (\ref{md_grad_bound}).
\end{Proof}

We now state some results on the error due to the approximation of the surface
\begin{Lem}[Geometric perturbation errors]\label{lem:geom_pert_errors}
For any $(\Psi_h(\cdot,t),\Phi_h(\cdot,t))\in\Sc(t)\times\Sc(t)$ with corresponding lifts $(\psi_h(\cdot,t),\varphi_h(\cdot,t))\in\Scl(t)\times\Scl(t),$
the following bounds hold:
\begin{align}
\lv \mbil{\psi_h}{\varphi_h}-\mhbil{\Psi_h}{\Phi_h}\rv&\leq ch^2\ltwon{\psi_h}{\Gt}\ltwon{\varphi_h}{\Gt}\label{eqn:pert_m}\\
\lv \abil{\psi_h}{\varphi_h}-\ahbil{\Psi_h}{\Phi_h}\rv&\leq ch^2\ltwon{\nabla_{\Gt}\psi_h}{\Gt}\ltwon{\nabla_{\Gt}\varphi_h}{\Gt}\label{eqn:pert_a}\\
\lv \gbil{\psi_h}{\varphi_h}{\vec v^a_h}-\ghbil{\Psi_h}{\Phi_h}{\vec V^a_h}\rv&\leq ch^2\Hiln{\psi_h}{1}{\Gt}\Hiln{\varphi_h}{1}{\Gt}\label{eqn:pert_g}\\
\lv \bbil{\psi_h}{\varphi_h}{\vec t^a_h}-\bhbil{\Psi_h}{\Phi_h}{\vec T^a_h}\rv&\leq ch^2\ltwon{\psi_h}{\Gt}\ltwon{\nabla_\Gt\varphi_h}{\Gt},\label{eqn:pert_b}
\end{align}
with $\vec V^a_h,\vec T^a_h,\vec v^a_h$ and $\vec t^a_h$ as defined in \S \ref{sec:fe_disc}.
\end{Lem}
\begin{Proof}
A proof of (\ref{eqn:pert_m}), (\ref{eqn:pert_a}) and (\ref{eqn:pert_g}) is given in  \citep[Lemma 5.5]{dziuk2010l2}.
We now prove (\ref{eqn:pert_b}). We start by introducing some notation.
We denote by $\delta_h$ the quotient between the discrete and smooth surface measures which satisfies \cite[Lemma 5.1]{dziuk2007finite}
\beq\label{surface_element}
\sup_{t\in(0,T)}\sup_{\Gct}\lv 1 -\delta_h\rv\leq ch^2
\eeq
We  introduce $\Proj,\Proj_h$
the projections onto the tangent planes of $\Gt$ and $\Gc$ respectively. We denote by ${\vec{\mathcal H}}$ the Weingarten map 
$(\mathcal{H}_{ij}=\partial_{x_j}\nu_i)$. 
\begin{align}\label{pf_bbil_1}
\lv \bbil{\psi_h}{\varphi_h}{\vec t^a_h}-\bhbil{\Psi_h}{\Phi_h}{\vec T^a_h}\rv
=
\lv \int_\Gt\psi_h\vec t^a_h\cdot\nabla_\Gt\varphi_h-\int_{\Gct}\Psi_h\vec T^a_h\cdot\nabla_\Gct\Phi_h\rv
\end{align}
From \cite{dziuk2007finite} we have
\beq\label{eqn:nab_gc_nab_gt}
\nabla_\Gc\eta=\vec B_h\nabla_\G\eta^l,
\eeq
where $\vec B_h=\Proj_h(\vec I-d\vec{\mathcal{H}})$.  We have with $\vec p,\vec x$ as in (\ref{eqn:x_gct_p_gt}),
\begin{align}
\vec T^a_h(\vec x,\cdot)\cdot\nabla_\Gc \Phi_h(\vec x,\cdot)&=\Proj_h\vec T^a_h(\vec x,\cdot)\cdot\nabla_\Gc \Phi_h(\vec x,\cdot)\\
\notag
&=\Proj_h\vec t^a_h(\vec p,\cdot)\cdot \Proj_h(\vec I-d\vec{\mathcal{H}})\Proj\nabla_\G \varphi_h(\vec p,\cdot)\\
\notag
&=\Proj_h\vec t^a_h(\vec p,\cdot)\cdot \Proj_h\Proj(\vec I-d\vec{\mathcal{H}})\nabla_\G \varphi_h(\vec p,\cdot)\\
\notag
&=(\vec I-d\vec{\mathcal{H}})\Proj\Proj_h \vec t^a_h(\vec p,\cdot)\cdot\nabla_\G \varphi_h(\vec p,\cdot)\\
\notag
&={\vec{\mathcal Q}_h}\vec t^a_h(\vec p,\cdot)\cdot\nabla_\G \varphi_h(\vec p,\cdot)
\end{align}
where the last equality defines $\vec{\mathcal Q}_h$. We denote the lifted version by  $\vec{\mathcal Q}_h^l$  Thus we may write (\ref{pf_bbil_1}) as 
\begin{align}
\lv \bbil{\psi_h}{\varphi_h}{\vec t^a_h}-\bhbil{\Psi_h}{\Phi_h}{\vec T^a_h}\rv
&=
\lv \int_\Gt\psi_h\vec t^a_h\cdot\nabla_\Gt\varphi_h-\int_\Gt\frac{1}{\delta_h^l}\psi_h\vec{\mathcal Q}^l_h\vec t^a_h\cdot\nabla_\Gt\varphi_h\rv\\
\notag
&\leq
\lv \int_\Gt\psi_h(\vec I-\vec{\mathcal Q}^l_h)\Proj\vec t^a_h\cdot\nabla_\Gt\varphi_h\rv+ch^2,
\end{align}
where we have used (\ref{surface_element}). 
We now apply the following result from \cite[Lem 5.1]{dziuk2007finite} 
\beq
\sup_{t\in(0,T)}\sup_{\Gct}\lv\left( \vec I-\vec {\mathcal Q}_h\right)\Proj\rv\leq ch^2,
\eeq
which yields the desired bound.
\end{Proof}

\section{Ritz projection estimates}
It proves helpful in the analysis  to introduce the Ritz projection $\Ritz:\Hil{1}(\G)\to \Scl$ defined as follows: for $z\in\Hil{1}(\G)$ with $\int_\G z=0$,
\begin{equation}\label{eqn:RP_definition}
\abil{\Ritz z}{\varphi_h}=\abil{z}{\varphi_h}\myall\varphi_h\in\Scl,
\end{equation}
with $\int_\G \Ritz z=0$.
\begin{Lem}[Ritz projection estimates]\label{Lem:RP_bounds}
\changes{
\margnote{ref  2. pt. 56.}
We recall the following estimates  proved in \citep[Thm. 6.1 and Thm. 6.2]{dziuk2010l2} that hold for the mesh-size $h$ sufficiently small 
}
\begin{align}
\label{eqn:Ritz_bound}
\ltwon{z-\Ritz z}{\G}&+h\ltwon{\nabla_\G\left(z-\Ritz z\right)}{\G}\leq ch^2\Hiln{z}{2}{\G}.\\
\label{eqn:MD_Ritz_bound}
\ltwon{\mdth{\vec v^a_h}\left(z-\Ritz z\right)}{\G}&+h\ltwon{\nabla_\G\mdth{\vec v^a_h}\left(z-\Ritz z\right)}{\G}\leq\\
&ch^2\left(\Hiln{z}{2}{\G}+\Hiln{\mdt{\vec v_a}z}{2}{\G}\right).\notag
\end{align}
\end{Lem}
 \section*{Acknowledgments}
This research has been supported by the UK Engineering and Physical Sciences Research Council (EPSRC), Grant EP/G010404. The research of CV has been supported by the EPSRC grant EP/J016780/1. This research was finalised while CME and CV were participants in the Newton Institute Program: Free Boundary Problems and Related Topics. \changes{Both authors would like to express their thanks to the anonymous reviewers for their careful reading of the manuscript and helpful suggestions.}

  \end{document}